\documentclass[a4paper]{article}

\usepackage{amsfonts,amssymb,amsthm}

\usepackage[totalwidth=17cm,totalheight=24cm]{geometry}

\makeatletter
\newif\iflabel\labelfalse \let\w@label=\label
\def\label{\global\labeltrue\w@label}
\let\w@eqn=\equation
\def\equation{\global\labelfalse\w@eqn}
\def\endequation{\iflabel\@eeq\else\@eeqw\fi$$\global\@ignoretrue}
\def\@eeq{\eqno \@eqnnum}
\def\@eeqw{\addtocounter{equation}{-1}}
\newskip\bw@\newskip\bws@
\def\dc@@pq{\global \bw@=\belowdisplayskip \global\bws@=\belowdisplayshortskip
   \global\belowdisplayshortskip=3pt plus -3pt
   \global\belowdisplayskip=\belowdisplayshortskip
   \hskip 500pt minus 500pt\relax}
\def\dc@pq{\dc@@pq$$}
\def\fc@pq{$$\hskip 500pt minus 500pt\global
   \belowdisplayshortskip=\bws@\global\belowdisplayskip=\bw@}
\def\coupeq{\dc@pq\fc@pq}
\def\coupepas{\par\noindent\begin{minipage}{\textwidth}}
\def\jusquela{\end{minipage}}
\makeatother

\catcode`\@=11\relax
\newwrite\@unused
\def\typeout#1{{\let\protect\string\immediate\write\@unused{#1}}}
\def\@nnil{\@nil}
\def\@empty{}
\def\@psdonoop#1\@@#2#3{}
\def\@psdo#1:=#2\do#3{\edef\@psdotmp{#2}\ifx\@psdotmp\@empty \else
    \expandafter\@psdoloop#2,\@nil,\@nil\@@#1{#3}\fi}
\def\@psdoloop#1,#2,#3\@@#4#5{\def#4{#1}\ifx #4\@nnil \else
       #5\def#4{#2}\ifx #4\@nnil \else#5\@ipsdoloop #3\@@#4{#5}\fi\fi}
\def\@ipsdoloop#1,#2\@@#3#4{\def#3{#1}\ifx #3\@nnil 
       \let\@nextwhile=\@psdonoop \else
      #4\relax\let\@nextwhile=\@ipsdoloop\fi\@nextwhile#2\@@#3{#4}}
\def\@tpsdo#1:=#2\do#3{\xdef\@psdotmp{#2}\ifx\@psdotmp\@empty \else
    \@tpsdoloop#2\@nil\@nil\@@#1{#3}\fi}
\def\@tpsdoloop#1#2\@@#3#4{\def#3{#1}\ifx #3\@nnil 
       \let\@nextwhile=\@psdonoop \else
      #4\relax\let\@nextwhile=\@tpsdoloop\fi\@nextwhile#2\@@#3{#4}}
\def\psdraft{
        \def\@psdraft{0}
}
\def\psfull{
        \def\@psdraft{100}
}
\psfull
\newif\if@prologfile
\newif\if@postlogfile
\newif\if@noisy
\def\pssilent{
        \@noisyfalse
}
\def\psnoisy{
        \@noisytrue
}
\psnoisy
\newif\if@bbllx
\newif\if@bblly
\newif\if@bburx
\newif\if@bbury
\newif\if@height
\newif\if@width
\newif\if@rheight
\newif\if@rwidth
\newif\if@clip
\newif\if@verbose
\def\@p@@sclip#1{\@cliptrue}
\def\@p@@sfile#1{
                   \def\@p@sfile{#1}
}
\def\@p@@sfigure#1{\def\@p@sfile{#1}}
\def\@p@@sbbllx#1{
                \@bbllxtrue
                \dimen100=#1
                \edef\@p@sbbllx{\number\dimen100}
}
\def\@p@@sbblly#1{
                \@bbllytrue
                \dimen100=#1
                \edef\@p@sbblly{\number\dimen100}
}
\def\@p@@sbburx#1{
                \@bburxtrue
                \dimen100=#1
                \edef\@p@sbburx{\number\dimen100}
}
\def\@p@@sbbury#1{
                \@bburytrue
                \dimen100=#1
                \edef\@p@sbbury{\number\dimen100}
}
\def\@p@@sheight#1{
                \@heighttrue
                \dimen100=#1
                \edef\@p@sheight{\number\dimen100}
}
\def\@p@@swidth#1{
                \@widthtrue
                \dimen100=#1
                \edef\@p@swidth{\number\dimen100}
}
\def\@p@@srheight#1{
                \@rheighttrue
                \dimen100=#1
                \edef\@p@srheight{\number\dimen100}
}
\def\@p@@srwidth#1{
                \@rwidthtrue
                \dimen100=#1
                \edef\@p@srwidth{\number\dimen100}
}
\def\@p@@ssilent#1{ 
                \@verbosefalse
}
\def\@p@@sprolog#1{\@prologfiletrue\def\@prologfileval{#1}}
\def\@p@@spostlog#1{\@postlogfiletrue\def\@postlogfileval{#1}}
\def\@cs@name#1{\csname #1\endcsname}
\def\@setparms#1=#2,{\@cs@name{@p@@s#1}{#2}}
\def\ps@init@parms{
                \@bbllxfalse \@bbllyfalse
                \@bburxfalse \@bburyfalse
                \@heightfalse \@widthfalse
                \@rheightfalse \@rwidthfalse
                \def\@p@sbbllx{}\def\@p@sbblly{}
                \def\@p@sbburx{}\def\@p@sbbury{}
                \def\@p@sheight{}\def\@p@swidth{}
                \def\@p@srheight{}\def\@p@srwidth{}
                \def\@p@sfile{}
                \def\@p@scost{10}
                \def\@sc{}
                \@prologfilefalse
                \@postlogfilefalse
                \@clipfalse
                \if@noisy
                        \@verbosetrue
                \else
                        \@verbosefalse
                \fi
}
\def\parse@ps@parms#1{
                \@psdo\@psfiga:=#1\do
                   {\expandafter\@setparms\@psfiga,}}
\newif\ifno@bb
\newif\ifnot@eof
\newread\ps@stream
\def\bb@missing{
        \if@verbose{
                \typeout{psfig: searching \@p@sfile \space  for bounding box}
        }\fi
        \openin\ps@stream=\@p@sfile
        \no@bbtrue
        \not@eoftrue
        \catcode`\%=12
        \loop
                \read\ps@stream to \line@in
                \global\toks200=\expandafter{\line@in}
                \ifeof\ps@stream \not@eoffalse \fi
                \@bbtest{\toks200}
                \if@bbmatch\not@eoffalse\expandafter\bb@cull\the\toks200\fi
        \ifnot@eof \repeat
        \catcode`\%=14
}       
\catcode`\%=12
\newif\if@bbmatch
\def\@bbtest#1{\expandafter\@a@\the#1
\long\def\@a@#1
\long\def\bb@cull#1 #2 #3 #4 #5 {
        \dimen100=#2 bp\edef\@p@sbbllx{\number\dimen100}
        \dimen100=#3 bp\edef\@p@sbblly{\number\dimen100}
        \dimen100=#4 bp\edef\@p@sbburx{\number\dimen100}
        \dimen100=#5 bp\edef\@p@sbbury{\number\dimen100}
        \no@bbfalse
}
\catcode`\%=14
\def\compute@bb{
                \no@bbfalse
                \if@bbllx \else \no@bbtrue \fi
                \if@bblly \else \no@bbtrue \fi
                \if@bburx \else \no@bbtrue \fi
                \if@bbury \else \no@bbtrue \fi
                \ifno@bb \bb@missing \fi
                \ifno@bb \typeout{FATAL ERROR: no bb supplied or found}
                        \no-bb-error
                \fi
                \count203=\@p@sbburx
                \count204=\@p@sbbury
                \advance\count203 by -\@p@sbbllx
                \advance\count204 by -\@p@sbblly
                \edef\@bbw{\number\count203}
                \edef\@bbh{\number\count204}
}
\def\in@hundreds#1#2#3{\count240=#2 \count241=#3
                     \count100=\count240        
                     \divide\count100 by \count241
                     \count101=\count100
                     \multiply\count101 by \count241
                     \advance\count240 by -\count101
                     \multiply\count240 by 10
                     \count101=\count240        
                     \divide\count101 by \count241
                     \count102=\count101
                     \multiply\count102 by \count241
                     \advance\count240 by -\count102
                     \multiply\count240 by 10
                     \count102=\count240        
                     \divide\count102 by \count241
                     \count200=#1\count205=0
                     \count201=\count200
                        \multiply\count201 by \count100
                        \advance\count205 by \count201
                     \count201=\count200
                        \divide\count201 by 10
                        \multiply\count201 by \count101
                        \advance\count205 by \count201
                     \count201=\count200
                        \divide\count201 by 100
                        \multiply\count201 by \count102
                        \advance\count205 by \count201
                     \edef\@result{\number\count205}
}
\def\compute@wfromh{
                \in@hundreds{\@p@sheight}{\@bbw}{\@bbh}
                \edef\@p@swidth{\@result}
}
\def\compute@hfromw{
                \in@hundreds{\@p@swidth}{\@bbh}{\@bbw}
                \edef\@p@sheight{\@result}
}
\def\compute@handw{
                \if@height 
                        \if@width
                        \else
                                \compute@wfromh
                        \fi
                \else 
                        \if@width
                                \compute@hfromw
                        \else
                                \edef\@p@sheight{\@bbh}
                                \edef\@p@swidth{\@bbw}
                        \fi
                \fi
}
\def\compute@resv{
                \if@rheight \else \edef\@p@srheight{\@p@sheight} \fi
                \if@rwidth \else \edef\@p@srwidth{\@p@swidth} \fi
}
\def\compute@sizes{
        \compute@bb
        \compute@handw
        \compute@resv
}
\def\psfig#1{\vbox {
        \ps@init@parms
        \parse@ps@parms{#1}
        \compute@sizes
        \ifnum\@p@scost<\@psdraft{
                \if@verbose{
                        \typeout{psfig: including \@p@sfile \space }
                }\fi
                \special{ps::[begin]    \@p@swidth \space \@p@sheight \space
                                \@p@sbbllx \space \@p@sbblly \space
                                \@p@sbburx \space \@p@sbbury \space
                                startTexFig \space }
                \if@clip{
                        \if@verbose{
                                \typeout{(clip)}
                        }\fi
                        \special{ps:: doclip \space }
                }\fi
                \if@prologfile
                    \special{ps: plotfile \@prologfileval \space } \fi
                \special{ps: plotfile \@p@sfile \space }
                \if@postlogfile
                    \special{ps: plotfile \@postlogfileval \space } \fi
                \special{ps::[end] endTexFig \space }
                \vbox to \@p@srheight true sp{
                        \hbox to \@p@srwidth true sp{
                                \hss
                        }
                \vss
                }
        }\else{
                \vbox to \@p@srheight true sp{
                \vss
                        \hbox to \@p@srwidth true sp{
                                \hss
                                \if@verbose{
                                        \@p@sfile
                                }\fi
                                \hss
                        }
                \vss
                }
        }\fi
}}
\catcode`\@=12\relax

\newcommand{\vs}{\vspace{0.3cm}}
\newcommand{\dr}{\partial}
\newcommand{\C}{{\mathbb C}}
\newcommand{\N}{{\mathbb N}}
\newcommand{\R}{{\mathbb R}}
\newcommand{\Z}{{\mathbb Z}}
\newcommand{\II}{I\hspace{-0.1cm}I}
\newcommand{\III}{I\hspace{-0.1cm}I\hspace{-0.1cm}I}
\newcommand{\tr}{\mbox{\rm tr}}
\newcommand{\ric}{\mbox{\rm ric}}
\newcommand{\cotg}{\mbox{\rm cotg}}
\newcommand{\SO}{\mbox{\rm SO}}
\newcommand{\ricb}{\overline{\ric}}
\newcommand{\Rc}{\check{R}}
\newcommand{\Rr}{\stackrel{\circ}{R}}
\newcommand{\diam}{\mbox{\rm diam}}
\newcommand{\area}{\mbox{\rm area}}

\def\pointir{\unskip  {. --- \ignorespaces }\hskip0cm}

\newcommand{\deltab}{\overline{\delta}}

\newcommand{\gba}{\overline{g}}
\newcommand{\fb}{\overline{f}}
\newcommand{\hb}{\overline{h}}
\newcommand{\kb}{\overline{k}}
\newcommand{\Db}{\overline{D}}
\newcommand{\Rb}{\overline{R}}
\newcommand{\Kb}{\overline{K}}
\newcommand{\Sb}{\overline{S}}
\newcommand{\Omegab}{\overline{\Omega}}

\newcommand{\pt}{\tilde{p}}
\newcommand{\ut}{\tilde{u}}
\newcommand{\yt}{\tilde{y}}
\newcommand{\Bt}{\tilde{B}}
\newcommand{\Ft}{\tilde{F}}
\newcommand{\Kt}{\tilde{K}}
\newcommand{\ept}{\tilde{\epsilon}}
\newcommand{\gat}{\tilde{\gamma}}
\newcommand{\kappat}{\tilde{\kappa}}
\newcommand{\thetat}{\tilde{\theta}}

\newcommand{\rhot}{\tilde{\rho}}
\newcommand{\taut}{\tilde{\tau}}
\newcommand{\phit}{\tilde{\phi}}
\newcommand{\sit}{\tilde{\sigma}}
\newcommand{\Phit}{\tilde{\Phi}}
\newcommand{\Sigmat}{\tilde{\Sigma}}
\newcommand{\Omt}{\tilde{\Omega}}

\newcommand{\kad}{\stackrel{\bullet}{\kappa}}
\newcommand{\gd}{\stackrel{\bullet}{g}}

\newtheorem{prop}{Proposition}[section]
\newtheorem{df}[prop]{Definition}
\newtheorem{lemma}[prop]{Lemma}
\newtheorem{thm}[prop]{Theorem}
\newtheorem{cor}[prop]{Corollary}
\newtheorem{asser}[prop]{Assertion}
\newtheorem{remark}[prop]{Remark}

\newenvironment{thn}[1]{\vskip 0.2cm \noindent{\bf Theorem #1.} \it}{\rm
\vspace{0.2cm}} 
\newenvironment{lmn}[1]{\vskip 0.2cm \noindent{\bf Lemma #1.} \it}{\rm
\vspace{0.2cm}} 

\newcommand{\btm}{\begin{thm}}
\newcommand{\etm}{\end{thm}}
\newcommand{\blm}{\begin{lemma}}
\newcommand{\elm}{\end{lemma}}
\newcommand{\bcr}{\begin{cor}}
\newcommand{\ecr}{\end{cor}}
\newcommand{\bdf}{\begin{df}}
\newcommand{\edf}{\end{df}}
\newcommand{\bprop}{\begin{prop}}
\newcommand{\eprop}{\end{prop}}
\newcommand{\bas}{\begin{asser}}
\newcommand{\eas}{\end{asser}}
\newcommand{\beq}{\begin{equation}}
\newcommand{\eeq}{\end{equation}}
\newcommand{\bpv}{\begin{proof}}
\newcommand{\epv}{\end{proof}}
\newcommand{\bit}{\begin{itemize}}
\newcommand{\eit}{\end{itemize}}
\newcommand{\bpn}{\begin{pfn}}
\newcommand{\epn}{\end{pfn}}
\newcommand{\btn}{\begin{thn}}
\newcommand{\etn}{\end{thn}}
\newcommand{\bln}{\begin{lmn}}
\newcommand{\eln}{\end{lmn}}

\newenvironment{pfn}[1]{\vskip 0.2cm \noindent{\it Proof #1.}}{$\square$
\vspace{0.2cm}} 

\newcommand{\Omb}{\overline{\Omega}}
\newcommand{\Sib}{\overline{\Sigma}}
\newcommand{\gb}{\overline{g}}
\newcommand{\Ub}{\overline{U}}
\newcommand{\Wb}{\overline{W}}
\newcommand{\db}{\overline{\partial}}

\newcommand{\Met}{\mathcal{M}et}
\newcommand{\Imm}{\mathcal{I}mm}
\newcommand{\CMet}{\mathcal{CM}et}
\newcommand{\cA}{\mathcal{A}}
\newcommand{\cC}{\mathcal{C}}
\newcommand{\cD}{\mathcal{D}}
\newcommand{\cE}{\mathcal{E}}
\newcommand{\cM}{\mathcal{M}}
\newcommand{\cS}{\mathcal{S}}
\newcommand{\cP}{\mathcal{P}}
\newcommand{\cV}{\mathcal{V}}
\newcommand{\cW}{\mathcal{W}}
\newcommand{\CImm}{\mathcal{CI}mm}
\newcommand{\gab}{\overline{\gamma}}
\newcommand{\hyp}{\mathbf{H}^3}
\newcommand{\dhyp}{\partial\hyp}
\newcommand{\cL}{\mathcal{L}}
\newcommand{\isom}{\mathrm{Isom}}
\newcommand{\im}{\mathrm{Im}}
\newcommand{\Na}{\nabla}
\newcommand{\Nat}{\tilde{\nabla}}
\newcommand{\Sit}{\tilde{\Sigma}}

\newcommand{\eps}{\epsilon}
\newcommand{\ga}{\gamma}
\newcommand{\si}{\sigma}
\newcommand{\om}{\omega}
\newcommand{\Ga}{\Gamma}
\newcommand{\La}{\Lambda}
\newcommand{\Si}{\Sigma}
\newcommand{\Om}{\Omega}

\begin{document}

\title{Circle patterns on singular surfaces}

\date{February 2006 (v2)}

\author{Jean-Marc Schlenker\thanks{
Laboratoire Emile Picard, UMR CNRS 5580,
Institut de Math{\'e}matiques, 
Universit{\'e} Paul Sabatier,
31062 Toulouse Cedex 9,
France.
\texttt{schlenker@math.ups-tlse.fr; http://www.picard.ups-tlse.fr/\~{
}schlenker}. }
\thanks{The author would like to thank the RIP program at Oberwolfach, where
  part of the research presented here was conducted. Partially supported by
  the ``ACI Jeunes Chercheurs'' {\it M\'etriques privil\'egi\'es sur les
    vari\'et\'es \`a bord}, 2003-06.}}

\maketitle

\begin{abstract}

We consider ``hyperideal'' circle patterns, i.e. patterns of disks appearing
in the definition of the Delaunay decomposition associated to a set of
disjoint disks, possibly with cone singularities at the center of those disks.
Hyperideal circle patterns are associated to hyperideal hyperbolic
polyhedra. We describe the possible intersection angles and singular
curvatures of those circle patterns, on Euclidean or
hyperbolic surfaces with conical singularities. 
This is related to results on the dihedral angles of
ideal or hyperideal hyperbolic polyhedra. The results presented here extend
those in \cite{hcp}, however the proof is completely different (and more
intricate) since \cite{hcp} used a shortcut which is not available here.

\bigskip

\begin{center} {\bf R{\'e}sum{\'e}} \end{center}

On consid\`ere des motifs de cercles ``hyperid\'eaux'', i.e. ceux qui
apparaissent dans la d\'efinition de d\'ecomposition Delaunay associ\'e \`a
un ensemble de disques disjoint, \'eventuellement avec des singularit\'es
coniques aux centre de ces disques. Les motifs de cercles hyperid\'eaux
sont associ\'es aux poly\`edres hyperboliques hyperid\'eaux. 
On d\'ecrit les angles d'intersections et les courbures singuli\`eres
possibles de ces motifs de cercles, sur les surfaces
euclidiennes ou hyperboliques \`a singularit\'es coniques. 
C'est li\'e \`a des r\'esultats sur les angles di\`edres des
poly\`edres hyperboliques id\'eaux ou hyperid\'eaux. Les r\'esultats
pr\'esent\'es ici \'etendent ceux de \cite{hcp}, mais les preuves sont
compl\`etement diff\'erentes (et plus \'elabor\'es) car \cite{hcp} prenait un
racour\c{c}i qui n'est pas utilisable ici.

\end{abstract}

\maketitle

\vspace{0.4cm}


\vspace{0.4cm}


\bigskip

\tableofcontents

\section{Introduction and results}

\subsection{Motivations.}

\paragraph{Delaunay decompositions.}

Let $x_1, \cdots, x_n\in \R^2$ be points, and suppose that the convex hull 
$C$ of those points has non-empty interior. There is then a unique 
{\it Delaunay decomposition} of $C$: a decomposition of $C$ as the
union of convex polygons 
with vertices at the $x_i$ such that, for each face $f$, there is a circle
$C_f$ which contains exactly the $x_i$ which are vertices of $f$ but whose
interior contains none of the $x_i$. 

The same construction can be made in the sphere (with at least 3 points) or in 
a hyperbolic surface. The circles appearing in this definition, with their
peculiar types of intersections, are forming what is defined below as an {\it
ideal 
circle pattern}. Ideal circle patterns are related in a natural way to 
ideal hyperbolic polyhedra.

The possible intersection angles of ideal circle patterns can be described
in terms of some linear equalities and inequalities. This was done for the
sphere (in the context of ideal hyperbolic polyhedra) by Andreev
\cite{andreev-ideal} and Rivin \cite{rivin-annals}. Related results were
obtained on the sphere by Luo
\cite{luo-sphere}, and on higher genus surfaces by different authors
\cite{bowditch,garrett,leibon1,Ri2,rivin-combi};  
in particular Bobenko and Springborn \cite{bobenko-springborn} recently gave
a statement, for hyperbolic surfaces, which is closely related to Theorem
\ref{tm:hyperb} below. A similar result is also stated in \cite{ideal}.

\paragraph{Delaunay decompositions with singular points.}

Suppose now that instead of being simply points on an Euclidean (or
hyperbolic) surface $S$, $x_1,
\cdots, x_n$ are conical singularities. It remains true that there exists a
unique Delaunay decomposition of $S$ with vertices at the $x_i$. The circles
which appear in the definition of this Delaunay decomposition still constitute
an ``ideal circle pattern'', but now there are conical singularities at the
intersection points of the circles. Again it is of interest to understand the
possible intersection angles of the circles, and results on this question can
be found in \cite{bowditch,garrett,leibon1,bobenko-springborn}.

\paragraph{Delaunay decomposition associated to disks.}

Another generalization of the notion of Delaunay decomposition can be obtained
by replacing the points $x_1, \cdots, x_n$ by disjoint disks $D_1, \cdots,
D_n$. To such a finite set of disks, say in an Euclidean surface $S$, one can
still associated in a unique way an embedded graph $\Gamma$, with vertices
$v_1, \cdots, v_n$ corresponding to the $D_i$, such that to each face $f$ of
$\Gamma$ is associated a circle $C_f$ which is orthogonal to the $D_i$
corresponding to the vertices of $f$,  but either does not intersect, or
intersects with an angle less than $\pi/2$, all the other disks $D_i$.

The circles $C_f$ involved in this definition constitute what is defined below
as a ``hyperideal circle pattern''. Hyperideal circle patterns are naturally
associated to hyperideal hyperbolic polyhedra. 
The $D_i$ are actually allowed in the
definition to be points, so that ideal circle patterns are special cases
of hyperideal circle patterns.

As for ideal circle patterns, it is interesting to describe the possible
intersection angles between the circles in a hyperideal circle pattern. 
For the sphere, a complete description is 
hidden behind a theorem of Bao and Bonahon \cite{bao-bonahon} on
hyperideal polyhedra, while statements concerning surfaces of genus at least
$2$ are similarly consequences of results \cite{rousset1} on equivariant
hyperideal polyhedra. Further results can be found in \cite{hphm,hcp}.

It is also possible to consider the Delaunay decomposition associated to a set
of disks containing conical singularities at their centers. The main goal of
this paper is to give simple statements describing the possible intersection
angles of the circles in this context, when the underlying metric is either
Euclidean (this is Theorem \ref{tm:euclid}) or hyperbolic (Theorem
\ref{tm:hyperb}). Those statements are new even for ideal circle patterns with
conical singularities. 

\paragraph{Circle patterns as descriptions of singular surfaces.}

Consider an Euclidean (or possibly hyperbolic) surface with conical
singularities, for instance a polyhedral surface in Euclidean space with its
induced metric. Such surfaces occur in many practical applications and it is
interesting to find practical ways to describe them. The most natural way is
of course to choose a triangulation, with the conical points as vertices, and
to describe the surface by the combinatorics of the triangulation and the
lengths of the edges. 

There is another possible way,  however, based on the Delaunay decomposition
associated to the conical singularities. It provides an ideal circle pattern
which is uniquely determined by the singular curvature at the cone
points of the metric, as well as some
intersection angles between the circles. It follows from Theorem
\ref{tm:euclid} below that the singular metric can be uniquely reconstructed
from this data. This idea is used in practical manner in \cite{KSS}. 
This description of singular metrics has some interesting characteristics:
\begin{itemize}
\item it gives a direct access to the curvature, since the total curvature of
  a domain is simply the sum of the curvatures of the singular points it
  contains, which is part of the description.
\item it is sensitive to the conformal structure, 
in particular it leads directly to a ``discrete conformal map'' to a domain in
the plane, obtrained by constructing the circle 
patterns with the same intersection
angles but with no singular curvature. 
\item the complete description of the possible data (for a given
  combinatorics) is ``simple'', since it is given by a set of linear
  inequalities between the possible intersection angles and singular
  curvatures. 
\item the computations which are needed involve the maximization of a 
  functional
  under linear constraints, so that they should be algorithmically simple.
\end{itemize}

It it also possible to apply the same idea, taking the Delaunay decomposition
associated to disjoints disks centered at the singular points. Such a set of
disks uniquely determines a hyperideal circle pattern, and, again by Theorem
\ref{tm:euclid}, the surface and the radii of the disks can be recovered from
the combinatorics, the singular curvatures and the intersection angles of the
circles. Compared to a description by an ideal circle pattern, this
description has the advantage of being more flexible, since the possible
intersection angles between the circles are required to satisfy affine
inequalities rather than equalities.

\paragraph{Relations with geometric topology.} 

Circle packings, and ``ideal'' circle patterns, are related to 3-dimensional
hyperbolic manifolds, orbifolds and cone-manifolds. Given an ideal circle
pattern, one can associate to it a non-complete 3-dimensional hyperbolic
manifold with polyhedral boundary. Gluing two copies of this manifolds to along
the boundary yields a 3-dimensional hyperbolic cone-manifold of finite volume,
which is an orbifold when the intersection angles of the original circle
pattern are of the form $\pi/k$, $k\geq 2$. In the same manner, ``hyperideal''
circle patterns natural yield, by a similar construction, compact 3-dimensional
cone-manifolds (resp. orbifolds). This point
will not be developed much here but it appears in section 7 as a key tool in
the proof of the main results.

\paragraph{A short comment on the proof.}

The proof is based on the deformation method: for a given combinatorics, one
considers the natural map sending a ``hyperideal'' circle pattern, with the
right incidence graph, to the set of its intersection angles and singular
curvatures. One has to prove that this map is a homeomorphism. The key point
in the proof is the fact that the map has injective differential at each
point, in other terms that hyperideal circle patterns are ``rigid''. 

There are several different methods to prove this kind of rigidity results. In
\cite{bao-bonahon}, it was done using an argument going back to Legendre
\cite{legendre} and Cauchy \cite{Cauchy} (see Sabitov's illuminating paper on
the subject \cite{sabitov-legendre}). In some situations involving equivariant
polyhedra (see e.g. \cite{rousset1}), the infinitesimal rigidity can be
obtained, through a transformation due to Pogorelov \cite{Po}, as consequences
of similar statements for equivariant polyhedra in the Minkowski space, see
\cite{polygones} (similar ideas were used for smooth surfaces in \cite{iie}). 
The results in \cite{hcp} are based on results (from
\cite{otal,bonahon-otal}) on the geometry of the convex core of hyperbolic
3-manifolds, which in turn rely on a local rigidity result of Hogdson and
Kerckhoff \cite{HK} for this setting. 

The proof of the infinitesimal rigidity used here, by contrast, is based on a
volume argument going back to \cite{bragger,Ri2} for ideal hyperbolic
polyhedra. The key point, already used in \cite{hphm,rcnp} is that those
argument also work for hyperideal polyhedra, because their volume is, as for
ideal polyhedra, a strictly concave function of their dihedral angles. 

\paragraph{Another possible approach.}

Some of the results presented here, in particular those concerning hyperideal
circle patterns on Euclidean surfaces, can also be obtained in a different
way, replacing the deformation approach favored here by a direct study of the
critical points of some volume-based functionals. This is much more explicit,
or even constructive, than the approach followed here, and it is better suited
to a computer implementation. On the other hand, the proof given here is
perhaps more flexible.

The second method of proof, based more directly on singular points of
volume-based functionals, was followed recently by
Boris Springborn \cite{springborn-hcp}, leading to different 
but strongly related
results on Euclidean circle patterns. The results presented here are
actually to a large extend the result of a collaboration with Boris
Springborn, whose contribution is important. 
We decided however to write two papers, each one specializing in
the approach towards which he felt more inclined.

\subsection{Definitions}

\paragraph{Ideal circle packings.}

We now move to a more precise description of the results. 
We define ideal circle patterns in $S^2$, although the definition
could be given in the Euclidean or the hyperbolic plane. Ideal circle 
patterns are naturally associated to ideal hyperbolic polyhedra, a point which
should become clear below.

\begin{df} \label{df:ideal}
A {\bf circle pattern} on $S^2$ is a finite family of oriented circles
$C_1, \cdots, C_N$. Given a circle pattern, an {\bf interstice} is a connected
component of the complement of the union of the open disks bounded by the
circles. If $C_1, \cdots, C_N$ is a circle pattern, it is {\bf ideal} if:
\begin{itemize}
\item each interstice is a point,
\item each circle contains at least 3 interstices,
\item if $D$ is an open disk in $S^2$, containing no interstice, but such that
  its closure contains at least 3 of the interstices, then $D$ is the open
  disk bounded by one of the $C_i, 1\leq i\leq N$.
\end{itemize}
Given an ideal circle pattern, its {\it incidence graph} is the graph,
embedded in the sphere which has:
\begin{itemize}
\item one vertex for each circle,
\item one edge between two vertices, when the corresponding circles intersect
  at two interstices. 
\end{itemize}
\end{df}

Clearly, each ideal circle pattern is the pattern of circles appearing in the
definition of the Delaunay decomposition of a set of points, namely its
interstices. Notice that, with the definition given here, two circles can 
intersect while the corresponding vertices of the incidence graph are not
adjacent. 

Again, the definition can clearly be extended from the sphere to any
hyperbolic or Euclidean surface.
Results concerning the possible intersection angles of ideal circle
patterns can be found in \cite{bobenko-springborn,rivin-combi,leibon1,ideal}.

\paragraph{Hyperideal circle patterns.}

There is a related notion of ``hyperideal'' circle patterns, which are
related to hyperideal hyperbolic polyhedra as ideal circle patterns are
related to ideal hyperbolic polyhedra. 
Given a circle pattern, the {\it intersection angle} between two 
intersecting circles is $\pi$ minus the angle measured in the intersection
of the disks bounded by the two circles. In other terms, this angle
is measured in the complement of one of the disks in the other. This
angle is defined to be equal to $\pi$ when the circles are tangent. 

\begin{df} \label{df:hyper}
Let $C_1, \cdots, C_N$ be a circle pattern in $S^2$, with interstices
$I_1,\cdots, I_M$. It is {\bf hyperideal} if:
\begin{itemize}
\item each interstice either is a point, or is topologically a disk,
\item for each $j\in \{ 1, \cdots, M\}$, 
  corresponding to an interstice which is 
  not a point, there is an oriented circle $C'_j$,
  containing $I_j$, which is orthogonal to all the circles $C_i$ adjacent to
  $I_j$. If $I_j$ is a point we set $C'_j:=I_j$,
\item for all $i\in \{1,\cdots, N\}$ and all $j\in \{ 1,\cdots, M\}$, if $C_i$
  is not adjacent to $I_j$, then either the interior of 
  $C_i$ is disjoint from
  $C'_j$, or $C_i$ intersects $C'_j$ (which is not a point) and  
  their intersection angle is strictly
  larger than $\pi/2$, 
\item if $D$ is an open disk in $S^2$ such that:
  \begin{enumerate}
  \item for each $j\in \{ 1,\cdots, M\}$, either $D$ is disjoint from
  $C'_j$, or $\dr D$ has an intersection angle at least $\pi/2$
  with $C'_j$,
  \item there are at least three $C'_j$ which are either points contained in
  $\dr D$ or circles orthogonal to $\dr D$,
  \end{enumerate}
then $\dr D$ is one of the $C_i$.
\end{itemize}
Such a pattern is {\it strictly hyperideal} if no interstice is
reduced to a point. The circles $C_i$ are called {\it principal} circles,
while the circles (or points) $C'_j$ are the {\it dual} circles.
\end{df}

\begin{figure}[ht]
\centerline{\psfig{figure=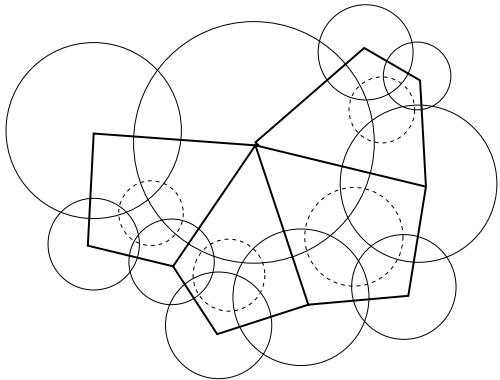,height=5cm}}
\caption{A hyperideal circle pattern and its incidence graph (the dual circles
  are dashed).}
\end{figure}

The incidence graph of a hyperideal circle pattern is defined as for
ideal circle patterns in Definition \ref{df:ideal} above. Hyperideal circle
patterns have a simple relation to circle packings: those correspond simply to
the limit when all intersection angles go to $\pi$.

Clearly the notion of circle pattern -- and the notion of angles between the
circles -- is not limited to spherical metrics; one could also consider
Euclidean or hyperbolic metrics. Actually one could simply consider a complex
projective structure, also called a $\C P^1$-structure. This viewpoint is
often interesting when thinking about circle patterns, however it will not
appear much here.

Some results concerning the possible intersection angles of hyperideal
circle patterns can be found in \cite{hcp,hphm} and (in the related
setting of hyperideal hyperbolic polyhedra) in \cite{bao-bonahon,rousset1}.
Our main goal here is to extend them to singular surfaces.

Given a hyperideal circle pattern, there are actually two families of circles
that one can consider. The first is made of the principal circles, 
appearing directly in
the circle pattern; the other contains the dual circles, 
associated to the faces of
the incidence graph, which are orthogonal to the circles of the first
family. Those dual circles are disjoint, but they also provide a
circle pattern of another kind, which is also of some interest, although 
this point will not be pursued here.

\paragraph{Circle patterns on singular surfaces}

We have already mentioned that the object of this paper is to extend results
on hyperideal circle patterns to Euclidean or hyperbolic surfaces with conical
singularities. But let $\Sigma_0$ be a complete Euclidian or hyperbolic
surface with one conical singularity. It is then not difficult to check that 
the interior of a
circle on $\Sigma_0$ can not contain the singular point, except at its center
-- otherwise the circle would not ``close up''. 

It follows that hyperideal circle patterns can only be considered on surfaces
with conical singularities if the singularities are at the centers of
the principal or the dual circles. We consider
here only the second possibilities, i.e. singularities at the centers of the
dual circles. More general results could perhaps be obtained by considering
also singularities at the centers of the principal circles,
however it would be at the cost of more complicated notations and of some
added technical difficulties. It is interesting to note, however, that some of
the technical points discussed in the proofs of the main results extend to
this more general situation.

\subsection{Main results.}

\paragraph{Admissible domains.}

We need one more definition before giving the first result of this text. Let
$\Gamma$ be a graph embedded in a closed surface $\Sigma$. 

\begin{df} \label{df:15}
An {\bf admissible domain} in $(\Sigma,\Gamma)$ is a connected open
domain $\Omega$, which contais a face of $\Gamma$, such that $\dr \Omega$ is a
finite union of segments which: 
\begin{itemize}
\item have as endpoints vertices of $\Gamma$,
\item either are edges of $\Gamma$ or are contained (except for their
  endpoints) in an open face of $\Gamma$.
\end{itemize}
To each such admissible domain, we can associate two numbers: its Euler 
characteristic, $\chi(\Omega)$,
and the number of boundary segments contained in open faces of $\Gamma$,
$m(\Omega)$. The {\bf boundary} of $\Omega$ will not be understood in the
usual way, but rather as the parametrized polygonal curve immersed in 
$\Sigma$ as the boundary of $\Omega$, considered as the image of an embedding
of a surface with boundary. 
For instance, edges of $\Gamma$ which are not contained in
$\Omega$ but are bounding two faces contained in $\Omega$ are considered as
contained {\bf twice} in $\dr\Omega$.
\end{df}

In other terms, the boundary of $\Omega$ is contained in $\Gamma$, except for
some segments of $\dr\Omega$ which are contained in faces of $\Gamma$.
The definition of the boundary of $\Omega$ should be compared with similar
considerations in \cite{bobenko-springborn}.

\paragraph{Main results on closed surfaces.}

The first result describes hyperideal circle patterns on a closed surface
with a singular Euclidean metric. Given a graph $\Gamma$, we 
denote by $\Gamma_1$ the set of its edges, and by $\Gamma_2$ the set of 
its faces.

\begin{thm} \label{tm:euclid}
Let $\Sigma$ be a closed orientable surface, and let 
$\Gamma$ be the 1-skeleton of a cellular decomposition of $\Sigma$. Let
$\kappa:\Gamma_2\rightarrow (-\infty, 2\pi)$ and
$\theta:\Gamma_1\rightarrow (0,\pi)$ be two functions. There exists a flat
metric $h$ with conical singularities on $\Sigma$, with a hyperideal circle
pattern $\sigma$ 
with incidence graph $\Gamma$, intersection angles given by $\theta$,
and singular points of curvatures given by $\kappa$ at the centers of the dual
circles, if and only if:
\begin{enumerate}
\item $\sum_{f\in \Gamma_2} \kappa(f)=2\pi\chi(\Sigma)$,
\item for any admissible domain $\Omega\subset \Sigma$ : 
$$ \sum_{e\in \Gamma_1, e\subset\dr\Omega} \theta(e) \geq (2\chi(\Omega) - 
m(\Omega))\pi - \sum_{f\in \Gamma_2, f\subset \Omega} \kappa(f)~, $$
with strict inequality except perhaps when $\Omega$ is a face of $\Gamma$.
\end{enumerate}
The metric $h$ is then unique up to homotheties, and $\sigma$ is unique given
$h$. 
\end{thm}

The sum over the edges of $\Gamma$ which are in $\dr\Omega$ should be
understood as mentioned above, $\dr\Omega$ is considered as a parametrized
immersed polygonal curve, and the edges of $\Gamma$ in $\dr\Omega$
which bound two faces in $\Omega$ are thus counted {\it twice}.

Note that, in Theorem \ref{tm:euclid}, the flat metric on $\Sigma$ is not
fixed: it is ``chosen'' by the combinatorics of $\Gamma$ and by the functions
$\theta$ and $\kappa$. The cone singularities mentioned here are at the 
centers of the dual circles, with a cone singularity of singular curvature 
$\kappa(f)$ at the center of the dual circle corresponding to the 
face $f$ of $\Gamma$. It is interesting to note that this result can be
formulated in a simpler manner under some more restricted conditions on the
conical singularities, in particular when there is no singularity at all,
or when $\kappa(f)\geq 0$ for all faces $f$ of $\Gamma$; then it is sufficient
to consider admissible domains $\Omega$ such that $\chi(\Omega)\geq 1$, 
i.e. disks, and that $m(\Omega)\leq 1$. The case with no singular points
is treated in \cite{hcp}, while \cite{bobenko-springborn} considers
singular surfaces, but in the restricted context of 
``ideal'' circle patterns, and contains a slightly different statement.

Even under the general conditions under which they are stated, the 
hypothesis are not as complicated as they appear at first sight. It is 
not too difficult, computationally speaking, to enumerate the possible
domains on which condition (2) has to be checked; for most of those 
domains, the right-hand side of the equation in condition (2)
is likely to be negative, so that the condition is trivially satisfied.

A similar result holds for higher genus orientable 
closed surfaces, with singular hyperbolic metrics. 

\begin{thm} \label{tm:hyperb}
Let $\Sigma$ be a closed orientable surface, and let 
$\Gamma$ be the 1-skeleton of a cellular decomposition of $\Sigma$. Let
$\kappa:\Gamma_2\rightarrow (-\infty, 2\pi)$ and 
$\theta:\Gamma_1\rightarrow (0,\pi)$ be two functions. There exists a
hyperbolic metric $h$ with conical singularities on $\Sigma$, with a hyperideal
circle pattern $\sigma$ 
with incidence graph $\Gamma$, intersection angles given by $\theta$,
and singular curvatures given by $\kappa$, if and only if:
\begin{enumerate}
\item $\sum_{f\in \Gamma_2} \kappa(f)>2\pi\chi(\Sigma)$,
\item for any admissible domain $\Omega\subset \Sigma$ : 
\beq \label{eq:cond-omega} 
\sum_{e\in \Gamma_1, e\subset\dr\Omega} \theta(e) \geq (2\chi(\Omega) -
m(\Omega))\pi - \sum_{f\in \Gamma_2, f\subset \Omega} \kappa(f)~, \eeq
with strict inequality except perhaps when $\Omega$ is a face of $\Gamma$.
\end{enumerate}
$h$ and $\sigma$ are then unique. 
\end{thm}

\paragraph{Simpler examples.}

As in the case of flat metrics, the statements become much simpler when one
considers non-singular metrics, or more generally only conical points
with positive singular curvature. 
Indeed under this hypothesis the right-hand side in the
second condition is positive only for disks with at most one boundary
component in a face of $\Gamma$. This explains the simpler form of the
statements obtained e.g. in \cite{bao-bonahon} or in \cite{bobenko-springborn}.

\paragraph{Surfaces with geodesic boundary.}

There is natural extension of Theorem \ref{tm:euclid} and Theorem
\ref{tm:hyperb} to flat or hyperbolic surfaces with geodesic boundary. 
To state it is necessary to extend slightly the definitions given above.
We consider a compact surface with boundary, $\Sigma$, along with a
graph $\Gamma$ embedded into $\Sigma$ so that the union of the closures
of the faces of $\Gamma$ covers $\Sigma$. In other terms, the boundary
of $\Sigma$ is covered by edges and vertices of $\Sigma$. 
We also extend the definition of an admissible domain, as in Definition 
\ref{df:15} above.

\begin{df}
Let $\Gamma$ be a graph embedded in a compact surface with boundary $\Sigma$,
such that $\dr \Sigma\subset \Gamma$. 
An {\it admissible domain} 
in $(\Sigma, \Gamma)$ is a domain $\Omega\subset \Sigma$, containing a face of
$\Gamma$, such that 
$\dr \Omega$ is a finite union of segments which: 
\begin{itemize}
\item have as endpoints vertices of $\Gamma$,
\item either are edges of $\Gamma$, or have interior contained in a face of
  $\Gamma$.
\end{itemize}
We define $\chi(\Omega)$ as the Euler characteristic of $\Omega$, 
and $m(\Omega)$ as the number of maximal
segments in $\dr \Omega$ which are contained either in a face of $\Gamma$ or in
$\dr\Sigma$ (i.e. this excludes connected components of $\dr\Sigma$ which
are boundary components of $\dr\Omega$, because they are closed curves 
rather than segments).
\end{df}

There is a natural notion of circle pattern with incidence graph $\Gamma$,
on $\Sigma$ endowed with a hyperbolic (resp. Euclidean) metric with 
geodesic boundary.  
Such a circle pattern has one circle for each interior vertex of $\Gamma$,
and, for each boundary vertex of $\Gamma$, a half-circle, which intersects
orthogonally the boundary of $\Sigma$.

\begin{thm} \label{tm:hyperb-geod}
Let $\Sigma$ be a compact orientable 
surface with boundary. 
Let $\Gamma$ be the 1-skeleton of a cellular decomposition of $\Sigma$, let
$\kappa:\Gamma_2\rightarrow (-\infty, 2\pi)$ and let
$\theta:\Gamma_1\rightarrow (0,\pi)$ be two functions. There exists a
hyperbolic metric $h$ with conical singularities on $\Sigma$, with a hyperideal
circle pattern $\sigma$ 
with incidence graph $\Gamma$, intersection angles given by $\theta$,
and singular curvatures given by $\kappa$, if and only if:
\begin{enumerate}
\item $\sum_{f\in \Gamma_2} \kappa(f)>2\pi \chi(\Sigma)$.
\item for any admissible domain $\Omega\subset \Sigma$ : 
\beq \label{eq:2-higher-geod}
\sum_{e\subset \dr\Omega\setminus \dr\Sigma} 
\theta(e) \geq \pi(2\chi(\Omega) - m(\Omega)) 
- \sum_{f\subset \Omega} \kappa(f)~, \eeq
with strict inequality except perhaps when $\Omega$ is a face of $\Gamma$.
\end{enumerate}
$h$ and $\sigma$ are then unique. 
\end{thm}

This theorem is a special case --- which is more readily understandable --- of
Theorem \ref{tm:framed-hyperb} below. A similar statement can 
be given for Euclidean metrics with geodesic boundary, based on
Theorem \ref{tm:framed-euclid}, we leave this point to the interested
reader.

\paragraph{Framed circle patterns.}

Theorem \ref{tm:hyperb-geod} above is not completely satisfactory,
in particular it does not apply to circle patterns in the plane 
(even without singularities). This paragraph contains more flexibles
statements, which remain quite simple but for which some additional
definitions are necessary. 

We first give the definition for a disk, so as to simplify the notations a
little. A more general definition is used below, but the only difference is in
the number of boundary components.

\begin{df}
A {\it framed circle pattern}, on a Euclidean disk with conical 
singularities, is a hyperideal circle pattern, with singularities at the 
center of the dual circles, along with a closed polygonal line 
made of geodesics segments $g_1, \cdots, g_n$, such that:
\begin{itemize}
\item the boundary of the exterior face is made of segments of the
circles $C_1, \cdots, C_n$ (in this order), 
\item for each $i=1, \cdots, n$, the circle $C_i$ intersects the 
interior of the segment $g_i$ in two points,
\item for each $i\in \Z/n\Z$, there is a  circle centered at the intersection
  point of $g_i$ and $g_{i+1}$ which is orthogonal to both $C_i$ and
  $C_{i+1}$.
\end{itemize}
\end{df}

\begin{figure}[ht]
\centerline{\psfig{figure=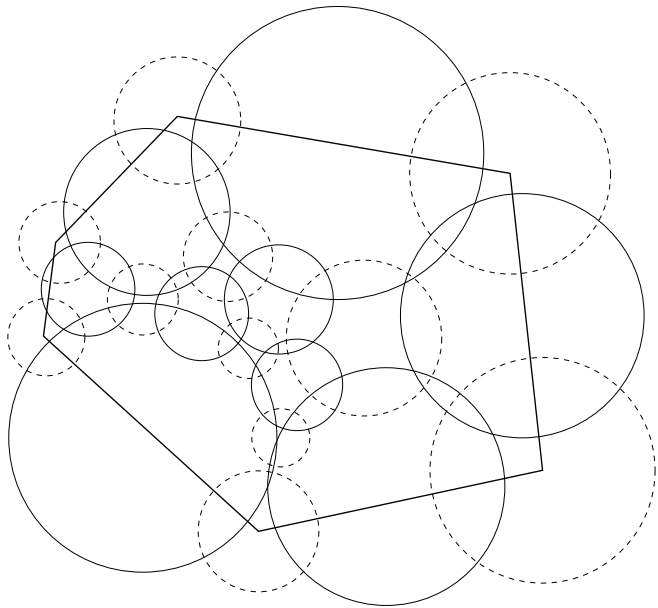,height=8cm}}
\caption{A framed hyperideal circle pattern.}
\end{figure}

Note that in the second condition the intersection of the boundary segments
with the corresponding circles are not at the endpoints of the segments, in
other terms the ``dual'' circles centered at the vertices of the boundary
polygon are not reduced to points.

Given a framed circle pattern, the intersection angles between the oriented
  geodesic segments $g_i$ and $g_{i+1}$, for $i=1,\cdots, n$, are called its
  {\it polygonal angles}, they are equal to $\pi$ minus the angle at the
  intersection point of the interior of the polygon with edges the $g_i$. 
The intersection angles between $g_i$ and $C_i$, for
$i=1,\cdots, n$, are called the {\it boundary angles}, they are defined as the
angle at the intersection point of the domain of the interior of $C_i$ which
is {\it outside} the polygon with edges the $g_i, i=1,\cdots,n$. 

\paragraph{The extended graph of a framed circle pattern.}

It is also convenient to consider a combinatorial data slightly more elaborate
than the incidence graph $\Gamma$ of the circle pattern. The {\it extended
  incidence graph} $\Gamma'$ of a framed circle pattern has, in addition to
the incidence graph $\Gamma$:
\begin{itemize}
\item one vertex for each of the $g_i, i=1,\cdots, n$. Those are called the 
{\it boundary vertices} of $\Gamma'$, and they form a set denoted by
$\Gamma'_{0,\dr}$, 
\item one edge going from the vertex corresponding to $g_i$ to the vertex
  corresponding to $g_{i+1}$, for each $i=1,\cdots, n$, those edges also form
  a set denoted by $\Gamma'_{1,\dr}$,
\item one edge going from the vertex corresponding to $C_i$ to the vertex
  corresponding to $g_i$, for each $i=1,\cdots, n$. 
\end{itemize}
Clearly $\Gamma'$ depends only on $\Gamma$ (and not on other data from the
framed circle pattern) so we will call $\Gamma'$ the {\it extended graph of
  $\Gamma$}. 

Since $\Gamma'$ is not embedded in a closed surface but in a ``surface with
boundary'', it is relevant to define what is meant by the dual graph
$\Gamma'^*$: it is the graph which has one vertex for each face of $\Gamma'$,
one edge between two vertices of $\Gamma'^*$ when the corresponding faces of
$\Gamma'$ are adjacent (in particular there is no edge corresponding to the
boundary edges of $\Gamma'$), and one face for each vertex of
$\Gamma$ (i.e. not for the boundary vertices of $\Gamma'$).

Just as the edges of $\Gamma$ carry naturally a ``weight'', which is the
intersection angle between the corresponding circles of the pattern, the
additional edges of $\Gamma'$ carry natural weights:
\begin{itemize}
\item For the edge going from the vertex corresponding to $g_i$ to the vertex
  corresponding to $g_{i+1}$, it is the ``polygonal'' angle between $g_i$ and
  $g_{i+1}$, as defined above.
\item For the edge going from the vertex corresponding to $C_i$ to the vertex
  corresponding to $g_i$, it is the ``boundary'' angle between $C_i$ and
  $g_i$. 
\end{itemize}
We call this data the ``extended intersection angles'' of the circle pattern.
The considerations made here are of course not limited to topological 
disks, one can
also define in the same way framed circle patterns on any compact surface with
boundary. 

One additional notation is necessary: given a compact surface with
boundary $\Sigma$ and a graph $\Gamma$ embedded in $\Sigma$, and
given an admissible domain $\Omega$ in $(\Sigma, \Gamma')$, we denote by
$n(\Omega)$ the number of segments in $\dr\Omega\cap \dr\Sigma$. Note
that the term ``segment'' excludes boundary components which are 
topologically circles.

We can now state an analog of Theorem \ref{tm:euclid} 
for framed circle patterns
on singular Euclidean surfaces with boundary.

\begin{thm} \label{tm:framed-euclid}
Let $\Sigma$ be compact orientable surface with boundary.
Let $\Gamma$ be the 1-skeleton of a cellular decomposition of $\Sigma$,
and let $\Gamma'$ be the extended graph of $\Gamma$.
Let $\kappa:\Gamma_2\rightarrow (-\infty, 2\pi)$, let
$\theta:\Gamma'_1\rightarrow (0,\pi)$ be functions. There exists an
Euclidean metric $h$ with conical singularities on $\Sigma$, with a hyperideal
circle pattern $\sigma$ 
with incidence graph $\Gamma$, extended intersection angles given by $\theta$,
and singular curvatures given by $\kappa$, if and only if:
\begin{enumerate}
\item $\sum_{f\in \Gamma_2} \kappa(f)= 2\pi \chi(\Sigma)
- \sum_{v\in \Gamma'_{1,\dr}} \theta(v)$,
\item for any admissible domain $\Omega$ in $(\Sigma, \Gamma')$: 
$$ \sum_{e\in \Gamma'_1, e\in \dr\Omega} 
\theta(e) \geq (2\chi(\Omega) 
- m(\Omega) - n(\Omega))\pi - 
\sum_{f\in \Gamma_2, f\subset \Omega} \kappa(f)~, $$
with strict inequality except that:
\begin{itemize}
\item equality is possible when $\Omega$ is a face of $\Gamma$, 
\item equality is required when $\Omega$ is the union of all faces of $\Gamma$
  (except the exterior face). 
\end{itemize}
\end{enumerate}
$h$ is then unique up to homotheties, and $\sigma$ is unique.
\end{thm}

Note that the second exception in condition (2) is a reformulation 
of condition
(1), we keep the two conditions for the sake of unity with the other similar
statements. 
There is a similar result for framed hyperbolic circle patterns. 

\begin{thm} \label{tm:framed-hyperb}
Let $\Sigma$ be a compact orientable surface with non-empty boundary.
Let $\Gamma$ be the 1-skeleton of a cellular decomposition of $\Sigma$,
and let $\Gamma'$ be the extended graph of $\Gamma$.
Let $\kappa:\Gamma_2\rightarrow (-\infty, 2\pi)$, let
$\theta:\Gamma'_1\rightarrow (0,\pi)$ be functions. There exists a
hyperbolic metric $h$ with conical singularities on $\Sigma$, with a hyperideal
circle pattern $\sigma$ 
with incidence graph $\Gamma$, extended intersection angles given by $\theta$,
and singular curvatures given by $\kappa$, if and only if:
\begin{enumerate}
\item $\sum_{f\in \Gamma_2} \kappa(f)> 2\pi - \sum_{v\in \Gamma_{1,\dr}}
  \theta(v)$,
\item for any admissible domain $\Omega$ in $(\Sigma, \Gamma')$: 
\beq \label{eq:cond-framed}
\sum_{e\in \dr\Omega} \theta(e) \geq 
(2\chi(\Omega) - m(\Omega) - n(\Omega))\pi - 
\sum_{f\subset \Omega} \kappa(f)~, \eeq
with strict inequality except perhaps when $\Omega$ is a face of $\Gamma$.
\end{enumerate}
$h$ and $\sigma$ are then unique.
\end{thm}
 
When $\kappa \equiv 0$ and $\Sigma$ is a disk, 
both Theorem \ref{tm:framed-euclid} and Theorem
\ref{tm:framed-hyperb} are direct consequences of a theorem of Bao and
Bonahon \cite{bao-bonahon} on the dihedral angles of hyperideal polyhedra.


\section{Outline of the proof}

Since the proofs of the main results are somewhat intricate and
have some technical parts, it appears helpful to give first a 
brief outline of the way they proceed.

Theorem \ref{tm:framed-euclid} follows from Theorem \ref{tm:euclid}, 
and Theorem \ref{tm:framed-hyperb} from Theorem \ref{tm:hyperb}, so
we concentrate on those statements here.

\subsection{A deformation method}

The proof of the two main results follow a so-called deformation
method, which has been classical at least since the work of Aleksandrov
\cite{alex} on convex polyhedra. 
First we choose a graph $\Gamma$ embedded in a closed surface. 
Thenn we consider a space
$\cC(\Gamma)$ of circle patterns with combinatorics given by $\Gamma$, and a
space $\cD(\Gamma)$ of intersection data and singular curvatures at the
singular points satisfying the hypothesis of the theorem.
There is also a natural map
$\Phi_\Gamma$ defined on $\cC(\Gamma)$ sending
a circle pattern to its intersection angles and singular curvatures. It is then
necessary to check the following points:
\begin{enumerate}
\item the image of $\Phi_\Gamma$ is contained in $\cD(\Gamma)$, i.e. the
  conditions in the theorem are  necessary,
\item the spaces $\cC(\Gamma)$ and $\cD(\Gamma)$ are differentiable manifolds
  of the same dimension,
\item the map $\Phi_\Gamma$ has injective differential at each point of
  $\cC(\Gamma)$. This 
  can be formulated as an {\it infinitesimal rigidity} statement: it is not
  possible to deform infinitesimally a circle pattern without changing either
  the intersection angles or the singular curvatures.
\item $\Phi_\Gamma$ is proper. 
  This translates as a {\it compactness} property: if a
  sequence of circle patterns is such that the intersection angles and the
  singular curvatures converge to a ``good'' limit, than it has a converging
  subsequence. 
\item $\cD(\Gamma)$ is contractible (actually $\cD(\Gamma)$ is the interior of
  a polytope). 
\end{enumerate}
Point (1) is mostly elementary, and is proved in section 3.
Point (3) is perhaps the key part of the
proof, and is done -- for each of the situations considered -- in section
5, points (2) and (5) are proved as consequences of the methods
developed there. 
Point (4) also demands some efforts, it is treated in section 6. 

It is well known, and the consequence of simple topological arguments,
that those points imply that $\Phi_\Gamma$ is a covering of $\cD(\Gamma)$ by
$\cC(\Gamma)$. In other terms, for each choice of $\Gamma$, there is an integer
$N_\Gamma \in \N$ such that each point of $\cD(\Gamma)$ has exactly $N_\Gamma$
inverse images by $\Phi_\Gamma$. It follows from the arguments used to prove
the infinitesimal rigidity that $N_\Gamma$ is at most $1$. 
It remains to prove that $N_\Gamma=1$ for 
any choice of $\Gamma$. 

We show in section 7 that, for Theorem
\ref{tm:hyperb}, this follows from 
important results in geometric topology, in particular the Orbifold
Hyperbolization Theorem \cite{boileau-porti}. To show the same result
for Euclidean surfaces, we use in section 7 an approximation argument, 
based on another use of our compactness statement, Lemma 
\ref{lm:compact}, which
permits us to obtain a circle pattern on an Euclidean surface with
conical singularities as a limit of circle patterns on hyperbolic surfaces
(after multiplication of the hyperbolic metrics by some coefficients).

We can now consider with a little more details the main points of the proof,
in particular points (3) and (4) of the outline above. 

\subsection{Compactness}

This part can be found in section 6. There are different ways to prove the
kind of compactness property needed here, in particular in the context of
hyperbolic polyhedra (see e.g. \cite{RH,CD,shu,cpt}) or ideal circle patterns
on singular surfaces (see e.g. \cite{Ri2,bobenko-springborn}). However the 
condition appearing in the main statements here, concerning linear
inequalities for each admissible domain, is more elaborate than the
condition in earlier references, which only concerned topological disks,
and this is partly reflected in the fact that the compactness argument 
is also a little more elaborate.
The approach followed here, however, is elementary, without reference to the 
underlying 3-dimensional
hyperbolic geometry, and only in terms  of circle patterns.

\subsection{Hyperbolic polyhedra and circle patterns}

As the reader might have already well understood, the results on circle
patterns presented here rely heavily on tools from 3-dimensional hyperbolic
geometry, and in particular on some properties of the volume of hyperideal
hyperbolic polyhedra. This is developed in section 4, we give a short outline
here. The ideas go back to hyperbolic geometry construction for circle
packings (in particular \cite{thurston-notes,CdeV,bragger}) 
and for ``ideal'' circle patterns (see in particular
\cite{Ri2,leibon1,bobenko-springborn}). 

Let's first describe the topological aspects of the construction. Given an
embedded graph on a surface $\Sigma$ 
we define a 3-dimensional cell complex $\cS$, 
which is a cone over the graph $\Gamma^*$ dual to $\Gamma$: $\cS$ has one
3-cell for each face of $\Gamma^*$ (that is, for each
vertex of $\Gamma$),
one 2-face for each face of $\Gamma^*$ 
and one for each edge of $\Gamma^*$
and one edge for each edge of $\Gamma^*$ 
and one for each vertex of $\Gamma^*$. 
Its vertices are the vertices of $\Gamma^*$
plus one, which will be called here the ``central'' vertex of $\cS$. 

To a  hyperideal circle pattern is naturally associated a non-complete
hyperbolic metric on $\cS$, which might have conical singularities along the
edges of $\cS$ corresponding to the vertices of $\Gamma^*$.
This metric is obtained by giving to each 3-cell of $\cS$
the hyperbolic metric on a hyperideal polyhedron, such that the central vertex
is either ideal -- when considering Euclidean circle patterns -- or strictly
hyperideal -- when considering hyperbolic circle patterns. It is such that
there is a well-defined notion of exterior dihedral angle at the edges of
$\Gamma^*$ and this dihedral angle is equal to the
intersection angle between the corresponding circles of the pattern.

The ``hyperideal'' circle pattern on the Euclidean (resp. hyperbolic) surface
is then recovered by considering a family of horosphere centered at the
``central'' vertex (resp. the intersection with each 3-cell of the hyperplane
dual to the ``central'' vertex) and projecting on it, in the direction of the
``central'' vertex, the circle which is the boundary at infinity 
of the face opposite the ``central'' vertex. (This is explained in section
4). 

The construction works both ways: given a hyperideal circle pattern on a
singular surface, one can recover a ``polyhedral'' object as described
above. Moreover, the total angles around the singularities on the Euclidean
(resp. hyperbolic) surface are equal to the total angles around the edges
going to the ``central'' vertex, and the intersection angles between the
circles are equal to the exterior dihedral angles at the edges which do not
contain the ``central'' vertex. 
So proving Theorem \ref{tm:euclid} and Theorem \ref{tm:hyperb}
is equivalent to proving the
existence and uniqueness of those polyhedral objects with given combinatorics,
conical singularities along some edges, and dihedral angles.

\subsection{The infinitesimal rigidity}

A key idea of the proof -- as in some of the references cited above on circle
packings or ``ideal'' circle  patterns -- is to consider a wider range of
hyperbolic metrics on the 3-cells of $\cS$, including some which do not allow
isometric gluings of the cells or for which this gluing leads to ``bad''
singularities on the edges going to the central vertex. In the deformation
method used here, however, it is only used in the proof of the infinitesimal
rigidity of circle patterns. 

The main point is that, among all those more general choices of hyperbolic
metrics on the 3-cells of $\cS$, those which correspond to circle patterns are
characterized as maxima of a function -- which is simply the sum of the
hyperbolic volumes of the 3-cells -- under some linear constraints on the
dihedral angles. Since this function is strictly concave, the maxima are
isolated, which proves the infinitesimal rigidity.

\subsection{The results on surfaces with boundary}

All the results concerning surfaces with boundary -- Theorems
\ref{tm:framed-euclid} and \ref{tm:framed-hyperb} -- are consequences
of Theorems \ref{tm:euclid} and \ref{tm:hyperb} 
using a doubling argument to pass from a surface with boundary to a closed
surface. 


\section{Necessary conditions}

\subsection{Elementary properties of circle patterns}

It is the object of this section to prove that the conditions appearing in the
statements of the main theorems given in the introduction are necessary. We
first state some preliminary properties of some simple patterns of circles 
on a surface with at most one singularity, and the next subsection will
show how similar arguments indicate that the conditions in the 
two theorems concerning circle patterns on 
closed surfaces are necessary. 

\paragraph{Circles intersecting at a point.}

The first statement describes a very simple situation, with a set of circles
meeting at the only singular point of an Euclidean surface.

\begin{lemma} \label{lm:bouquet}
Let $C_E(\kappa)$ (resp. $C_H(\kappa)$) be the complete Euclidean 
(resp. hyperbolic) metric on the plane with one 
conical singularity, with singular curvature $\kappa$ ($\kappa\in (-\infty,
2\pi)$), and let $x_0$ be its singular point. 
Let $C_1, \cdots, C_N$ be a sequence of circles on $C_E(\kappa)$,
containing $x_0$, such that the oriented tangents to the $C_i$ at $x_0$ 
appear in cyclic order. 
Let $\theta_i$ be the angle between $C_i$ and $C_{i+1}$ 
(with the convention that $C_{N+1}=C_1$). Then
$$ \sum_{i=1}^N \theta_i = 2\pi - \kappa~. $$
\end{lemma}

\begin{figure}[ht]
\centerline{\psfig{figure=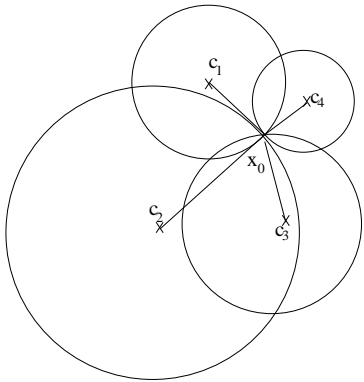,height=5cm}}
\caption{Lemma 3.1}
\end{figure}

\bpv
Let $c_1, \cdots, c_N$ be the centers of the $C_i$. For $i=1,\cdots, N$, 
let $s_i$ be the oriented geodesic segment starting at $x_0$ and ending
at $c_i$. Then the oriented angle at $x_0$ between $C_i$ and $s_i$
is equal to $\pi/2$, so that the oriented angle between $s_i$ and 
$s_{i+1}$ is equal to $\theta_i$. So $\sum_{i=1}^N \theta_i$ 
is equal to the total angle at the conical point $x_0$, which is equal 
to $2\pi-\kappa$.
\epv

\subsection{Necessary conditions}

\paragraph{The conditions of Theorem 1.4 are necessary.}

We now check that the conditions appearing in Theorem
\ref{tm:euclid}
are necessary. Condition (1) is a restatement of the Gauss-Bonnet
theorem. We will see that condition (2) follows from arguments similar
to those used in the
proof of Lemma \ref{lm:bouquet}, with the equality case a direct consequence
of that statement. 

Suppose that $\Gamma$ is the incidence graph of a circle pattern on a surface
$\Sigma$ with
a flat metric with conical singularities, as in Theorem \ref{tm:euclid}, and
let $\Omega$ be an admissible domain in $(\Sigma, \Gamma)$. Suppose first that
$\dr \Omega$ is contained in the 1-skeleton of $\Gamma$, so that $\dr\Omega$
is a disjoint union of closed curves $S_1, \cdots, S_k$, 
each of which is made of a sequence of
edges of $\Gamma$. To each connected component $S_i$ is associated a 
sequence of circles, $C^i_1, \cdots, C^i_{N_i}$, 
such that $C^i_j$ intersects $C^i_{j+1}$
(including $j=N_i$, then $C^i_{N_i}$ intersects $C^i_1$). 

\begin{figure}[ht]
\centerline{\psfig{figure=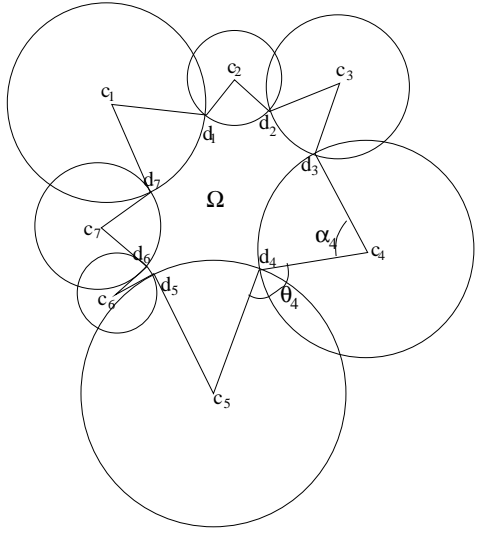,height=7cm}}
\caption{Chains of intersecting circles}
\end{figure}

For each $j=1,\cdots,N$, let $d^i_j$ be the intersection point 
between $C^i_j$ and $C^i_{j+1}$ which is on the same side of 
$\gamma^i_1\cup\cdots\cup \gamma^i_{N_i}$ as
$\Omega$. We call $\theta^i_j$ the angle at 
$d^i_j$ between $[d^i_j,c^i_j]$ and
$[d^i_j,c^i_{j+1}]$, and $\alpha^i_j$ the 
angle at $c^i_j$ between $[c^i_j,d^i_{j-1}]$ and
$[c^i_j,d^i_j]$. 

Let $K$ be the sum of the singular curvatures of the metric 
in the domain bounded by the polygonal curves just described, i.e. 
$K$ is the sum of the singular curvatures attached by the 
function $\kappa$ to the faces of $\Gamma$ which are contained in 
$\Omega$. 
The Gauss-Bonnet theorem, applied to the polygonal curve made of the
segments $[c^i_j,d^i_j]$ and $[d^i_j,c^i_{j+1}]$, yields that:
\beq \label{eq:gb} 
\sum_{i=1}^k\sum_{j=1}^{N_i} 
(\theta^i_j-\pi) + (\pi-\alpha^i_j) = 2\pi \chi(\Omega) - K~, \eeq
so that:
$$ \sum_{i=1}^k\sum_{j=1}^{N_i} 
\theta^i_j = 2\pi \chi(\Omega) -K + 
\sum_{i=1}^k\sum_{j=1}^{N_i} \alpha^i_j~, $$
and the result follows. 

Suppose now that the equality is attained in condition (2) of Theorem 
\ref{tm:euclid}. It follows quite directly from the argument given 
above that, for each polygonal curve $S_i$ corresponding to one of the
boundary component of $\Omega$, all the angles $\alpha^i_j$ are zero,
so that the situation is exactly the one described in 
Lemma \ref{lm:bouquet} -- in particular $\Omega$ is a face of $\Gamma$,
as claimed.

Consider now the more general case where $\Omega$ is an admissible
domain with boundary containing, in addition of edges of $\Gamma$,
$m$ curves which are contained in faces of $\Gamma$. The same procedure
can then be applied, associating a circle to each vertex of $\dr\Omega$.
In addition we also associate a circle to each segment of $\dr\Omega$
which is contained in a face $f$ of $\Gamma$: namely, the dual circle
corresponding to $f$. This circle intersects
orthogonally the circles corresponding to the two vertices of $\dr\Omega$
which are adjacent to it. Therefore, following the proof given above
yields the same formula, except that to each segment of $\dr\Omega$
which is contained in a face of $\Gamma$ correspond two angles equal
to $\pi/2$ in equation (\ref{eq:gb}), so that this left-hand side
is increased by $m(\Omega)\pi$. 
This proves that the conditions in Theorem \ref{tm:euclid} are necessary.

\paragraph{Hyperbolic surfaces}

The arguments showing that the conditions of Theorem \ref{tm:hyperb} and of
Theorem \ref{tm:framed-hyperb} are necessary then proceed exactly as in the
proof of the corresponding Euclidean statements. 
Condition (1) is a consequence
of the Gauss-Bonnet theorem, since the area of the surfaces appears in
the formula. The fact that the second condition is also necessary 
can also be proved as in the Euclidean case, with only one difference,
namely that the area of the domain which is considered -- bounded by a 
disjoint union of polygonal curves -- also comes in the Gauss-Bonnet
formula, but with a sign which does not disturb the proof.


\section{Some hyperbolic geometry}

\paragraph{A short introduction.}

As already outlined in section 2, it will be useful to consider some
3-dimensional hyperbolic geometric objects, of a polyhedral nature, associated
to the circle patterns on which our attention is focussed. Those are obtained
by gluing a hyperideal pyramid for each face of the graph $\Gamma^*$ dual to
$\Gamma$, with a vertex ``outside'' $\Gamma^*$ which is either ideal -- when
considering Euclidean circle patterns -- or strictly hyperideal -- when
considering hyperbolic circle patterns. 

It is first necessary to recall the definitions of hyperideal polyhedra, and a
basic result of Bao and Bonahon on their possible dihedral angles, and then
some further (rather elementary) definitions concerning their volume and edge
lengths. The next subsection recalls the Schl\"afli formula for hyperideal
polyhedra, and the concavity property for their volume. The third and last
subsection describes precisely the relation between a hyperideal circle
pattern and the underlying 3-dimensional cell complex with its hyperbolic
metric, as well as the more general geometric objects which will appear in the
proof of the infinitesimal rigidity of circle patterns.

\subsection{Ideal and hyperideal polyhedra}

\paragraph{Hyperideal polyhedra.}

When considering hyperideal polyhedra, it is helpful to make use of the Klein,
or projective, model of $H^3$; it is a map from $H^3$ to the unit ball in
Euclidean 3-space, which sends hyperbolic geodesics to geodesic segments (see
e.g. \cite{GHL}). Compact hyperbolic polyhedra correspond in this model to
polyhedra in the open unit ball, while {\it ideal} polyhedra correspond to
polyhedra in the unit ball with all their vertices on the unit sphere, and
{\it hyperideal} polyhedra correspond to polyhedra with all their vertices
outside the open unit ball, but with all their edges intersecting this ball. 
A vertex of a hyperideal polyhedron is {\it ideal} if it sits on the unit
sphere, {\it strictly hyperideal} otherwise. A hyperideal polyhedron is {\it
  strictly hyperideal} if all its vertices are strictly
hyperideal.

To each edge of a compact/ideal/hyperideal polyhedron, we can associate its
{\it dihedral angle}, which is the angle between the two faces of the
polyhedron which meet at that edge. We will always consider here the {\it
  exterior dihedral angle}, which is defined as $\pi$ minus the angle measured
in the interior of the polyhedron. 

The possible exterior dihedral angles of hyperideal hyperbolic polyhedra have
been described recently by Bao and Bonahon, extending previous results of
Andreev \cite{andreev-ideal} and Rivin \cite{rivin-annals} concerning ideal
polyhedra. It can also be obtained as a consequence of (an extension of)
a result of Rivin and
Hodgson \cite{Ri,RH} on compact hyperbolic polyhedra, see
\cite{rousset1}. 
The statement uses the notion of {\it admissible open path} in a graph, which
is defined as a simple path $\gamma$, beginning and ending at vertices of a
face $f$, and not contained in the boundary of $f$.

\begin{thm}[Bao, Bonahon \cite{bao-bonahon}] \label{tm:bb}
Let $\Gamma$ be the 1-skeleton of a polytopal cellular decomposition
of $S^2$. Let $w:\Gamma_1\rightarrow (0,\pi)$. There exists a hyperideal
polyhedron with combinatorics dual to $\Gamma$ and (exterior) dihedral
angles given by $w$ if and only if:
\begin{itemize}
\item for each simple closed curve $\gamma$ in $\Gamma$, the sum of the values
  of $w$ on the edges of $\gamma$ is at least $2\pi$, with strict equality
  unless $\gamma$ bounds a face,
\item for each admissible open path $\gamma$ in $\Gamma$, the sum of the
  values of $w$ over the edges of $\gamma$ is strictly larger than $\pi$.
\end{itemize}
This hyperideal polyhedron is then unique (up to global hyperbolic
isometries). 
\end{thm}

The equality case in the first condition corresponds exactly to the faces of
$\Gamma$ corresponding to ideal vertices of the polyhedron.

\paragraph{Edge lengths of ideal and hyperideal polyhedra.}

There is a natural projectively defined duality between points which are
outside the closed unit ball in $\R^3$ and planes intersecting the unit ball
(which are themselves identified with hyperbolic planes by the projective
model of $H^3$). It is related to a duality between the hyperbolic and the de
Sitter space, see
e.g. \cite{coxeter-de-sitter,gelfand5,thurston-notes,Ri,RH,shu}, but
this aspect will not be used here. It associates to each point $x$ outside the
unit ball the plane $p$ containing the intersection with the unit sphere of the
cone with vertex at $x$ which is tangent to this unit sphere. It has a
striking property: the intersection of $p$ with the unit ball, considered as a
hyperbolic plane, is orthogonal to the intersections with the unit ball of all
lines going through $x$, also considered as hyperbolic geodesics. The
intersection of $p$ with the unit ball, considered as a hyperbolic plane, is
called the hyperbolic plane {\it dual} to $x$.

Note that, given a hyperideal polyhedron, the hyperbolic planes dual to the 
strictly hyperideal vertices are disjoint. This is a direct consequence, using
the projective definition of the duality, of the fact that the line going
through two vertices always intersects the unit ball in $\R^3$.

Given a hyperideal hyperbolic polyhedron $P$, the associated {\it truncated
polyhedron} $P_t$ 
(defined in \cite{bao-bonahon}) is obtained by cutting off each
strictly hyperideal vertex by its dual plane. Its vertices are ``usual''
points of $H^3$, except for the ideal vertices of $P$, which remain vertices
of $P_t$. The faces of $P_t$ are either truncated faces of $P$, or ``new''
faces, which are orthogonal to all adjacent faces. 

Given an edge of $P$ which has as endpoints two strictly hyperideal
vertices, its {\it length} is defined as the length of the corresponding edge
of $P_t$. To define the lengths of the edges with one or two ideal vertices as
their endpoints, it is necessary to choose, for each ideal vertex $v$ of $P$,
a horosphere $h$ ``centered'' at $v$; one then defines the length of an edge
of $P$, with endpoints $v_1$ and $v_2$ as the oriented distance between the
plane dual to $v_1$ (if $v_1$ is strictly hyperideal) or the horosphere chosen
for $v_1$ (if $v_1$ is ideal) to the plane dual to $v_2$ or the horosphere
chosen for $v_2$. Replacing a horosphere centered at a vertex $v$ by another
one clearly adds a constant -- the oriented distance between the two
horospheres -- to all the lengths of the edges with endpoint $v$, so that, if
$P$ has $n$ ideal vertices, the
lengths of the edges of $P$ are defined up to the addition of $n$ constants,
one for each ideal vertex (the constants being added to the lengths of the
edges containing the respective vertex).

Obviously, what has been said of hyperideal polyhedra is also true of
hyperideal triangles, which are defined in the same way using the projective
model of $H^2$. This leads to the following elementary proposition.

\begin{prop} \label{pr:gluing}
Two hyperideal triangles are isometric if and only if they
have the same number of ideal vertices and the same edge
lengths. Any non-trivial first-order deformation of a hyperideal triangles
induces a non-zero first-order variation of its edge lengths.
\end{prop}

\begin{proof}
Another formulation of the first statement is that hyperideal triangles are
uniquely determined, up to global isometries, by their edge lengths. Consider
first a strictly hyperideal triangle (i.e. with no ideal vertex). Using the
projective model of $H^2$, one can choose arbitrarily the position of the
first vertex, $v_1$; the fact that the edge going from $v_1$ to $v_2$ has given
length is equivalent, in the projective model, to the fact that $v_2$ is on
a segment of an ellipse tangent to the unit circle, with endpoints at the
tangency points. Different choices of
points on that segment are equivalent up to an isometry.
\begin{figure}[ht]
\centerline{\psfig{figure=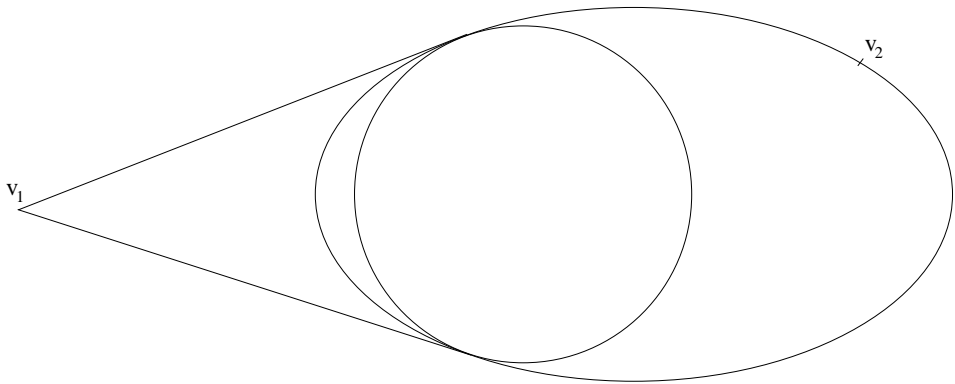,height=5cm}}
\caption{Possible positions of $v_2$ given $v_1$ and the edge length.}
\end{figure}

Then, given the position of $v_1$ and $v_2$, the position of $v_3$ is
uniquely determined (up to an isometry) 
by the condition that it is contained in two ellipses, one
determined by the length of the edge between $v_1$ and $v_3$, the other by the
length of the edge from $v_2$ to $v_3$. 

The same argument can be used if $v_1$ is an ideal vertex, using the fact
that, given $v_1$ and a horosphere $h_1$ ``centered'' at $v_1$, the 
possible positions
of $v_2$ such that the length between $v_1$ and $v_2$ (measured with respect
to the choice of $h_1$) is given, are on a circle tangent to the unit circle
at $v_1$. The same argument works when $v_1$ and $v_2$ are ideal points, then,
given the choice of $v_1$ and $h_1$, the position of $h_2$ can be chosen
arbitrarily on the unit circle but determines the horosphere $h_2$
``centered'' at $v_2$. Finally, if all three vertices are ideal, the statement
is trivial, since all ideal triangles are isometric, and there is only one
possible set of edge lengths (up to the addition of three arbitrary constants,
corresponding to the three idea vertices).

The proof of the second part of the 
statement can be done exactly along the same lines, by
considering a first-order variation of the edge lengths and the implied
conditions on the first-order displacement of the vertices.
\end{proof}

\subsection{The volume of hyperideal polyhedra}

It follows from the definition of the truncated polyhedron $P_t$ 
associated to a hyperideal polyhedron $P$ that $P_t$ always has finite volume,
leading us to define the {\it volume} of $P$ as the volume of $P_t$ (of course
the intersection of hyperbolic 3-space with the interior of $P$, in the
projective model of $H^3$, does not have finite volume unless all vertices of
$P$ are ideal).

\paragraph{The Schl\"afli formula.}

Consider first a compact polyhedron $P$.
The classical Schl\"afli formula (see e.g. \cite{milnor-schlafli}) gives the
first-order variation of the volume of $P$ under a first-order deformation:
\beq \label{eq:schlafli} 
dV = \frac{1}{2}\sum_e l_e d\theta_e~, \eeq 
where the sum is over the edges of $P$, $l_e$ is the length of edge $e$, and
$\theta_e$ is its (exterior) dihedral angle. Note that no minus sign is needed
in the formula because we consider the exterior dihedral angles. 

This formula remains valid for polyhedra having some ideal vertices (see
\cite{milnor-schlafli,rivin-annals}) if one uses the definition of the edge
lengths given above, defined up to the addition of one constant for each ideal
vertex, and under the condition that the ideal vertices remain ideal in the
first-order deformation considered. Note that formula (\ref{eq:schlafli}) does
not depend on the choice of the horospheres at the vertices, because the sum
of the exterior dihedral angles of the edges containing an ideal vertex is
always equal to $2\pi$ (this is seen by considering the link of the ideal
vertex, which is an Euclidean polygon with edge lengths equal to those exterior
dihedral angles).

Since both the volume and the edge lengths of a hyperideal polyhedron are
defined by reference to the associated truncated polyhedron, it follows that
(\ref{eq:schlafli}) also holds for hyperideal polyhedra.

\paragraph{A technical lemma on convex functions.}

We need below the following elementary statement on properties of concave
function. It is taken from \cite{hphm}, and we leave the proof, which is
simple, to the reader.

\begin{remark} \label{rk:concavity}
Let $\Omega\in \R^N$ be a convex subset, and let $f:\Omega\rightarrow
\R$ be a smooth, strictly concave function. Let
$\rho:\R^N\rightarrow \R^p$ be a 
linear map, with $p<N$, and let $\Omegab:=\rho(\Omega)$. Define a
function:
$$
\begin{array}{cccc}
\fb: & \Omegab & \rightarrow & \R \\
& y & \mapsto & \max_{x\in \rho^{-1}(y)} f(x) 
\end{array}
$$
Then $\Omegab$ is convex, and $\fb$ is a smooth, strictly
concave function on $\Omegab$.
\end{remark}

\paragraph{A key convexity property.}

The key point of this section is the fact, already pointed out in \cite{hphm}
and used in \cite{rcnp}, that the volume of hyperideal polyhedra is a
strictly concave function of the dihedral angle. 
We give here a very short
outline of the proof for the sake of the curious reader, but refer to
\cite{hphm} for a complete proof. The ideas used here are close to those in
\cite{bragger,rivin-annals,hphm}. 

\begin{lemma} \label{lm:concave}
Let $\Gamma$ be the 1-skeleton of a (polytopal) cellular decomposition of the
sphere, and let $V\subset \Gamma_2$ be a subset of the set of its faces. Let
$P_\Gamma$ be the space of hyperideal polyhedra with combinatorics dual to
$\Gamma$, such that the ideal vertices correspond exactly to the faces of
$\Gamma$ which are in $V$. The volume, considered as a function on $P_\Gamma$,
is a strictly concave function of the dihedral angles. 
\end{lemma}

\begin{proof}[Main ideas of the proof]
We consider here only the case of strictly hyperideal polyhedra, otherwise the
proof is similar but slightly more complicated, as should be clear from
section 5. Recall that, by Theorem
\ref{tm:bb}, those hyperideal polyhedra are parametrized by their
possible dihedral angles.

It is necessary to consider first the simplest case, the hyperideal
simplices. One can prove, along the ideas of the (direct) proof of Proposition
\ref{pr:gluing}, that hyperideal simplices are rigid: any non-trivial
first-order deformation of a hyperideal simplex induces a non-zero first-order
variation of its edge lengths. 

Using the Schl\"afli formula, this implies that the Hessian of the volume
function with respect to the dihedral angles 
-- which by (\ref{eq:schlafli}) is the matrix of the derivatives of the edge
lengths with respect to the dihedral angles -- is non-degenerate, since a
non-zero vector in its kernel would precisely correspond to a first-order
variation of the dihedral angles inducing no first-order variation of the edge
lengths. 

Therefore, the signature of the Hessian of the volume (with respect to the
dihedral angles) is constant over the space of hyperideal simplices (which is
connected by Theorem \ref{tm:bb}). To prove that the volume is a
strictly concave function of the dihedral angles, it is therefore sufficient
to check, by an explicit computation, that it is true for a well chosen
simplex (for instance one with maximal symmetry), which is done in the
appendix of \cite{hphm}. 

We now consider a hyperideal polyhedron $P$. We choose a vertex $v$ of $P$,
and add edges to the faces of $P$ so as to subsdivide all faces of $P$ into
triangles, with all faces containing $P$ subdivided by adding only edges
containing $v$. This yields a new polyhedron $P'$, with the same edges as $P$
plus additional edges where the exterior dihedral angle is $0$. 
Then we define a triangulation of the interior of $P'$ into simplices, by
simplices with one vertex at $v$ and three vertices which are the vertices of
a face of $P'$ which does not contain $v$.
Let $S_1, \cdots, S_N$ be
the simplices in this decomposition. Then the $S_i$ are the cells of a natural
cell complex, with two kind of edges: the ``exterior'' edges, which are 
either edges of $P$ or contained in faces of $P$, and the ``interior'' edges,
which are in the interior of $P$ and have $v$ as one of their endpoints.

For each $i=1,\cdots,N$, let $\cA_i$ be the space of possible {\it interior}
dihedral angles
defined on $S_i$ by identifying it to a strictly  hyperideal simplex. Then let
$\cA:=\Pi_{i=1}^N\cA_i$. There is a natural function $\cV:\cA\rightarrow \R$,
sending an element $a\in \cA$ to the sum of the volumes assigned to the $S_i$
by $a$.

Each point $a\in \cA$ determines an ``extended''
hyperbolic metric corresponding to a hyperideal simplex 
on each of the $S_i, i=1,\cdots,N$, but in general not an ``extended''
hyperbolic metric corresponding to a hyperideal polyhedron on $P$, because:
\begin{itemize}
\item It is in general not possible to glue the $S_i$ isometrically along
  their faces. 
\item Even if this gluing is possible, the 
  resulting hyperbolic metric might have
  non-trivial holonomy along the ``interior'' edges: the total angle around
  those edges could be different from $2\pi$.
\end{itemize}
Actually it is not difficult to check that those conditions are sufficient, so
that, if they are realized, then $a$ determines an identification of $P$ with
a hyperideal polyhedron, which is isometric on each of the $S_i, i=1,\cdots,
N$. There is an element $a_0\in \cA$ associated to the ``original'' situation
of $P$, since it determines an ``extended'' hyperbolic metric on each of the
$S_i, i=1,\cdots, N$. 

Let $E'$ be the set of edges of $P'$, and let $E_i$ be the set of ``interior''
edges of the cellular decomposition of $P'$ as the union of the $S_i, i=1,
\cdots, N$. Let $\cW$ be the vector space of maps $w:E'\cup E_i\mapsto
\R$. There is a natural affine map $\Phi:\cA\mapsto \cW$, sending an
element $a\in \cA$ -- i.e. an assignment of dihedral angles to each of the
edges of the $S_i, i=1,\cdots, N$ -- to a map $w(a):E'\cup E_i\mapsto \R$,
defined by:
\begin{itemize}
\item For each $e\in E'$, $w(a)(e)$ is $\pi$ minus the sum of the values of
  $a$ on the edge $e$ for all the $S_i$ which contain $e$.
\item For each $e_i\in E_i$, $w(a)(e_i)$ is 
  $2\pi$ minus the sum of the values of
  $a$ on the edge $e$ for all the $S_j$ which contain $e_i$.
\end{itemize}
For instance, $w(a_0)$ is zero on all interior edges, and is equal to the
(exterior) dihedral angle of $P'$ on all edges which are edges of $P'$.

Let $a\in \cA$, remark that the lengths of each edge is the same for
all the $S_i$ containing it if and only if $a$ is a critical point of $\cV$ on
$w^{-1}(w(a))$. This is a direct consequence of the Schl\"afli formula
(\ref{eq:schlafli}). Indeed, if the lengths assigned to each edge are the same
for all
simplices containing it, then it is quite clear that any first-order variation
of $a$ which does not change the sum of the angles at each edge leaves $\cV$
unchanged (at first order). 
Conversely, if two simplices, say $S_i$ and $S_j$, both contain an
edge $e\in E'\cup E_i$ and $a$ assigns different lengths to $e$ as an edge of
$S_i$ and of $S_j$, than there is a first-order deformation of $a$, tangent to
$w^{-1}(w(a))$, which increases the angle of $e$ in $S_i$ and decreases the
edge length of $e$ in $S_j$ by the same amount, so that $\cV$ varies.

Since $w$ is an affine function, $\cA$ is foliated by the level sets of $w$,
which are affine submanifolds. But $\cV$ is a sum of strictly concave
functions, so it is strictly concave. It follows that, in a neighborhood of
$a_0$ -- which is a local maximum of $\cV_{|w^{-1}(w(a_0))}$ -- each level set
of $w$ contains a unique local maximum of the restriction of $\cV$, so that
each small variation of the dihedral angles of $P'$ (resp. $P$)
is obtained uniquely by a small deformation of $P'$ (resp. $P$). So
Remark \ref{rk:concavity} shows that the volume of a strictly concave function
on hyperideal polyhedra in the neighborhood of $P$.
\end{proof}

This argument is very close to the one that is used in section 5 to prove the
infinitesimal rigidity of hyperideal circle patterns, and can be used as a
``toy model'' for it. The proof in section 5 is a little more complicated
since it has to take into account the possibility of ideal vertices.

\subsection{The 3-dimensional cone-manifold associated to a circle pattern.}

\paragraph{Euclidean circle patterns.}

Let $\Gamma$ be a graph embedded in a closed surface
$\Sigma$, as in the setting of Theorem
\ref{tm:euclid}. Following the outline given in section 2, we define a cell
complex $\cS_{\Sigma,\Gamma}$ associated to $\Gamma$ as the cone over
$\Gamma^*$. It has:
\begin{itemize}
\item One vertex for each vertex of $\Gamma^*$, and an additional one, $v_0$,
  common to all the 3-cells.
\item One 3-cell for each face $f$ of $\Gamma^*$. If $f$ has $p$ vertices, then
  the corresponding 3-cell of $\cS_{\Sigma,\Gamma}$ has $p+1$ vertices,
  corresponding to the vertices of $f$ and to the ``central'' vertex, $v_0$,
  which is common to all 3-cells.
\item One 2-face for each face $f$ of $\Gamma^*$ --- corresponding to the face
  of the 3-cell associated to $f$ which is opposite to $v_0$ --- 
  and one for each
  edge of $\Gamma^*$. We call ``horizontal'' the faces of the first kind, and
  ``vertical'' the others.
\item One edge for each edge of $\Gamma^*$, and one for each vertex of
  $\Gamma^*$. Again we call ``horizontal'' the edges of the first type, and
  ``vertical'' the others.
\end{itemize}
We call $S_1, \cdots, S_N$ the 3-cells of $\cS_{\Sigma,\Gamma}$. 

The definition of $\cS_{\Sigma,\Gamma}$ up to this point is entirely
topological. However it is also desirable to put on this cell complex a
hyperbolic metric, or more precisely the structure of a (non-compact)
hyperbolic cone-manifold, with cone singularities along the verticale edges,
and with a boundary which is composed of parts 
of the horizontal faces and of the
horizontal edges. To this end, it will be useful to consider two simple
propositions. The first is an elementary statement of plane geometry, which
will clear up the following considerations.

\begin{prop} \label{lm:4-cercles}
Let $C_1, C_2$ be two circles in the Euclidean (resp. hyperbolic) plane which
intersect in two points. Let $C'_1, C'_2$ be two non-intersecting circles
which are both orthogonal to both $C_1$ and $C_2$. Then the intersection
points of $C_1$ and $C_2$ lie on the line containing the centers of
$C'_1$ and $C'_2$.
\end{prop}

\begin{figure}[ht]
\centerline{\psfig{figure=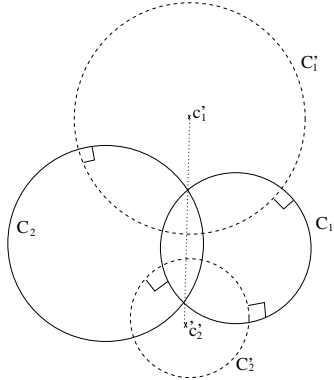,height=5cm}}
\caption{Proposition 4.5.}
\end{figure}

\begin{proof}
We consider first the Euclidean case. It is possible to identify isometrically
the Euclidean plane with a horosphere $H_0$ in $H^3$, ``centered'' at a point
$v_0\in \dr_\infty H^3$. Let $N$ be the unit normal vector field to $H_0$, in
the direction opposite to the direction of $v_0$. There is natural map
$G:H_0\rightarrow \dr_\infty H^3\setminus \{v_0\}$,
sometimes called the ``hyperbolic Gauss map'' for $H_0$, sending a point $x\in
H_0$ to the endpoint of the geodesic ray starting from $x$ in the direction of
$N$. This map is conformal: indeed, the pull-back to a surface in $H^3$ of the
conformal structure at infinity by the hyperbolic Gauss map is equal to the
conformal structure of $I+2\II+\III$, where $\II$ and $\III$ are the second
and third fundamental forms of the surface (see e.g. \cite{horo}). 
For a horosphere, $I=\II=\III$, so
that $I+2\II+\III=4I$, and $G$ is conformal. 

The images by $G$ of $C_1, C_2, C'_1$ and $C'_2$ are 
circles in $\dr_\infty H^3$
(for the natural $\C P^1$-structure on $\dr_\infty H^3$). Let $F_1$ and $F_2$
be the hyperbolic planes with boundary at infinity $G(C_1)$ and $G(C_2)$,
respectively, and also (with an abuse of notation) their extension as planes
in $\R^3$, using the projective model of $H^3$. Let also $F'_1, F'_2$ be the
hyperbolic planes with boundary at infinity $G(C'_1)$ and $G(C'_2)$, and let
$v'_1, 
v'_2$ be the dual points (again using the projective model of $H^3$). The
definition of the projective duality shows that $v'_1$ is contained in the
extension as an Euclidean line of the hyperbolic geodesic starting from the
center of $C'_1$ in the orthogonal direction (this can be checked in a special
case, with $C'_1$ in a simple position relative to $H_0$, and then extended to
the general case using the projective invariance of the duality).

Since $C_1$ is orthogonal to $C'_1$, $v'_1$ is contained in $F_1$. For the
same reason, $v'_1\subset F_2$, $v'_2\subset F_1$, and $v'_2\subset F_2$. So
both $v'_1$ and $v'_2$ are contained in the intersection of $F_1$ and $F_2$,
which is the line defined by the intersection of $G(C_1)$ and $G(C_2)$. The
result therefore follows by projecting all four points orthogonally to $H_0$.

The argument is the same in the hyperbolic case, using instead of $H_0$ a
totally geodesic plane $P_0\subset H^3$. The hyperbolic Gauss map considered
now goes from $P_0$ to a connected component of $\dr_\infty H^3\setminus
\dr_\infty P_0$, and it is still conformal since, for a hyperbolic plane,
$\II=\III =0$. The proof can therefore proceed as in the Euclidean setting.
\end{proof}

The next proposition describes a hyperbolic polyhedron associated to a simple
pattern of circles in the Euclidean plane. It will serve below as a ``building
block'' for the construction of a hyperbolic cone-manifold structure
associated to a hyperideal circle pattern. A similar statement
holds in the hyperbolic case, but it is kept separate so as to avoid
complicated notations. 

Consider the following situation. 
$C_0$ is an oriented circle in the Euclidean plane, and $C'_1, \cdots,
C'_N$ are oriented circles bounding disjoint disks, orthogonal to $C_0$,
appearing in cyclic order along it, with $C'_i$ having center $c'_i$ and
radius $r_i$, $i=1, \cdots, N$. Let $d_i$ be the Euclidean 
distance between $c'_i$ and $c'_{i+1}$. The Euclidean plane is identified, as
above, with a horosphere $H_0$ ``centered'' at a point $v_0\in \dr_\infty
H^3$, and we consider the hyperbolic Gauss map
$G:H_0\rightarrow \dr_\infty H^3$. Let $v'_1, \cdots, v'_N$ be the points dual
to the hyperbolic planes with boundary at infinity $G(C'_1),\cdots,G(C'_N)$. 

It is then possible to consider the polyhedron $P$, first in $\R^3$ using the
projective model of $H^3$, with vertices $v_0, v'_1, \cdots, v'_N$. 

\begin{figure}[ht]
\centerline{\psfig{figure=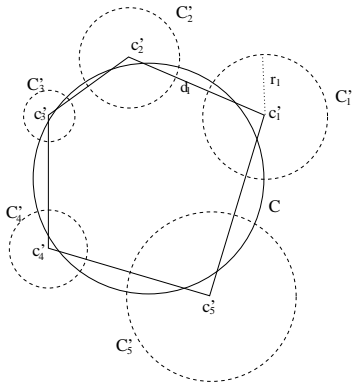,height=5cm}}
\caption{Proposition 4.6.}
\end{figure}

\begin{prop} \label{pr:pyramide}
$P$ is a hyperideal polyhedron, with ideal vertex $v_0$ and strictly
hyperideal vertices at $v'_1, \cdots, v'_N$, it is combinatorially a pyramid
with vertex $v_0$. Moreover, defining the edge
lengths involving $v_0$ with respect to the horosphere $H_0$:
\begin{enumerate}
\item The length of the edge $[v_0, v'_i]$ is $-\log(r_i)$.
\item The length of the edge $[v'_i, v'_{i+1}]$ depends only in $d_i, r_i$
  and $r_{i+1}$.
\item The (exterior) dihedral angle at the edge $[v'_i,v'_{i+1}]$ is equal to
  the angle between $[c'_i,c'_{i+1}]$ and $C_0$, measured in the interior of
  the intersection of the interior of $C_0$ with the interior of the polygon
  $p$ with vertices the $c'_i, i=1,\cdots, N$.
\item The (exterior) dihedral angle at the edge $[v_0,v'_i]$ is equal to the
  exterior angle of $p$ at $c'_i$.
\end{enumerate}
\end{prop}

Point (2) could be made more explicit at the cost of a fairly classical
computation, but a precise expression will not be necessary here. 

\begin{proof}
It is simplest to consider the first point in the Poincar\'e half-space model,
taking as $H_0$ the plane $\{ z=1\}$. Then the Gauss map is simply the
vertical projection on the plane $\{ z=0\}$, which corresponds to $\dr_\infty
H^3\setminus \{ v_0\}$. So the hyperbolic plane with boundary at infinity
$G(C'_i)$ corresponds to a half-sphere of radius $r_i$, and its hyperbolic
distance to $H_0$ is obtained by integrating $dz/z$ from $r_i$ to $1$, so it
is equal to $-\log(r_i)$. This, along with the definition of the length of the
edges of hyperideal polyhedra, proves the first point.

The second point can be proved in the same manner; still in the Poincar\'e
half-space model, the hyperbolic planes with boundary at infinity $G(C'_i)$
and $G(C'_{i+1})$ correspond to half-spheres of radii $r_i$ and $r_{i+1}$,
respectively, and with centers at Euclidean distance $d_i$ on the plane $\{
z=0\}$. Computing the distance between them is thus a simple exercice (which
we leave to the reader) and the result is a function only of $r_i, r_{i+1}$
and $d_i$. 

The third and last points can be proved using the same model. In the Poincar\'e
half-space model, the faces of $P$ correspond to the half-sphere with boundary
the circle $G(C_0)$, which is contained in the plane $\{ z=0\}$, and the
vertical strips intersecting the plane $\{ z=0\}$ at the edges of the polygon
$p$ with vertices the images by $G$ of the $c'_i$. By the conformality of the
Poincar\'e half-space model, the hyperbolic angles between those faces are the
same as the Euclidean angles, and the result follows directly.
\end{proof}

\begin{remark}
Suppose that some of the $C'_i$ are replaced by points, so that some of the
vertices of $P$ other than $v_0$ are ideal vertices. Choose for each of those
ideal vertices the horosphere which is tangent to $H_0$. Then points (2), 
(3) and (4) of the previous proposition still hold. Point (1) holds as well
when $v'_i$ is strictly hyperideal, while the length of the edge going from
$v_0$ to $v'_i$ is zero when $v'_i$ is ideal.    
\end{remark}

The proof of this remark uses exactly the same arguments as the proof of the
previous proposition. Proposition \ref{pr:pyramide} and this remark, along with
Proposition \ref{pr:gluing}, make it possible to define a hyperbolic metric
associated to a hyperideal circle pattern.

\begin{df} \label{df:metrique-tore}
Let $C$ be a hyperideal circle pattern embedded in $\Sigma$, with incidence
graph $\Gamma$. We call $H(C)$ the non-complete metric defined on the
complement of the vertical edges of $\cS$,
defined from $C$ by gluing one hyperideal
hyperbolic pyramid for each face $f$ of $\Gamma^*$, as described in Proposition
\ref{pr:pyramide}, taking as $C_0$ the circle corresponding to $f$ (or more
precisely to the dual vertex of $\Gamma$) and as $C'_1, \cdots, C'_N$ the dual
circles corresponding to the vertices of $f$ (or more precisely to the dual
faces of $\Gamma$). Those pyramids are glued isometrically along their
``vertical'' faces, which correspond to the edges of $\Gamma^*$.
\end{df}

Note that it is indeed possible to glue the pyramids along their vertical
faces since, by Proposition \ref{pr:pyramide}, corresponding faces in pyramids
corresponding to adjacent faces of $\Gamma^*$ are triangles with the same edge
lengths, so that, by Proposition \ref{pr:gluing}, they are isometric.

\begin{prop} \label{pr:cone-metric-eucl}
The completion of $H(C)$ is a hyperbolic metric on $\cS$ with conical
singularities along the ``vertical edges''. For each such edge $e$, 
corresponding to a
vertex $v$ of $\Gamma^*$ (that is, to a face $v^*$ of $\Gamma$) the total
angle around $e$ is equal to $2\pi$ minus the singular curvature $\kappa(v^*)$
of $C$ at $v^*$. 
\end{prop}

\begin{proof}
The completion of $H(C)$ is obtained by adding a line for each vertical edge
of $\cS$. The holonomy of $H(C)$ around each such edge is the composition of a
translation along the edge and a rotation with axis equal to the edge. 
Proving that the completion is a metric with conical singularities is
equivalent to proving that the translation component vanishes. 

On each 3-cell $S_i$ of $\cS$, one can define a ``Busemann function'' $B_{S_i}$
corresponding to the point at infinity $v_0$, as the ``oriented distance'' to
the horosphere $H_0$ -- that is,  $B_{S_i}$ is equal to the distance to $H_0$
for points which are on the side of $H_0$ opposite to $x_0$, and to minus the
distance to $H_0$ for points which are on the same side of $H_0$ as
$x_0$. Proving that the translation component of the holonomy along a
``vertical'' edge $e$ vanishes is
clearly equivalent to proving that the holonomy acts trivially on $B_{S_i}$,
for each 3-cell $S_i$ adjacent to $e$. But when one considers the 3-cells
adjacent to $e$, in cyclic order, the point of $e$ where the functions
$B_{S_i}$ vanish remains the same -- it is the intersection of $e$ with
$H_0$. So, after completing a turn around $e$, the function $B_{S_i}$ remains
the same, so that the translation component of the holonomy remains the same,
and the completion of $H(C)$ is indeed a hyperbolic metric with conical
singularities. 

The total angle around a ``vertical'' edge $e$ of $\cS$ is obtained by summing
the {\it interior} 
dihedral angles at $e$ of the 3-faces adjacent to it. Since those angles
are equal to the angles at the vertex of $\Gamma^*$ corresponding to $e$ of
the polygons corresponding to the faces of $\Gamma^*$, the interior dihedral
angles at $e$ of the 3-cells adjacent to it sum to the total angle around the
singular point at the center of the circle corresponding to $v$, i.e. to
$2\pi-\kappa(v^*)$. 
\end{proof}

The metric which is constructed in this manner is not complete; in addition to
the conical singularities, it has a boundary, which is ``polyhedral'',
i.e. each of its points has a neighborhood in which the boundary looks either
like a plane or like the neighborhood of an edge in a hyperbolic polyhedron.

\paragraph{A property of those hyperbolic metrics.}

The hyperbolic metrics constructed on $\cS$ above have one important property 
that is not quite apparent yet: it is possible to choose for each cell in $\cS$
a horosphere, centered at the ``central'' vertex, so that the horospheres
``match'' along the vertical faces, as in the next definition.

\begin{df} \label{df:consistent}
Let $h$ be a hyperbolic metric on $\cS$, with conical singularities along the 
vertical edges, which is obtained by gluing the hyperbolic metrics on the cells
corresponding to the identification of each cell with a hyperideal polyhedron.
We say that $h$ is {\bf consistent} if it is possible to choose, for each cell 
$C_i$ of $\cS$, a horosphere $H_i$ in $C_i$ centered at the ``central'' vertex,
so that, if $C_i$ and $C_j$ share a 2-dimensional face $f$, then 
$H_i\cap f=H_j\cap f$.
\end{df}
 
The way the metric $H(C)$ was constructed implies that:

\begin{remark}
$H(C)$ is consistent.
\end{remark}

\begin{proof}
This follows directly from the construction of $H(C)$, which was made in each
cell $C_i$ of $\cS_{\Sigma, \Gamma}$ with respect to a horosphere $H_0$, 
which we can rename as $H_i$. Then the $H_i$  have precisely the requested
property. 
\end{proof}

Conversely, any hyperbolic metric which is consistent is related to a
hyperideal circle pattern. 

\begin{remark} \label{rk:converse}
Let $h$ be a hyperbolic metric on $\cS_{\Sigma,\Gamma}$ with conical
singularities on the vertical edges. Suppose that the restriction of $h$ to
each 3-dimensional cell of $\cS_{\Sigma,\Gamma}$ is the metric induced on a
hyperideal polyhedron, with an ideal vertex at $v_0$, and that $h$ 
is consistent. Suppose moreover that the boundary of $\cS_{\Sigma,\Gamma}$ for
$h$ is locally convex, i.e. the exterior dihedral angle is positive at each
``horizontal'' edge of $\cS_{\Sigma,\Gamma}$. For each cell $C_i$
of $\cS_{\Sigma,\Gamma}$, let $H_i$ be the part of 
horosphere appearing in Definition
\ref{df:consistent}. Then the $H_i$ are isometric to Euclidean polygons which
can be glued along their edges to obtain a Euclidean metric with conical
singularities, which has a hyperideal circle pattern with
incidence graph $\Gamma$ and intersection angles equal to the exterior
dihedral angles between the ``horizontal'' faces of $(\cS_{\Sigma,\Gamma}, h)$.
\end{remark}

\begin{proof}
The fact that the $H_i$ are isometric to Euclidean polygons follows from the
elementary properties of horospheres. They can be glued isometrically along
their edges because, if two cells $C_i$ and $C_j$ share a 2-dimensional face
$f$, then $H_i\cap f=H_j\cap f$, so that those edges of $H_i$ and $H_j$ can be
glued isometrically, to obtain an Euclidean metric with conical singularities.

Each of the $H_i$ is the orthogonal projection of the corresponding cell
$C_i$, so that its vertices lie on an Euclidean circle $c_i$ 
-- the projection of
the boundary at infinity of the ``horizontal'' face of $C_i$. The convexity
condition in the statement means precisely that the interior of 
$c_i$ intersects no vertical edge of $\cS_{\Sigma,\Gamma}$, so that the $c_i$
indeed constitute a hyperideal circle pattern. It is then clear that the
intersection angle between $c_i$ and $c_j$ is equal to the exterior dihedral
angle between the ``horizontal'' faces of $C_i$ and $C_j$.
\end{proof}

\paragraph{Hyperbolic circle patterns.}

It is also possible to use the same construction for hyperideal circle
patterns on a closed surface endowed with a singular hyperbolic
metric. The main difference is that the ``central'' vertex $v_0$ is now
strictly hyperideal rather than ideal. 
So the situation we consider is now as follows. 
$C_0$ is an oriented circle in the hyperbolic plane, and $C'_1, \cdots,
C'_N$ are oriented circles bounding disjoint disks, orthogonal to $C_0$,
appearing in cyclic order along it, with $C'_i$ having center $c'_i$ and
radius $r_i$, $i=1, \cdots, N$. $d_i$ is the hyperbolic 
distance between $c'_i$ and $c'_{i+1}$. The hyperbolic plane is identified
with a totally geodesic plane $P_0$ with dual point $v_0$, 
and we consider the hyperbolic Gauss map
$G:H_0\rightarrow \dr_\infty H^3$. Let $v'_1, \cdots, v'_N$ be the points dual
to the hyperbolic planes with boundary at infinity $G(C'_1),\cdots,G(C'_N)$,
and let $P$ be the polyhedron with vertices $v_0, v'_1, \cdots, v'_N$.
The analog of Proposition \ref{pr:pyramide} is as follows.

\begin{prop} \label{pr:pyramide-hyp}
$P$ is a hyperideal polyhedron, with strictly
hyperideal vertices at $v_0, v'_1, \cdots, v'_N$, 
it is combinatorially a pyramid
with vertex $v_0$. Moreover:
\begin{enumerate}
\item the length of the edge $[v_0, v'_i]$ depends only on $r_i$,
\item the length of the edge $[v'_i, v'_{i+1}]$ depends only in $d_i, r_i$
  and $r_{i+1}$,
\item the (exterior) dihedral angle at the edge $[v'_i,v'_{i+1}]$ is equal to
  the angle between $[c'_i,c'_{i+1}]$ and $C_0$, measured in the interior of
  the intersection of the interior of $C_0$ with the interior of the polygon
  $p$ with vertices the $c'_i, i=1,\cdots, N$,
\item the (exterior) dihedral angle at the edge $[v_0,v'_i]$ is equal to the
  exterior angle of $p$ at $c'_i$.
\end{enumerate}
\end{prop}

Points (1) and (2) could be made more explicit at the cost of a fairly
classical computation (which is even easy for point (1)), 
but a precise expression will not be necessary here. 

\begin{proof}
For the first point note that, by definition of the edge lengths of a
hyperideal polyhedron, the length of the edge between $v_0$ and $v'_i$ is
equal to the distance between $P_0$ and the hyperbolic plane with boundary at
infinity $G(C'_i)$. Clearly this only depends on the radius of the circle in
$P_0$ which is the orthogonal projection of $C'_i$. 

The second point can be checked in the same manner, now the length is equal to
the hyperbolic distance between the hyperbolic planes with boundaries at
infinity $G(C'_i)$ and $G(C'_{i+1})$, and again it only depends on the radii
of $C'_i$ and $C'_{i+1}$ and on the distances between their centers.

Finally the third and last points can be proved as in the Euclidean setting.
\end{proof}

It follows, as in the Euclidean case, that it is possible to glue
isometrically hyperideal pyramids, corresponding to the faces of $\Gamma^*$,
to obtain a hyperbolic cone-manifold with polyhedral boundary canonically
associated to a hyperideal circle pattern on a hyperbolic surface.

Conversely, it is possible to reconstruct a hyperideal circle pattern
on a hyperbolic surface with conical singularities (at the centers of 
the dual circles) from such a gluing of hyperideal pyramids. 
This follows from the same argument as the one used in 
Remark \ref{rk:converse}, with the difference that the projection
is now on the planes dual to $v_0$ in each 3-cell of $\cS(\Sigma,\Gamma)$.


\section{Infinitesimal rigidity}

\subsection{Introduction and outline}

This section contains the proof, in the different situations which are
considered, of the infinitesimal rigidity of hyperideal circle patterns: any
first-order deformation of such a pattern -- including a deformation of the
underlying singular Euclidean or hyperbolic metric -- induces a first-order
deformation of either the intersection angles between the circles or the
singular curvatures at the cone singularities. 

The main idea used here is not original, it originates in papers of Rivin
\cite{Ri2} and Br\"agger \cite{bragger} on ideal polyhedra, and was
also related to or used in
\cite{CdeV,bobenko-springborn,rivin-combi,leibon1,leibon2,ideal}. Its
extension to hyperideal polyhedra was already used to some extend in 
\cite{hphm,rcnp}. 

The next subsection describes some relevant spaces of assignments of angles to
the cells of $\cS$, which are then used in the infinitesimal rigidity proof. 
The argument is slightly simpler when one considers hyperbolic circle
patterns, so this case is considered first. 

\subsection{Spaces of angle assignments}

To each graph $\Gamma$ on a closed surface, we have associated the
cellular complex $\cS(\Gamma)$. 
We can then associate to $\cS(\Gamma)$ a polytope, which is the space of
possible assignments of dihedral angles to the 3-dimensional simplices of
$\cS(\Gamma)$ so that, for each cell, the angles are compatible with Theorem
\ref{tm:bb}. The definition takes into account which vertices of $\Gamma^*$
are intended to be ideal vertices. It is first written for an
Euclidean metric with conical singularities.

\begin{df} \label{df:sigma}
Let $\Gamma$ be the 1-skeleton of a cellular decomposition of a closed surface
$\Sigma$,
and let $V\subset \Gamma_2$ be a subset of the set of faces of $\Gamma$. 
We call $\cA_E(\Gamma)$ the space of assignments, to each of the cells
of $\cS(\Gamma)$ 
associated to the vertices of $\Gamma$, of a set of dihedral angles,
satisfying the conditions given in Theorem \ref{tm:bb}, so that, for each
cell:  
\begin{itemize}
\item the sum of the dihedral angles of the vertical edges is equal to
  $2\pi$, 
\item for each $v$ vertex other than the central vertex, the sum of the angles
  assigned to the adjacent edges is equal to $2\pi$ if $v\in V$, and strictly
  larger than $2\pi$ otherwise.
\end{itemize}
\end{df}

A natural complement is the definition of a space of assignments of angles to
the edges of $\cS(\Gamma)$ --- of course this space has much lower dimension
than $\cA_E(\Gamma)$, since each edge of $\cS(\Gamma)$ can be an edge of
several simplices of that same simplicial complex.

\begin{df}  \label{df:52}
Let $\Gamma$ be the 1-skeleton of a cellular decomposition of $\Sigma$. 
\begin{itemize}
\item We call $\cD'(\Gamma)$ the space of functions from the set of edges of 
$\cS(\Gamma)$ to $\R$.
\item We call $\sigma$ the linear map from $\cA_E(\Gamma)$ 
to $\cD'(\Gamma)$, such that, for each $\Theta\in \cA_E(\Gamma)$ and each
edge $e$ of $\cS(\Gamma)$, the value of $\sigma(\Theta)$ on $e$ is equal
to the sum of the angles at $e$ of the cells of $\cS(\Gamma)$ containing it.
\end{itemize}
\end{df}

Note that $\cA_E(\Gamma)$ is defined by linear
equalities and inequalities, so that it is the intersections of a finite
number of half-spaces in $\R^p$, for some $p$. Moreover $\cA_E(\Gamma)$ is
relatively compact, since the angles are in $(0, \pi)$, so it is a
polytope. It follows that its image by $\sigma$ is also a polytope in
$\cD'(\Gamma)$.  
Definition \ref{df:sigma} has a variant for surfaces with
singular hyperbolic metrics:

\begin{df}
Let $\Sigma$ be a closed surface, let 
$\Gamma$ be the 1-skeleton of a cellular decomposition of $\Sigma$,
and let $V\subset \Gamma_2$ be a subset of the set of faces of $\Gamma$. 
We call $\cA_H(\Sigma,\Gamma)$ the space of assignments, to each of the cells
of $\cS(\Gamma)$
associated to the vertices of $\Gamma$, of a set of dihedral angles,
satisfying the conditions given in Theorem \ref{tm:bb}, so that, for 
each cell:
\begin{itemize}
\item the sum of the dihedral angles of the vertical edges
is strictly larger than $2\pi$,
\item for each vertex $v$ other than the central vertex, the sum of the angles
  assigned to the adjacent edges is equal to $2\pi$ if $v\in V$, and strictly
  larger otherwise.
\end{itemize}
\end{df}

\subsection{Gluing conditions}

An angle assignment on the simplices of $\cS(\Gamma)$, for a graph
embedded in a surface, does not always define a hyperbolic metric on
$\cS(\Gamma)$. There are basically two reasons for this; the first is that it
might not be possible to glue isometrically the faces of the simplices, the
second is that, if such a gluing is possible, some nasty singularities might
appear along the edges of $\cS(\Gamma)$. 

\paragraph{The possibility of gluing.}

The next definition is quite natural, it should be compared to the definition
of term ``consistent'' above.

\begin{df}
Let $\Gamma$ be the 1-skeleton of a cellular decomposition of a closed
surface, and let
$\Theta\in \cA_E(\Gamma)$. $\Theta$ is {\bf compatible} 
if the hyperbolic metrics on
the simplices of $\cS(\Gamma)$ can be glued along the vertical faces to yield a
(non-complete) hyperbolic metric on the complement of the edges of
$\cS(\Gamma)$. 
\end{df}

The same definition will be used in the contex of circle patterns on closed
hyperbolic surfaces, as well as for Euclidean and hyperbolic circle patterns
on the disk; we do not repeat the definition. 

\paragraph{More on the possible singularities.}

A compatible angle assignment $\Theta$ does not, in general, define a
hyperbolic metric which can be extended to the edges of $\cS(\Gamma)$, since
the holonomy around those edges could be non-trivial. More precisely, each
edge $e$ of $\cS(\Gamma)$ corresponds to a 
segment of geodesic for the hyperbolic
metrics on each of the simplices which contain it, and the holonomy of a curve
going around
this edge can be decomposed as the sum of a translation along $e$ and a
rotation around $e$. The translation component of the holonomy around $e$
can be obtained as follows. Let $S_1, \cdots,
S_p$ be the 3-cells containing $e$, in the cyclic orders in which they
stand around $e$ (with $S_{p+1}:=S_1$). 
Let $j\in \{ 1, \cdots, p\}$, so that $S_j$ and $S_{j+1}$  
share a face $f_j$. Let $H_1$ be a totally geodesic plane in $S_1$ which
is orthogonal to $e$. Then $H_1$ can be extended uniquely across $f_1$ to a
totally geodesic plane $H_2$ in $S_2$, orthogonal to $e$. Repeating this
procedure, we obtain totally geodesic plane $H_3, \cdots, H_p$, all orthogonal
to $e$. The translation component along $e$ of the holonomy around $e$ is then
the oriented distance, along $e$, between $e\cap H_1$ and $e\cap H_p$. The
same holds with the totally geodesic planes replaced by horospheres. 

In addition, the holonomy on elements of the fundamental
group of the surface can also be ``complicated'', as should be apparent from
the definition of a ``consistent'' metric in Definition \ref{df:consistent}.

The case of circle patterns on hyperbolic surfaces is much simpler,
because, as will be apparent below, the translation component of the holonomy
is always zero, so that the non-complete metric obtained by gluing the
hyperideal simplices (when it is possible) can always be completed to a
hyperbolic metric with cone singularities along lines (and with polyhedral
boundary). 

\subsection{Hyperbolic circle patterns on closed surfaces}

The description just made in the hyperbolic case can be translated as the
following lemma.

\begin{lemma} \label{lm:compatib-h}
Let $\Sigma$ be a closed surface, and let $\Gamma$ be the
1-skeleton of a cellular decomposition of $\Sigma$, with trivalent vertices. 
Let $\Theta\in\cA_H(\Sigma,\Gamma)$. If $\Theta$ is compatible, then it
defines a hyperbolic metric on $\cS(\Gamma)$ with conical singularities along
the edges.
\end{lemma}

\begin{proof}
Let $S_1, \cdots, S_N$ be the simplices in $\cS(\Gamma)$. 
For each $i=1,\cdots, N$,
there exists a unique hyperbolic plane $H_i$ in $S_i$ which is dual to the
``central'' vertex of $\cS(\Gamma)$. It is characterized by the fact that 
it intersects orthogonally all the edges adjacent to the ``central'' vertex of
$S_i$. 

Suppose that $S_i$ and $S_j$ are two simplices in $\cS(\Gamma)$ which share a
2-dimensional face, $f$. The intersections of $H_i$ and of $H_j$ with $f$ are
two geodesic segments which both intersect orthogonally the two edges of $f$
adjacent to the ``central'' vertex. It follows that $H_i\cap f=H_j\cap f$. 

Let $e$ be a vertical edge of $\cS(\Gamma)$. Then the description
given above of the translation component of the holonomy around $e$
can be applied with each $H_i$ equal to the totally geodesic plane, in $S_i$,
dual to the ``central'' vertex $v_0$; indeed those planes are orthogonal to
$e$, and $H_i\cap e=H_{i+1}\cap e$, as needed. It follows immediately that the
translation component of the holonomy around $e$ vanishes, so that the
hyperbolic metric defined by $\Theta$ can be completed to a hyperbolic metric
with conical singularities along the vertical edges of $\cS(\Gamma)$.
\end{proof}

The conditions under which the simplices of $\cS(\Gamma)$ can be glued -- 
according to Proposition \ref{pr:gluing}, that the length of the edges of 
each ``vertical'' face is the same for both simplices containing it -- 
can be expressed quite simply in the hyperbolic case. Indeed, each 3-cell of
$\cS(\Gamma)$ contains a hyperbolic plane which is dual to $v_0$, so there
is a canonical choice of the horospheres centered at the ideal vertices --
they  can be chosen as the one which are tangent to this dual plane. 

There is a natural function defined over $\cA_H(\Gamma)$: the sum of 
the hyperbolic volumes of the simplices in $\cS(\Gamma)$, for the 
hyperbolic metrics determined by the dihedral angles which are 
assigned to each of them. 
Note that those hyperbolic metrics are the
metrics on hyperideal simplices, so that the volume which is used
has to be the volume of the corresponding truncated simplices.
We call $\cV$ this total volume.

\begin{lemma} \label{lm:critique-h}
Let $\Theta\in \cA_H(\Sigma, \Gamma)$ be an assignement of angles to 
the simplices in $\cS(\Gamma)$, and let $\alpha:=\sigma(\Theta)$. 
Then $\Theta$ is compatible if and only if 
$\Theta$ is a critical point of the restriction of $\cV$ to 
$\sigma^{-1}(\alpha)\subset \cA_H(\Sigma,\Gamma)$. 
\end{lemma}

\begin{proof}
Suppose first that $\Theta$ is a compatible, and choose for each ideal
vertex of each 3-cell of $\cS(\Gamma)$ a horosphere as suggested above, 
that is, tangent 
to the hyperbolic plane which is dual to the ``central'' vertex $v_0$. 

Consider a ``vertical'' edge $e$ of $\cS(\Gamma)$, such that the endpoint
of $e$ other than $v_0$ is an ideal vertex $v$ of $\cS(\Gamma)$.
Then the hyperbolic
planes dual to $v_0$ in each of the 3-cells containing $e$ intersect
$e$ at the same point, say $p$, where they are orthogonal to $e$. The 
horospheres centered at $v$ in each of those 3-cells are also 
orthogonal to $e$; since they are tangent to planes dual to $v_0$,
they have to intersect $e$ at $p$.
It follows that the length of $e$
is $0$ in each of the 3-cells that contain it (relative to the choice of 
the horospheres tangent to the planes dual to $v_0$). According to 
Proposition \ref{pr:gluing}, the length of each ``horizontal'' edges of 
$\cS(\Gamma)$, for the hyperbolic metrics on each of the two 3-cells
containing it (and for this choice of horospheres) have to be equal. 

Let $\Theta'$ be a first-order deformation of $\Theta$ in 
$\sigma^{-1}(\alpha)$. In other terms, to each 3-cell $S$ of $\cS(\Gamma)$
and each edge $e$ of $S$ is associated a number, $\Theta'(S,e)$, 
corresponding to the first-order variation of the dihedral angle of $S$
at $e$, and the $\Theta'(S,e)$ satisfy some linear constraints corresponding
to the ideal vertices (the vertices in $V)$. 
The fact that $\Theta'$ is tangent to $\sigma^{-1}(\alpha)$
then means that the sum of the numbers attached to an edge, for 
each of the 3-cells containing it, vanishes. 

It follows from the Schl\"afli formula (\ref{eq:schlafli}):
\begin{equation} \label{eq:sums} 
d\cV(\Theta') = \sum_{S\subset \cS(\Gamma)}\sum_{e\subset S} l_S(e) 
\Theta'(S,e)~, 
\end{equation}
where the first sum is over the 3-cells in $\cS(\Gamma)$ and the 
second over the edges in each 3-cells, and $l_S(e)$ is the 
length of the edge $e$ for the metric in the cell $S$. But we have 
seen that $l_S(e)$ does not depend on $S$ and can be written as $l(e)$,
so that the 
sum above can be written as:
$$ d\cV(\Theta') = \sum_{e} l(e) \sum_{S\supset e}  
\Theta'(S,e)~, $$
which makes it apparent that $d\cV(\Theta')=0$.

Conversely, suppose that $\Theta$ is not compatible, which means that
it is not possible to glue isometrically the simplices in $\cS(\Gamma)$ 
along their common faces. If there is an edge $e$ of $\cS(\Gamma)$ between two
strictly hyperideal vertices, which has different lengths in two simplices
$S$ and $S'$
which contain it, then we can consider the first-order variation of the
angles which increases at unit speed the angle of $S$ at $e$ and decreases
at unit speed the angle of $S'$ at $e$, while keeping all other angles 
fixed. The Schl\"afli formula shows that, under this deformation, the
first-order variation of $\cV$ is non-zero. So we can now suppose that 
the edges between strictly hyperideal vertices are the same in all the 
simplices containing them.

Choose horospheres centered at the ideal vertices 
of $\cS(\Gamma)$ as above, i.e. tangent to the hyperbolic planes dual to 
$v_0$. Then the same argument as above clearly shows that the lengths of
the ``vertical'' edges of $\cS(\Gamma)$ are zero in each of the cells that
contain them, so, by Proposition \ref{pr:gluing}, there is at least 
one ``horizontal'' edges of $\cS(\Gamma)$, $e_0$, which has different 
lengths in the two cells containing it, say $S_1$ and $S_2$. 

Let $v_1$ and $v_2$ be the endpoints of $e_0$. 
Consider the first-order variation $\Theta'$ of $\Theta$ which is 
defined as follows:
\begin{itemize}
\item the angle of $S_1$ at $e$ has derivative $1$, and the angle of $S_2$ 
at $e$ has derivative $-1$,
\item the angle of $S_1$ on the ``vertical'' edge of $\cS(\Gamma)$ ending
at $v_1$ has derivative $-1$, as the angle of $S_1$ on the ``vertical'' edge
of $\cS(\Gamma)$ ending at $v_2$,
\item the angle of $S_2$ on the ``vertical'' edge of $\cS(\Gamma)$ ending
at $v_1$ has derivative $1$, as the angle of $S_2$ on the ``vertical'' edge
of $\cS(\Gamma)$ ending at $v_2$,
\item the other angles remain constant.
\end{itemize}

Then $\Theta'$ is a tangent vector to $\cA_H(\Sigma,\Gamma)$, because the
only linear constraints on the first-order variations of the angles
of the simplices is that the sum of the angles in each cell of the edges 
adjacent to an ideal vertex remain constant, and it happens here for 
both $v_1$ and $v_2$ (which might or might not be ideal). 

Moreover, $\Theta'$ is tangent to $\sigma^{-1}(\alpha)$, because it
is easy to check that the sum of the angles remain fixed at $e_0$
and at the vertical edges ending at $v_1$ and at $v_2$. However the 
Schl\"afli formula shows that $d\cV(\Theta')\neq 0$, with the only
non-zero term in (\ref{eq:sums}) coming from $e_0$. This finishes the proof. 
\end{proof}

The key infinitesimal rigidity statement that we need follows from 
the previous lemma and from the strict concavity of $\cV$, considered
as a function of the dihedral angles of the simplices in 
$\cS(\Sigma,\Gamma)$. We say here that a dual circle is ``ideal'' if it is
reduced to a point.

\begin{lemma} \label{lm:rigidity-h}
Let $\Sigma$ be a closed surface, and let $C$ be a 
hyperideal circle pattern on $\Sigma$, with respect to a hyperbolic metric
with conical singularities at the center of the dual circles. The 
first-order deformations of $C$, among hyperideal circle patterns with
the same incidence graph and the same set of ideal dual circles, are
parametrized by the first-order variations of the intersection angles
between the circles and of the singular curvatures at the
centers of the dual circle $c'$, under the condition that, for each ideal dual
circle, the sum of the intersection angles between adjacent principal circles
varies as the singular curvature at $c'$.
\end{lemma}

\begin{proof}
The constructions made in section 4, and in particular Proposition 
\ref{pr:pyramide-hyp}, show that the statement of the lemma is 
equivalent to another statement on the angle assignments to the 
simplices of $\cS(\Sigma, \Gamma)$, where $\Gamma$ is the incidence
graph of $C$. Namely that, given $\Theta\in \cA_H(\Sigma,\Theta)$
which is compatible, the first-order deformations of $\Theta$
among compatible angle assignments are uniquely determined by the
corresponding first-order variations of their images in $\cD'(\Sigma,\Gamma)$ 
by $d\sigma$. The fact that each ideal dual circle $c'$ remains ideal is a
consequence of the fact that the sum of the intersection angles between the
adjacent principal circles varies as the singular curvature at $c'$.

Since $\sigma$ is a linear map, the inverse images by $\sigma$ of the
elements of $\cD'(\Sigma,\Gamma)$ define a foliation of $\cA_H(\Sigma,\Gamma)$ 
by affine submanifolds. Since $\Theta$ is compatible, it is a critical
point of $\cV$ restricted to the inverse image by $\sigma$ of 
$\sigma(\Theta)$ by Lemma \ref{lm:critique-h}. Since $\cV$ is strictly 
concave over $\cA_H(\Sigma,\Gamma)$ -- as a sum of the volumes of
the simplices, which are strictly concave functions by Lemma \ref{lm:concave} 
-- $\Theta$ is a local maximum of $\cV$, and each leaf close to $\Theta$,
i.e. the inverse image by $\sigma$ of each element $\beta$ 
of $\cD'(\Sigma,\Gamma)$ 
close to $\sigma(\Theta)$, contains a unique maximum of the restriction of
$\cV$. Moreover, the inverse function theorem, along with the strict concavity
of $\cV$, shows that this local 
maximum varies smoothly with $\beta$. The result follows.
\end{proof}

\subsection{Euclidean circle patterns and the holonomy}

We now turn to hyperideal circle patterns on a surface with 
a singular Euclidean metric. As mentioned above, the argument is 
basically similar to what has just been said on the case of
hyperbolic surfaces, but with some additional arguments 
which are necessary to understand the translation component of 
the holonomy around ``vertical'' edges of $\cS(\Gamma)$, as
well as when the hyperbolic is consistent (cf Definition 
\ref{df:consistent}). It is in particular helpful to consider a special
kind of first-order variation of the dihedral angles of the
simplices, as explained in the next paragraph. 

\paragraph{Some convenient deformations.}

\begin{figure}[ht] 
\centerline{\psfig{figure=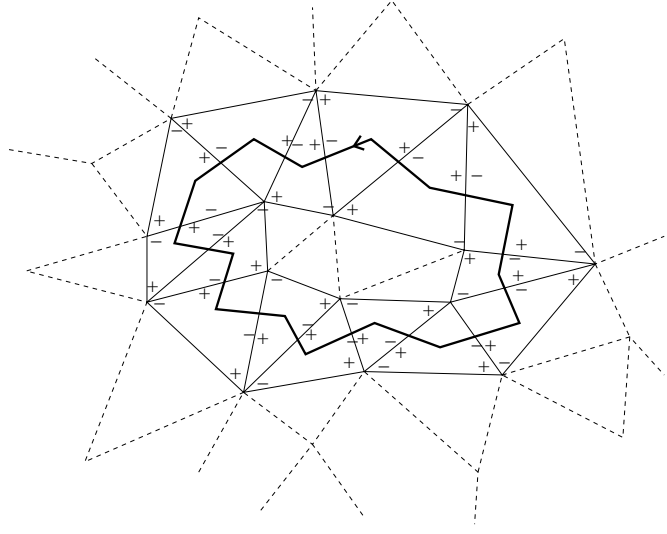,height=7cm}}
\caption{A useful deformation associated to a closed path in $\Gamma$.}
\label{fg:defo}\end{figure}

Let $\gamma$ be an oriented simple closed path in $\Gamma$. Consider
an angle assignment $\Theta \in \cA_E(\Sigma,\Gamma)$. We associate to 
$\gamma$ a first-order deformation $\Theta'_\gamma$ of $\Theta$ (i.e. a 
vector tangent to $\cA_E(\Sigma,\Gamma)$ at $\Theta$) as explained in
Figure \ref{fg:defo}; the closed curve
$\gamma$ appears as the thicker polygonal curve.
Each triangle corresponds to a simplex in
$\cS(\Sigma, \Gamma)$, each segment to a ``horizontal'' edge, and each
vertex to a ``vertical'' edge of $\cS(\Sigma, \Gamma)$; the signs ``$+$'' 
are attached to edges with an angle, in the simplex on which the sign stands,
which has differential $+1$, while the signs ``$-$'' are attached to 
edges with an angle with differential $-1$. The edges which have no
sign attached have an angle which remains constant. 

A more explicit description of $\Theta'_\gamma$
can be given as follows. Let $v_1, \cdots, v_N$ 
be the successive vertices of $\gamma$, which correspond to faces $v_1^*,
\cdots, v_N^*$ of $\Gamma^*$ and therefore to simplices $S_1,\cdots,
S_N$ of $\cS(\Gamma)$. The union of the interiors of the triangles
$v_i^*, i=1, \cdots, N$, along with the edges $e_i$ of $\Gamma^*$ between
$v_i^*$ and $v_{i+1}^*, i=1,\cdots, N$, is an open strip, which is called 
$S_\cup$ here. Its boundary $\dr S_\cup$ has two connected components,
$\dr_lS_\cup$ and $\dr_rS_\cup$, respectively on the left and on the
right of $\gamma$. We paramerized them with
the orientation coming from the orientation of $\gamma$, that is, 
the orientation of $\dr_rS_\cup$ is compatible with the orientation
of $S_\cup$, but the orientation of $\dr_lS_\cup$ is opposite to the
natural orientation of $\dr S_\cup$. Note that both $\dr_lS_\cup$
and $\dr_rS_\cup$ can, in some cases, be reduced to points.

Let $i\in \{ 1,\cdots, N\}$. The edge $e_i$ between $v_i^*$ and 
$v_{i+1}^*$ corresponds to a ``horizontal'' edge $e'_i$ of $\cS(\Gamma)$
which is contained in two simplices of $\cS(\Gamma)$, $S_i$ and $S_{i+1}$.
Under the deformation $\Theta'_\gamma$, the angle of $S_i$ at $e'_i$
has differential equal to $1$, and the angle of $S_{i+1}$ at $e'_i$
has differential equal to $-1$. Morever, the triangle $v_i^*$ has two
edges which are contained in the interior of 
$S_\cup$, namely $e_{i-1}$ and $e_i$, 
and exactly one edge, say $E_i$, which is either on $\dr_lS_\cup$
or on $\dr_rS_\cup$. In each case, we call $E_{i,-}$ and $E_{i,+}$
the endpoints of $E_i$, ordered according to the orientation of 
$\dr_lS_\cup$ and of $\dr_rS_\cup$ described above. Then $E_{i,-}$
(resp. $E_{i,+}$)
corresponds to a ``vertical'' edge $E'_{i,-}$ (resp. $E'_{i,+}$), 
which is contained in several simplices of $\cS(\Gamma)$, including
$S_i$. The angle of $S_i$ at $E'_{i,+}$ has differential equal to $1$,
and the angle of $S_i$ at $E'_{i,-}$ has differential equal to $-1$.
The same holds for all $i\in \{ 1, \cdots, N\}$, and 
the other angles remain constant. 

It is necessary to check that the first-order variation of the angles
of the simplices of $\cS(\Gamma)$ indeed defines a vector tangent 
to $\cA_E(\Gamma)$ at $\Theta$, i.e. that the linear constraints on
the angles of the simplices remain satisfied. Since those constraints
correspond to the ideal vertices of the simplices of $\cS(\Gamma)$, 
the result follows by considering 3 cases.
\begin{itemize}
\item At the ``central'' vertex $v_0$ (which is ideal since we're 
considering Euclidean metrics). Each $S_i, i=1,\cdots, N$, has 
exactly two edges ending at $v_0$ on which the angles vary, namely 
$E'_{i,-}$ and $E'_{i,+}$. Since the differentials of those angles
are opposite, the sum of the angles of the edges of $S_i$ adjacent
to $v_0$ remains equal to $2\pi$.
\item Let $v$ be a vertex of $\dr_lS_\cup$ (the same argument can be 
used for a vertex of $\dr_rS_\cup$). Then $v$ is the endpoint of an
even number $2p$ of edges of $\dr_lS_\cup$, say $E_{i_1},
E_{i_2}, \cdots, E_{i_{2p}}$. So:
$$ v = E_{i_1,+}=E_{i_2,-} = E_{i_3,+}=E_{i_4,-} = \cdots =
E_{i_{2p-1},+}=E_{i_{2p},-}~. $$
There are two kinds of simplices of $\cS(\Gamma)$ containing edges
on which the angle varies under the deformation $\Theta'(\gamma)$:
\begin{itemize}
\item the simplices $S_{i_{2k-1}}$ (resp. $S_{i_{2k}}$), $k=1,\cdots, p$,
on which the angle on the 
``vertical'' edge ending at $v$ has differential equal to
$1$ (resp. $-1$), while the angle on the ``horizontal'' edge 
$e'_{i_{2k-1}+1}$ (resp. $e'_{i_{2k}}$) has differential $-1$
(resp. $1$), 
\item the other simplices corresponding to triangles $v_i^*$
which contain $v$ but have no edge containing $v$ and contained
in $\dr_lS_\cup$. The angles of those simplices at edges containing
$v$ are constant except for two horizontal edges, on which the
differentials of the angles are opposite.
\end{itemize}
\item At all vertices of $\cS(\Gamma)$ which are neither $v_0$ nor 
a vertex of $\dr S_\cup$, no angle varies.
\end{itemize}
This shows that $\Theta'(\gamma)$ defines a vector tangent to
$\cA_E(\Gamma)$.

In addition, the first-order variation that we consider does not vary
the total angle around the ``vertical'' edges of $\cS(\Gamma)$ and 
the dihedral angles at the ``horizontal'' edges.

\begin{prop} \label{pr:gamma}
Let $\alpha:=\sigma(\Theta)$. Then
$\Theta'(\gamma)$ is tangent to $\sigma^{-1}(\alpha)$.
\end{prop}

\begin{proof}
The proof should be apparent from Figure \ref{fg:defo}. On the 
``horizontal'' edges of $\cS(\Gamma)$ which are not one of the 
$e'_i$, none of the angle varies, so the total angle remains constant.
On the $e'_i, i=1,\cdots, N$, the two angles vary in opposite ways,
so that again the total angle remains constant. 

Similarly, for ``vertical'' edges of $\cS(\Gamma)$, no angle varies except
when one of the endpoints is a vertex of $\dr S_\cup$. When one of the 
endpoints is on $\dr S_\cup$ (the other being $v_0$, as for all ``vertical''
edges) there is a even number of angles which vary, half of them with
differential $1$ and the other with differential $-1$, and the sum of the
angles remains constant, as claimed.
\end{proof}

A simple but already useful instance of a deformation is obtained by
taking as $\gamma$ the boundary of a face of $\Gamma$. Then $\dr_lS_\cup$
is reduced to a point. 

\paragraph{Critical points of the volume functional.}

We are now equiped to prove the Euclidean version of the 
key technical point of this subsection:
the critical points of the restriction of $\cV$ to the level sets of
$\sigma$ are the angle assignments for which the simplices of $\cS(\Gamma)$ 
can be isometrically glued along their faces, in such a way that it
is possible to choose horosphere centered at all ``ideal'' vertices
which ``match'' along the 2-dimensional faces. 

\begin{lemma} \label{lm:critique-e}
Let $\Theta\in \cA_E(\Gamma)$, and let $\alpha:=\sigma(\Theta)$.
Then $\Theta$ is a critical point of the restriction of $\cV$
to $\sigma^{-1}(\alpha)$ if and only if:
\begin{enumerate}
\item the hyperbolic metrics defined by $\Theta$ on the simplices of 
$\cS(\Gamma)$ can be glued isometrically along their faces,
\item the resulting hyperbolic metric can be completed on the ``vertical''
  edges of $\cS(\Gamma)$ to obtain a hyperbolic with conical singularities on
  $\cS(\Gamma)$ which moreover is consistent.
\end{enumerate}
\end{lemma}

\begin{proof}
Suppose first that conditions (1) and (2) hold. Let $h$ be the hyperbolic
metric with conical singularities obtained by gluing the hyperbolic metrics
on the simplices of $\cS(\Gamma)$. Then, by definition of a
consistent hyperbolic metric (with conical singularities) it is possible
to choose for each simplex $S$ of $\cS(\Gamma)$ a horosphere $H(S)$ centered
at the ``central'' vertex $v_0$ in such a way that, if $S$ and $S'$ are two
simplices which share a face $f$, then $H(S)\cap f=H(S')\cap f$. It is then
also possible to choose, for each ideal vertex $v$ of $\cS(\Gamma)$ (other
than $v_0$) and each simplex $S$ adjacent to $v$ a horosphere $H'(S,v)$ 
centered at $v$, which is tangent to $H(S)$, and this condition defines
$H'(S,v)$ uniquely. Moreover, if $S$ and $S'$ are two simplices which share
a face $f$ containing $v$, then $H'(S,v)\cap f = H'(S',v)\cap f$, because
both curves are horocycles which are tangent at the edge $[v_0,v]$ to 
$H(S)\cap f=H(S')\cap f$.

Let $\Theta'$ be a first-order variation of $\Theta$ which is tangent to 
$\sigma^{-1}(\alpha)$. The same argument as in the proof of Lemma
\ref{lm:critique-h}, with this choice of horospheres, shows that
$d\cV(\Theta')=0$. So conditions (1) and (2) imply that $\Theta$ is a critical
point of the restriction of $\cV$ to $\sigma^{-1}(\alpha)$.

Conversely, suppose first that condition (1) does not hold. This means
that there is a ``vertical'' face $f$ of $\cS(\Gamma)$, between two simplices
which can not be isometrically glued. Since two ideal triangles are always
isometric, $f$ has either one or two ideal vertices (including $v_0$).
Figure \ref{fg:defo2} describes a simple first-order variation of the 
angles of the simplices containing $f$ which leads to a non-zero variation
of the sum of the volumes of the two simplices; the strictly hyperideal
vertices are marked with circles. 

\begin{figure}[ht] 
\centerline{\psfig{figure=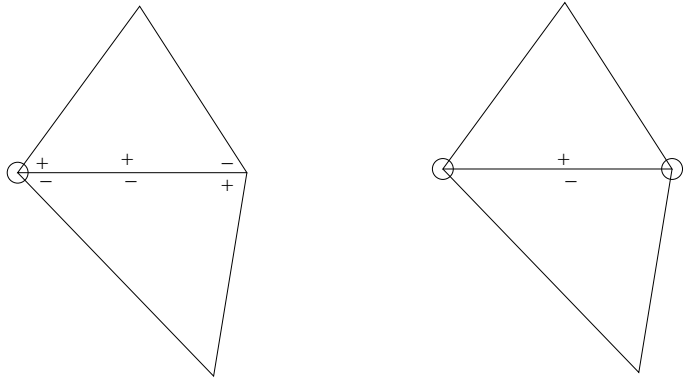,height=5cm}}
\caption{Deformations close to a face without gluing, with two (left) or one
  (right) ideal vertices.}\label{fg:defo2}
\end{figure} 

Checking that it is
indeed the case is quite obvious when $f$ has no ideal vertex beyond
$v_0$ (on the right in Figure \ref{fg:defo2}), since there is no 
linear constraint at the strictly hyperideal vertices, and the fact that
the total volume varies is a direct consequence of the Schl\"afli
formula, since the length for the two metrics of the ``horizontal''
edge in $f$ has to be different. 

Now suppose that $f$ has two ideal vertices at $v_0$ and $v_1$, and a strictly
hyperideal vertex at $v_2$ (pictured on the left in Figure
\ref{fg:defo2}). Note that the deformation pictured there is indeed tangent
to $\cA_E(\Gamma)$, because the linear constraints at $v_0$ and at 
$v_1$ hold. 
It is helpful to choose, in both simplices, horospheres centered at $v_0$
which are tangent to the hyperbolic plane dual to $v_2$, and horospheres
centered at $v_1$ which are tangent to the horospheres centered at $v_0$. This
way the length of the ``vertical'' edges of $f$ are zero, and it is clear 
that the deformations pictured in Figure
\ref{fg:defo2} induces a non-zero first-order variation of the sum of the 
volumes of the simplices. Finally, this deformation is tangent to
$\sigma^{-1}(\alpha)$, because the sum of the angles at all 3 edges of $f$
remains constant.

The next step is when $\Theta$ defines hyperbolic metrics on the simplices of
$\cS(\Gamma)$ which can be isometrically glued, but so that the translation
component of the holonomy around one of the ``vertical'' edges of
$\cS(\Gamma)$, say $e$, is non-zero. In other terms, it is possible to choose 
horosphere centered at $v_0$ for each of the simplices $S_1,\cdots, S_N$
containing $e$, so that the horospheres on both sides of $S_i\cap S_{i+1}$
match for $i=1,\cdots,N-1$, but not for $i=N$. 
The endpoints of $e$ are $v_0$ and a
vertex $v$ of $\Gamma^*$, corresponding to a face $v^*$ of $\Gamma$. Consider
the first-order variation $\Theta'_{\dr v^*}$. Proposition \ref{pr:gamma}
shows that it is tangent to $\sigma^{-1}(\alpha)$, while the Schl\"afli
formula shows that $d\cV(\Theta'_{\dr v^*})\neq 0$.

The last step is to suppose that the hyperbolic metrics on the simplices can
be glued, that the resulting metric on $\cS(\Gamma)$ minus its edges has zero
translation holonomy around the ``vertical'' edges, but that it is not
consistent. Definition \ref{df:consistent} then means that it is not possible
to choose horospheres centered at $v_0$, in all the simplices of
$\cS(\Gamma)$, in such a way that their intersections with the ``vertical''
faces match. This means that there exists a closed path $\gamma$ in $\Gamma$
(which has to be homotopically non-trivial), with vertices $v_1, \cdots, v_N$
corresponding to simplices $S_1, \cdots, S_N$ of $\cS(\Gamma)$, so that, if
one choose horospheres $H_1, \cdots, H_N$ 
centered at $v_0$ in the $S_i$ such that their
intersections with $S_1\cap S_2, S_2\cap S_3, \cdots, S_{N-1}\cap S_N$ match,
then the intersections of the horosphere $H_N$ 
in $S_N$ with $S_N\cap S_1$ does not
match with the intersection with that same face of the horosphere $H_1$ 
in $S_1$; let $\delta$ be the oriented distance between them.

For each ideal vertex $v$ of the $S_i$ other than $v_0$, we choose horospheres
such that for all $i=1, \cdots, N$, if  $v\in S_i\cap S_{i+1}$ then
the horospheres chosen at $v$ in $S_i$ and in $S_{i+1}$ match on $S_i\cap
S_{i+1}$. 

We now use the 
Schl\"afli formula to compute the variation $d\cV(\Theta'_\gamma)$
under the first-order deformation $\Theta'_\gamma$ of the angles assigned to
the simplices of $\cS$. With the choice of horospheres described it is quite
clear that the contributions corresponding to all the edges except the two
``vertical'' edges of $S_1\cap S_N$ cancel out, so that 
$d\cV(\Theta'(\gamma))=-2\delta\neq 0$. 
\end{proof}

\paragraph{Euclidean circle patterns.}

The results of the previous paragraph lead naturally to an analog, for
Euclidean surfaces with conical singularities, of Lemma \ref{lm:rigidity-h}.

\begin{lemma} \label{lm:rigidity-e}
Let $C$ be a hyperideal circle pattern on a closed surface
$\Sigma$, with respect to an Euclidean metric
with conical singularities at the center of the dual circles. The 
first-order deformations of $C$, among hyperideal circle patterns with
the same incidence graph and the same set of ideal dual circles 
(considered up to homothety) are 
parametrized by the first-order variations of the intersection angles
between the circles and of the singular curvatures at the
centers of the dual circles, under the condition that, for each ideal dual
circle $c'$, the sum of the intersection angles between the adjacent 
principal circles varies as the singular curvature at $c'$.
\end{lemma}

The proof follows exactly the same pattern as the proof of Lemma
\ref{lm:rigidity-h}, the basic block is now Lemma \ref{lm:critique-e}.


\section{Compactness}

This section deals with the properness of the map sending a hyperideal circle
pattern on a singular surface to its intersection angles and singular
curvatures. This can be stated in terms of compactness: a sequence of such
circle patterns, with given incidence graph and such that the intersection
angles and the singular curvatures converge to a ``reasonable'' limit, has a
subsequence which converges to a ``good'' circle pattern.

The methods used here are elementary, based solely on considering the possible
behavior of the circles relative to the underlying metric. Other approaches
are possible, in particular using the 3-dimensional geometry of the hyperbolic
simplicial complex associated to a circle pattern. 

\subsection{Statement}

We state here the compactness lemma which is
used in the proofs of the main theorems. 

\begin{lemma} \label{lm:compact}
Let $\Gamma$ be a graph embedded in the closed surface $\Sigma$, 
and let $(C_n)_{n\in \N}$ be a
sequence of hyperideal circle patterns with incidence graph
$\Gamma$, on $\Sigma$ with a sequence of metrics $(h_n)_{n\in \N}$ 
of constant curvature $K_n\leq 0$, of area equal to $1$, with
conical singularities at the centers of the dual circles. Suppose that 
$K_n\rightarrow K_\infty$, and that the
intersection angles and the singular curvatures are described by functions
$\theta_n:\Gamma_1 \rightarrow (0,\pi)$ and $\kappa_n:\Gamma_2\rightarrow
(-\infty, 2\pi)$, 
which are converging towards functions $\theta:\Gamma_1 \rightarrow (0,\pi)$
and $\kappa:\Gamma_2\rightarrow \R$, respectively, such that $\theta$ and
$\kappa$ still satisfy the hypothesis of Theorem \ref{tm:euclid} or
\ref{tm:hyperb} (depending on whether $K_\infty=0$ or $K_\infty<0$). Then,
perhaps after extracting a sub-sequence,
$(C_n)_{n\in \N}$ converges to a circle pattern $C$ on $\Sigma$, 
on a metric of constant curvature $K_\infty$ with conical 
singularities at the centers of the dual circles, with
intersection angles given by $\theta$ and singular curvatures given by
$\kappa$.  
\end{lemma}

The proof proceeds in two steps. The first is to prove that the radii of the
principal circles remain bounded, i.e. they can not become very
small or very large.  
This implies that the distances between the singular points can not become
too small. Direct geometric arguments then show that those
metrics with conical 
singularities have a converging sub-sequence.

\subsection{Collapsing circles}

The arguments given below use a graph which is canonically associated to
$\Gamma$ in a simple manner.

\begin{df}
Let $\Gamma^s$ the graph, embedded in $\Sigma$, which has: 
\begin{itemize}
\item one vertex for each vertex of $\Gamma$, and one for each 
vertex of $\Gamma^*$ (i.e. for each face of $\Gamma$),
\item one edge between two vertices if they correspond to 
adjacent vertices of $\Gamma$, or if they correspond to a 
face of $\Gamma$ and a vertex of $\Gamma$ contained in that face.
\end{itemize}
\end{df}

There is a simple interpretation of $\Gamma^s$, its vertices correspond
to the circles which are either principal or dual circles in the circle
pattern, and the edges correspond to pairs of adjacent circles. Clearly 
$\Gamma^s$ has triangular faces, with each face having
one vertex corresponding to a dual circle and two corresponding to 
principal circles.

The proof of the compactness lemmas above is based on a simple proposition.

\begin{prop} \label{pr:collapse}
Under the hypothesis of Lemma \ref{lm:compact}, the radii of the principal
circles remain bounded between two positive constants.
\end{prop}

\begin{proof} \label{pr:radii}
We have supposed that the area is always equal to $1$, it follows quite
directly that the radii of the circles can not go to infinity. Suppose that,
after taking a sub-sequence, the radii of some of the principal
circles go to $0$.
Then this also happens to one of the dual circle adjacent to each
``collapsing'' principal circle; indeed each principal circle is
adjacent to at least 3 dual circles (because the principal circles correspond
to the vertices of the embedded graph $\Gamma$). A ``very small''
principal circle has to be orthogonal to at least 3 dual circles, which
are disjoint, and this is possible only if at least one of those adjacent dual
circles is also very small.

Consider, in the graph $\Gamma^s$, the set $V_0$ of vertices corresponding
to circles which have a radius going to $0$. The remark above shows that 
$V_0$ contains at least one point corresponding to a principal circle, and
one corresponding to an adjacent dual circle.
Let $\Gamma^s_0$ be the
subgraph of $\Gamma^s$ which has as its vertices the vertices in $V_0$
and as edges the edges of $\Gamma^s$ with both endpoints in $V_0$, and 
let $\Gamma^s_1$ be a connected component of $\Gamma^s_0$ containing
at least one vertex corresponding to a principal circle.

Let $\dr\Gamma_1^s$ be the set of edges of $\Gamma^s$ which are, in a face of 
$\Gamma^s$ with exactly one vertex in $\Gamma_1^s$, the edge opposite 
to that vertex. $\dr\Gamma_1^s$ is by construction an union of closed
curves, which we call $\dr_1\Gamma_1^s, \cdots,  \dr_N\Gamma_1^s$.
By constructions all the circles corresponding to the vertices of 
$\Gamma_1^s$ ``collapse'' to a point in the limit surface, while
this happens to none of the circles correspondings to the vertices of 
$\dr \Gamma_1^s$. 

We associate to $\Gamma_1$ an admissible domain $\Omega$ in $(\Sigma, \Gamma)$,
as follows. For each edge $e$ of $\dr\Gamma_1^s$ with endpoints corresponding
to vertices of $\Gamma$, we consider the corresponding edge
$e'$ of $\Gamma$. For each vertex $v$ of $\dr\Gamma_1^s$ corresponding to a
face $f$ of $\Gamma$, the two vertices adjacent to $v$ in $\dr\Gamma_1^s$
correspond to vertices $v_1$ and $v_2$ of $\Gamma$, both contained in $f$.
Then we consider a segment $s$ in $f$ going from $v_1$ to $v_2$. 

Let $\gamma$ be the union of
those two types of segments, it is a finite union of closed curves $\gamma_1,
\cdots, \gamma_N$, with $\gamma_i$ 
corresponding to $\dr_i\Gamma_1^s$. It is also 
simple to check that $\gamma$ bounds a domain $\Omega$ in $\Sigma$, 
and that $\Omega$ is an admissible domain in
$(\Sigma, \Gamma)$. Below we set $\dr_i\Omega:=\gamma_i$.

All vertices of $\Gamma_1^s$ 
correspond to circles which ``collapse'' to a point. Since $\Gamma_1^s$ is
connected, all those circles actually ``collapse'' to the same point, which
we call $x_0$. Since the circles corresponding to vertices of 
$\dr\Gamma_1^s$ have
radii converging to positive numbers, the union of the disks bounded by those
circles converges to a flat surface $S_1$ with 
one singular point (which corresponds to $x_0$).
More precisely, the limit surface is obtained by considering, for each
connected component $\dr_i\Omega$ of $\dr\Omega$, 
a flat surface $S_i$ with boundary, with one
conical singularity, at the point $x_i$ corresponding to $x_0$ under the
gluing of those flat surfaces. Moreover, in each of the flat surfaces
$S_i, 1\leq i\leq N$, each of the circles contains the singular point. 

\begin{figure}[ht] 
\centerline{\psfig{figure=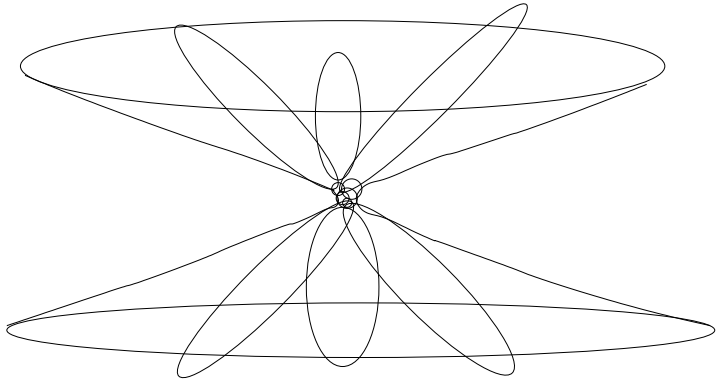,height=5cm}}\label{fg:annulus}
\caption{A collapsing annulus, and some of the adjacent circles.}
\end{figure} 

For each $n\in \N$, let $\Omega_n$ be an open domain in 
$(\Sigma, h_n)$ isotopic to $\Omega$, with
piecewise smooth boundary, which contains all the conical singularities
corresponding to the faces of $\Gamma$ contained in $\Omega$, but no
other conical singularity. Let $k_n$ be the integral of the geodesic 
curvature of $\dr\Omega_n$, the Gauss-Bonnet theorem shows that: 
$$ k_n = 2\pi \chi(\Omega) - \sum_{f\subset \Omega} \kappa_n(f) - 
K_nA(\Omega_n)~. $$

We can choose the $\Omega_n, n\in \N$, so that $(\Omega_n)_{n\in\N}$ 
converges, as
$n\rightarrow\infty$, to a limit domain $\Omega_\infty$ 
in the flat surface with
boundary and with one singular point which we already mentioned above. 
Let $k$ be the integral of the geodesic curvature of $\dr\Omega_\infty$.
Then, taking the limit in the previous equation:
$$ k = 2\pi \chi(\Omega) - \sum_{f\subset \Omega} \kappa(f) 
- K_\infty A(\Omega_\infty)~. $$

Let $k_1, \cdots, k_N$ the singular curvatures at $x_i$ of the flat surface
$S_i$, $1\leq i\leq N$. Then, for each of those surfaces, the 
integral of the geodesic curvature of the boundary of any domain 
containing the singular point is $2\pi-k_i$ minus $K_\infty$ times
the area of the corresponding disk, by the Gauss-Bonnet 
Theorem. It follows that: 
\beq \label{eq:egalemme}
\sum_{i=1}^N (2\pi-k_i) = 2\pi\chi(\Omega) - \sum_{f\subset \Omega}
\kappa(f)~. \eeq

To each $\dr_i\Omega, 1\leq i\leq N$ 
we associate a sequence of circles containing the cone
point $x_i$, namely the circles corresponding to the vertices of
$\dr_i\Gamma_1^s$.
Lemma \ref{lm:bouquet}, applied to those circles, shows that the sum
of their intersection angles in the limit
is equal to $2\pi-k_i$. But their intersections
angles are of two types: 
\begin{itemize}
\item each dual circle intersects orthogonally two principal circles, so that
  the total contribution of each dual circle is equal to $\pi$,
\item each edge of $\dr_i\Gamma_1^s$, with endpoints corresponding to two
  principal circles, corresponds to an edge of the connected component of
  $\dr\Omega$ corresponding to $\gamma_i$.
\end{itemize}
Consider equation (\ref{eq:egalemme}), replacing 
the left-hand side by the sum of the intersection angles between the angles, 
we obtain that: 
$$ \sum_{e\in \dr\Omega}\theta(e) + m(\Omega)\pi =  2\pi\chi(\Omega) -
\sum_{f\subset \Omega} \kappa(f)~, $$
which means precisely that the inequality in condition (2) of Theorem
\ref{tm:euclid} or Theorem \ref{tm:hyperb} 
is an equality. Note also that $\Gamma_1^s$ contains at least
one vertex corresponding to a principal circle, and it follows that $\Omega$
is not a face of $\Gamma$. So the equality above is a contradiction.
\end{proof}

\subsection{Proof of the compactness lemma}

Proposition \ref{pr:collapse}, along with its hyperbolic analog, 
is the main tool for the proof of Lemma \ref{lm:compact}. 
Once it is known, the compactness
results are obtained by using elementary geometric arguments concerning the
constant curvature metrics on the surfaces with boundary obtained by
considering a finite (but increasing) set of overlapping disks.

Given $\Gamma$, a hyperideal circle pattern with incidence graph $\Gamma$, and
its underlying constant curvature $K$ metric, are determined by a finite set
of data: 
\begin{itemize}
\item the curvature $K$,
\item the radii of the principal circles,
\item the singular curvatures at the centers of the dual circles,
\item an ``angle'' data describing the relative positions of the circles
  intersecting a given principal circle, say $C_i$: if $C_j$ and $C_k$
  intersect $C_i$, it is the angle under which the centers of $C_j$ and $C_k$
  are ``seen'' from the center of $C_i$.
\end{itemize}
Indeed given this information the radii and position of the dual circles
follow, and it is not difficult to check that this completely determines the
underlying metric and circle pattern.

Proposition \ref{pr:collapse} shows that the radii of the principal circles
remain bounded, and the singular curvatures at the centers of the dual circles
converge. The angle data are defined in a compact set. Therefore after taking
a subsequence all those data converge. Since the radii of the principal
circles remain bounded from below, the injectivity radii of the underlying 
sequence of metrics remain bounded from below. 

Notice that the combinatorics of the circle pattern remains the same in the
limit. Clearly two circles corresponding to adjacent vertices of $\Gamma$
remain so in the limit, because they still intersect and their intersection
points remain in interstices. Conversely, two circles not adjacent cannot
become adjacent in the limit, which would mean that either they become
tangent and their tangency point is in an interstice, or that their 
intersection points are contained, in the limit, in interstices. 
This would be forbidden by
the constraints on a hyperideal circle pattern. For instance if
$C_1$ and $C_2$ are two principal circles which are not adjacent in the
sequence, but which are both orthogonal to a dual circle $C'$, then, if $C_1$
and $C_2$ where tangent in the limit, it would imply that (still in the limit)
$C'$ would have to be on either side of $C_1\cup C_2$, and that it could not
be adjacent to some of the principal circles which corresponds to vertices of
$\Gamma$ contained in the face corresponding to $C'$ ``between'' the vertices
corresponding to $C_1$ and $C_2$.

Moreover, the distance between the centers of the
dual circles is also bounded from below, because, as is not difficult to see,
the centers of two dual circles can ``collapse'' only if the radii of those
circles go to $0$ and the intersection between two principal circles, 
both adjacent to those two dual circles, goes to $\pi$. It follows that the
limit metric is a constant curvature metric with conical singularities and the
corresponding circle pattern has combinatorics given by $\Gamma$ and 
intersection angles and singular curvatures prescribed by $\theta$ and
$\kappa$. 


\section{Proofs of the main results on closed surfaces}

\subsection{Outline}

The previous section contain results showing that the map sending a hyperideal
circle pattern on a singular surface to its intersection angles and the
singular curvatures at the singular points is a local homeomorphism, and also
that it is proper. It follows that it is a covering of the spaces of possible
intersection angles and singular curvatures --- in the different contexts
considered here --- by the corresponding spaces of hyperideal circle
patterns. However it remains necessary to prove that the number of inverse
images of each set of intersection angles and singular curvatures, which is
locally constant, is equal to $1$.

Two different arguments will be used. For the uniqueness --- the fact that
each set of intersection angles and singular curvatures has at most one
inverse image --- we use the fact that circle patterns correspond to critical
points of the volume over polytopes, while the volume is a strictly concave
function. This is done in subsection 7.2. 

It is also necessary to prove that, for any given combinatorics, there are
some assignments of intersection angles and singular curvatures which are
indeed realized by a hyperideal circle pattern. This is done first for
hyperbolic surfaces in subsection 7.3, using the Orbifold Hyperbolization
Theorem (see \cite{boileau-porti}) and a double doubling argument, which shows
that, to some strictly hyperideal circle pattern on a hyperbolic surface,
it is possible to associate a hyperbolic orbifold, and conversely.

Then, in subsection 7.4, we use the existence result for some hyperbolic
circle patterns to obtain a similar existence result for some Euclidean circle
patterns, using an approximation by hyperbolic circle patterns and Lemma
\ref{lm:compact}. 

\subsection{Uniqueness}

The proof of the uniqueness of circle patterns with a given combinatorics,
intersection angles and singular curvatures, follows quite directly from the
constructions made in section 5 to prove the infinitesimal rigidity. 

\paragraph{Closed hyperbolic surfaces.}

Consider first the case of case of circle patterns with 
incidence graph $\Gamma$ on a closed surface $\Sigma$ of
genus at least $2$ with a singular hyperbolic metric. 
We have seen in section 5 the definition of a polytope $\cA_H(\Sigma,
\Gamma)$, on which two functions are defined: 
\begin{itemize}
\item the volume function, $\cV$, which is strictly concave,
\item the affine function $\sigma:\cA_H(\Sigma,\Gamma)\rightarrow
  \cD(\Sigma,\Gamma)$. 
\end{itemize}
Moreover, the circle patterns are associated to the critical points of the
restriction of $\cV$ to the level sets of $\sigma$ (see Lemmas
\ref{lm:compatib-h} and \ref{lm:critique-h}). At each such
critical point, the value of $\sigma$ determines the intersection angles and
the singular curvatures of the circle pattern. 
Since $\sigma$ is affine, its level sets are affine subsets of $\cA_H(\Gamma,
\Sigma)$. Moreover, the restriction of $\cV$ to those level sets remains
strictly concave, so that $\cV$ has at most one critical point on each level
set. This shows that each family of intersection angles and singular
curvatures is obtained on at most one circle pattern, proving the uniqueness
in Theorem \ref{tm:hyperb}.

\paragraph{Euclidean surfaces.}

The same argument can be used for Euclidean circle patterns on closed
surfaces, 
i.e. Theorem \ref{tm:euclid}. The relevant polytope is now $\cA_{E}(\Gamma)$
(see Definition \ref{df:sigma}), but it remains true that circle patterns are
in one-to-one correspondence with the critical points of the restriction of
$\cV$ to the level sets of $\sigma$ (see Lemma \ref{lm:critique-e}). The
uniqueness of the circle patterns in Theorem \ref{tm:euclid} thus follows.

\subsection{Existence for hyperbolic metrics}

\paragraph{The Orbifold Hyperbolization Theorem.} 

We will use an important statement, discovered by Thurston (who gave an
outline of the proof) and which has recently been extensively proved
\cite{boileau-porti,BLP,CHK}. 

\begin{thm}[Thurston, see \cite{boileau-porti,BLP}]
Let $M$ be a compact, connected, orientable, irreducible 
3-dimensional orbifold with non-empty singular locus, which has either
non-empty singular locus or infinite fundamental group. If $M$ is
topologically atoroidal, then it admits a hyperbolic, Euclidean, or Seifert
fibered structure.
\end{thm}

It will be applied here for Haken orbifolds, for which a simpler proof has
been available for some time. 

To check that an orbifold with singularities along closed curves
has a hyperbolic
orbifold structure, it is sufficient (cf \cite{BLP}) to check that:
\begin{itemize}
\item every essential sphere intersects singular lines on which the sum of the
  singular curvatures is at least $4\pi$,
\item every essential torus intersects at least one singular line.
\end{itemize}
This statement will be applied to a 3-dimensional orbifold constructed 
by simple doubling constructions from a hyperbolic circle pattern.

\paragraph{Doubling constructions.}

We first consider a hyperideal circle pattern $C$ on a hyperbolic surface
$\Sigma$, with incidence graph $\Gamma$, singular curvatures given by
$\kappa:\Gamma_2\rightarrow (-\infty, 2\pi)$, and angle intersections given by
$\theta:\Gamma_1\rightarrow (0,\pi)$. We suppose that:
\begin{itemize}
\item no interstice of $C$ is reduced to a point,
\item for each $f\in \Gamma_2, \kappa(f)=2\pi-2\pi/k$, for some integer $k\geq
  2$, 
\item for each $e\in \Gamma_1, \theta(e)=\pi-\pi/k$, 
for some integer $k\geq 2$. 
\end{itemize}
We have seen in section 4 how to construct from $C$ a cellular complex $\cS$,
and a complete hyperbolic metric with cone singularities $H(C)$ on $\cS$. 
Moreover, each of the 3-dimensional cells of $\cS$ has one face corresponding
to a face of $\Gamma^*$, and another corresponding to a vertex $v_0$ which is
the same for all 3-cells of $\cS$. In the hyperbolic metric $H(C)$, $v_0$
corresponds to a strictly hyperideal vertex --- because the circle pattern $C$
is on a hyperbolic surface --- while all other vertices correspond to strictly
hyperideal vertices --- because we have supposed that the interstices of $C$
are ``strictly hyperideal'' (they are not reduced to points).

We define another cellular complex $\cS_t$ with a hyperbolic cone-metric
$H_t(C)$ from $(\cS, H(C))$ by truncating it at the dual planes of all its
vertices (which are all strictly hyperideal). The resulting cone-manifold is
compact, homeomorphic to the product of $\Sigma$ by an interval, 
with a boundary made of 3 different parts:
\begin{itemize}
\item $\dr_0\cS_t$, defined as the union of the faces dual to $v_0$ in all the
  simplices of $\cS$. It follows from the construction made in section 4 that
  $\dr_0\cS_t$ is connected and totally geodesic.
\item $\dr_t\cS_t$, the union of the faces dual to all the vertices of $\cS$
  other than $v_0$. Also by construction, $\dr_t\cS_t$ is a disjoint union of
  totally geodesic polygons, each corresponding to one of the vertices of
  $\cS$ other than $v_0$.
\item $\dr_p\cS_t$, which is the part of $\dr\cS_t$ which was already
  contained in $\dr \cS$. $\dr_p\cS_t$ is topologically $\Sigma$ with one
  disked removed for each face of $\Gamma$.
\end{itemize}
In addition, $\cS_t$ has one singular line for each face $f$ of $\Gamma$, with
total angle equal to $\kappa(f)$. $\dr_p\cS_t$ is polyhedral, it has
one edge for each edge $e$ of $\Gamma$, with exterior dihedral angle equal to
$\theta(e)$. 

We now define a second simplicial complex $\cS_2$ with a hyperbolic
cone-metric $H_2(C)$ by gluing two copies of $(\cS_t, H_t(C))$ isometrically 
along $\dr_0\cS_t$ and $\dr_t\cS_t$. 
The resulting cone-manifold is topologically
the product of $\Sigma$ by an interval, with, for each face of $\Gamma$, one
handle connecting one face with the other. It also has one singularity along a
closed line for each face $f$ of $\Gamma$, with total angle
$2\pi-\kappa(f)$. 
Its boundary is ``polyhedral'', with one face corresponding to
each vertex of $\Gamma$, one edge --- which is a closed curve --- for each
edge of $\Gamma$, and no vertex.

Finally we define a third cone-manifold, $\cS_4$, with a hyperbolic
cone-metric $H_4(C)$, by gluing two copies of $(\cS_2, H_2(C))$ isometrically 
along their boundary. By construction, the total angle around each singular
line --- which corresponds either to a face or to an edge of $\Gamma$ --- is
equal to $2\pi/k$, for some $k\geq 2$. Therefore, $(\cS_4, H_4(C))$ is a
hyperbolic orbifold.

\paragraph{Construction of hyperbolic orbifolds.}

We now consider the same constructions in a topological manner, with no
hyperbolic metric involved. Consider a graph $\Gamma$ embedded in 
$\Sigma$ and two
functions $\kappa:\Gamma_2\rightarrow (0,2\pi)$ and
$\theta:\Gamma_1\rightarrow (0,\pi)$, satisfying the hypothesis of Theorem
\ref{tm:hyperb}, such that $\kappa(f)=2\pi(1-1/k)$ for
some $k\geq 2$ and that $\theta(e)=\pi-\pi/k$ 
for some $k\geq 2$, for all face $f$
(resp. for all edge $e$) of $\Gamma$. By the same construction as above, it is
possible to construct a cone-manifold $\cS_4$, with singularities along closed
curves, with angle $2\pi/k$ (for some $k\geq 2$) along each of those curves,
so $\cS_4$ is an orbifold.

\begin{lemma} \label{lm:orbifold}
$\cS_4$ admits a hyperbolic orbifold structure.
\end{lemma}

\begin{proof}
$\cS_4$ was constructed above by gluing along their boundary two copies of
$\cS_2$, which we call $\cS_{2,+}$ and $\cS_{2,-}$. Each is topologically the
product of $\Sigma$ by an interval, with one handle for each face of
$\Gamma$. The boundary of $\cS_{2,+}$ and $\cS_{2,-}$ contain some singular
curves of $\cS_4$, those corresponding to an edge of $\Gamma$. To each such
curve is associated an ``exterior angle'', which is half the singular 
curvature of the corresponding curve in $\cS_4$.

Consider an essential torus $T^2\subset \cS_4$. 
Suppose first that $\Sigma$ is not a torus. Given the topology of
$\cS_{2,+}$ and $\cS_{2,-}$, they do not contain any essential torus, so the
intersection of $T^2$ with $\cS_{2,+}$ (resp. $\cS_{2,-}$) contains an
essential annulus $A_+$ (resp. $A_-$). The boundary of $A_+$ is made of two
curves on $\dr \cS_2$, each of which intersects at least one curve
corresponding to an edge of $\Gamma$. So, in $\cS_4$, $T^2$ intersects at
least one of the singular curves of $\cS_4$.

If $\Sigma$ is a torus, the same argument holds except in the case where 
$T^2$ is contained in $\cS_{2_+}$ (or in $\cS_{2,-}$); then it has to 
intersect at least one of the singular curves corresponding to the faces of
$\Gamma$.

Now let $S$ be an essential sphere in $\cS_4$. Since
$\cS_{2,+}$ and $\cS_{2,-}$ are irreducible, the intersection of $S$ with
$\cS_{2,+}$ (resp. $\cS_{2,-}$) contains at least one essential disk $D_+$
(resp. $D_-$). The boundary of $D_+$ is a closed curve in $\dr \cS_{2,+}$. 
Let $K_i(D_+)$ be the sum of the singular curvatures associated to the 
singular curves in $\cS_{2,+}$ which intersect the interior of $D_+$,
and let $\Theta_b(D_+)$ be the sum of the ``exterior angles'' associated
to the curves in $\dr \cS_{2,+}$ (corresponding to the edges of $\Gamma$)
which intersect $\dr D_+$. 
We will show that:
\begin{equation} \label{eq:disk}
K_i(D_+) + \Theta_b(D_+)>2\pi~. 
\end{equation}

$\cS_{2,+}$ is itself obtained by gluing two copies of $\cS_t$, which we
call $\cS_{t,+}$ and $\cS_{t,-}$, along one boundary component and along disks
corresponding to the faces of $\Gamma$. Let $D_{+,+}:=D_+\cap \cS_{t,+}$,
$D_{+,-}:=D_+\cap \cS_{t,-}$. After deformation of $D$, we can suppose
that $D_{+,+}$ and $D_{+,-}$ have a minimal number of connected components.
We consider three cases.
\begin{enumerate}
\item $D_{+,-}=\emptyset$ . 
Consider the projection of $D_{+,+}$ on $\dr\cS_{t,+}$, it is a disk
with boundary equal to $\dr D_{+,+}$ (which does not intersect the disks
in $\dr\cS_{t,+}$ corresponding to the faces of $\Gamma$). Consider
the sequence of edges of $\Gamma^*$ which intersect $\dr D_{+,+}$,
the corresponding edges of $\Gamma$ are the boundary of an admissible domain
$\Omega$ in $(\Sigma, \Gamma)$, which is topologically a disk, with
  $m(\Omega)=0$. Condition (2) in Theorem \ref{tm:hyperb} is that:
$$ \sum_{e\in \Gamma_1, e\subset\dr\Omega} \theta(e) \geq (2\chi(\Omega) -
m(\Omega))\pi - \sum_{f\in \Gamma_2, f\subset \Omega} \kappa(f)
= 2\pi - \sum_{f\in \Gamma_2, f\subset \Omega} \kappa(f)~, $$
which means precisely that $\Theta_b(D_+)>2\pi - K_i(D(_+)$, as needed.
\item $D_{+,+}=\emptyset$. The same argument can then be used with $D_{+,-}$ 
instead of $D_{+,+}$. 
\item $D_{+,-}$ and $D_{+,+}$ are both non-empty. 
There are then at least one connected component of $D_{+,+}$, 
called $D' _{+,+}$, with exactly one boundary
component contained in one of the disks, corresponding to a face $f$ 
of $\Gamma$,
along which $\cS_{t,+}$ is glued to $\cS_{t,-}$. Considering as
above the sequence of edges of $\Gamma$ corresponding to the edges of
$\Gamma^*$ intersected by $\dr D'_{+,+}$,
and adding a segment in $f$, we obtain an admissible domain $\Omega_{+,+}$ 
in $(\Sigma,\Gamma)$, which is a disk with $m(\Omega_{+,+})=1$. 
Condition (2) of Theorem \ref{tm:hyperb} now shows that the sum of the 
exterior angles associated to the edges intersected by $\dr D'_{+,+}$ is
larger than $\pi$ minus the sum of the singular curvatures associated
to the singular curves of $\cS_{2,+}$ intersected by $D'_{+,+}$.
The same argument can be use for a connected component $D'_{+,-}$ of 
$D_{+,-}$ with exactly one boundary component contained in one of the disks
on $\dr\cS_{t,-}$ corresponding to a face of $\Gamma$, and the sum of the
two contributions shows equation (\ref{eq:disk}).
\end{enumerate}
The same argument can also be applied to $D_-$, and it shows that, using the
same notations for $D_-$ as for $D_+$: $K_i(D_-) + \Theta_b(D_-)>2\pi$, while
by definition $\Theta_b(D_-)=\Theta_b(D_+)$. However the sum of the singular
curvatures of the singular curves which intersect $S$ is equal to $K_i(D_+) +
K_i(D_-) + 2\Theta_d(D_+)$, so it is larger than $4\pi$. This means 
that each essential sphere in $\cS_4$ intersects singular curves with a sum of
singular curvatures larger than $4\pi$.

It now follows from the Orbifold Hyperbolization Theorem that $\cS_4$
has a hyperbolic orbifold structure.
\end{proof}

Notice that condition (2) of Theorem \ref{tm:hyperb} is used here only in the
case where $\Omega$ is a disk.  This is quite normal since we consider
here circle patterns on hyperbolic surfaces with positive singular curvature
at the cone points, and we have already seen in section 1 that, in that case,
condition (2) is void except when $\Omega$ is a disk.

\paragraph{An existence result in some special cases.}

Let $\Gamma$ be the 1-skeleton of a cellular decomposition of a closed surface
$\Sigma$. We have seen above that, from a circle pattern of incidence graph
$\Gamma$, with intersection angles and singular curvatures satisfying 
some simple conditions, it is possible to ``built'' a hyperbolic orbifold
structure on $\cS_4$.
The argument given right above using the Orbifold Hyperbolization Theorem
shows that, conversely, given $\Gamma$ and two functions $\theta$ and
$\kappa$ satisfying the same conditions, there is a
unique hyperbolic orbifold structure on $\cS_4$ with total angles
around the singular curves given by $2\pi-\kappa$ and by $2(\pi-\theta)$
on singular curves corresponding to faces or to edges of $\Gamma$, 
respectively. 
By uniqueness, this metric has the same symmetry as $\cS_4$, so it
can be cut along a polyhedral surface to obtain two hyperbolic cone-manifolds
with boundary, $\cS_{2,+}$ and $\cS_{2,-}$, which are isometric. 

Then, again
by uniqueness, it is possible to ``cut'' $\cS_{2,+}$ along a totally
geodesic surface $\Sigma_0$ and a disjoint union of disks, to obtain 
a cone-manifold $\cS_t$. The boundary of $\cS_t$ has two connected
component, one corresponding to $\Sigma_0$ which is totally geodesic,
and another one which is polyhedral, with some faces corresponding to the 
faces of $\Gamma$ and other corresponding to the vertices of $\Gamma$.
Moreover two faces corresponding to the endpoints of an edge $e$ of
$\Gamma$ are adjacent, and the exterior dihedral angle between them is 
$\theta(e)$. It is possible to project orthogonally on $\Sigma_0$ the 
boundaries at infinity of the hyperbolic planes containing the faces of
$(\dr\cS_{2,+})\setminus \Sigma_0$ corresponding to vertices of $\Gamma$,
one obtains in this manner a hyperideal circle pattern on $\Sigma_0$,
with incidence graph $\Gamma$, singular curvatures given by $\kappa$, and
intersection angles given by $\theta$. 

This result can be stated as follows.

\begin{lemma} \label{lm:realisation}
Let $\Sigma$ be a closed surface, let $\Gamma$ be the 1-skeleton of 
a cellular decomposition of $\Sigma$, and let 
$\kappa:\Gamma_2\rightarrow (0,2\pi)$ and $\theta:\Gamma_1\rightarrow (0,\pi)$
be two functions such that: 
\begin{itemize}
\item $\Gamma, \kappa$ and $\theta$ satisfy the hypothesis of Theorem 
\ref{tm:hyperb}, 
\item the sum of the values of $\theta$ on the boundary of each face $f$
of $\Gamma$ is strictly larger than $2\pi - \kappa(f)$,
\item for each $f\in \Gamma_2$, $\kappa(f)=2\pi - 2\pi/k$, for some integer 
$k\geq 2$,
\item for each $e\in \Gamma_1$, $\theta(e)=\pi-\pi/k$, 
for some integer $k\geq 2$.
\end{itemize}
Then there is a strictly hyperideal circle pattern $C$ on $\Sigma$, with 
incidence graph $\Gamma$, singular curvatures given by $\kappa$, and
intersection angles given by $\theta$.
\end{lemma}

\paragraph{Proof of Theorem 1.5.}

Let again $\Sigma$ be a closed surface, and let $\Gamma$ be the 1-skeleton 
of a cellular decomposition of $\Sigma$. Following the notation from 
section 2, we call  
$\cD(\Gamma)$ the vector space of functions $\kappa, \theta$ satifying
the hypothesis of Theorem \ref{tm:hyperb}.

Notice that, when $\Sigma$ is a sphere, $\Gamma$ needs to have at least 3
faces for condition (1) of Theorem \ref{tm:hyperb} to apply.

\begin{lemma} \label{lm:speciaux}
There exist functions $\kappa,\theta$ which are in $\cD(\Gamma)$ and for 
which:   
\begin{enumerate}
\item the sum of the values of $\theta$ on the boundary of each face $f$ of 
$\Gamma$ is strictly larger than $2\pi-\kappa(f)$,
\item for each face $f\in \Gamma_2$, $\kappa(f)=2\pi-2\pi/k$, 
for some integer $k\geq 2$, 
\item for each $e\in \Gamma_1$, $\theta(e)=\pi-\pi/k$, 
for some integer $k\geq 2$.
\end{enumerate}
\end{lemma}

\begin{proof}
Let $k\in \N, k\geq 2$. Choose $\kappa(f)=2\pi(1-1/k)$ for all 
$f\in \Gamma_2$, and let $\theta(e)=\pi(1-1/k)$ for all $e\in \Gamma_1$. 
For each face $f$ of $\Gamma$, $f$ has at least $3$ edges, so the sum
of the values of $\theta$ on the boundary of $f$ is at least equal to
$3\pi-3\pi/k$. On the other hand, $2\pi- \kappa(f)=2\pi/k$. 
Since $k\geq 2$, condition (1) above is always satisfied.

We now have to prove that, if $k$ is chosen large enough, $\theta$ and
$\kappa$ correspond to a point in $\cD(\Gamma)$, i.e. that the conditions
in Theorem \ref{tm:hyperb} are satisfied. Condition (1) is simply that
the sum of $\kappa(f)$ over the faces of $f$ is larger than $\chi(\Sigma)$;
since $\kappa(f)>0$ for all $f\in \Gamma_2$, this is clearly true when 
$\Sigma$ is a torus or a higher genus surface. When $\Sigma$ is a sphere,
$\Gamma$ has at least 3 faces, so the condition is satisfied as soon as 
$k\geq 4$. 

Since $\kappa(f)>0$ for all faces of $\Gamma$, the right-hand side of
the second condition in Theorem \ref{tm:hyperb} is negative unless $\Omega$
is a disk with $m(\Omega)\in \{ 0,1\}$. Moreover, if $k\geq 4$, it
will also be negative as soon as $\Omega$ contains two faces of $\Gamma$,
since $2\pi-2\times 2\pi(1-1/k)<0$, and also if $\Omega$ contains one face 
of $\Gamma$ and $m(\Omega)=1$. So the only case left is when $\Omega$ is a face
of $\Gamma$ and $m(\Omega)=0$, and then, since each face has at least
3 edges, the condition applies if $k\geq 4$. 
\end{proof}

The proof of Theorem \ref{tm:hyperb} is now clear. Given $\Sigma$ and
$\Gamma$, we consider the map $\Phi$ on $\cC(\Gamma)$ sending a hyperideal
circle pattern to the data given by its intersection angles and singular 
curvatures. We have seen in section 3 that $\Phi_\Gamma$ takes its values
in $\cD(\Gamma)$. Moreover Lemma \ref{lm:rigidity-h} shows that $\Phi_\Gamma$
is a local homeomorphism, and Lemma \ref{lm:compact} shows that it is
proper, so that $\Phi_\Gamma$ is a covering of $\cD(\Gamma)$ by 
$\cC(\Gamma)$. It was also proved in subsection 7.1 that the number of 
inverse images of each point in $\cD(\Gamma)$ is at most equal to $1$.
But Lemma \ref{lm:speciaux} shows that $\cD(\Gamma)$ contains
at least one ``special'' point which, by Lemma \ref{lm:realisation}, has 
one inverse image. So $\Phi_\Gamma$ is a homeomorphism.

\subsection{Euclidean metrics}

\paragraph{Outline of the proof.}

The proof of Theorem \ref{tm:euclid} follows
the same path as the proof of Theorem \ref{tm:hyperb} except for the
last step. As for Theorem \ref{tm:hyperb}, given a closed surface $\Sigma$
and a graph $\Gamma$ embedded in $\Sigma$ which is the 1-skeleton of a
cellular decomposition, we can consider the space $\cC(\Gamma)$ of 
hyperideal circle patterns on $\Sigma$ with incidence graph $\Gamma$,
and the map $\Phi_\Gamma$ sending a hyperideal circle pattern to its
singular curvatures --- a function defined on the faces of $\Gamma$ ---
and its intersection angles --- a function defined on the edges of 
$\Gamma$. It follows from section 3 that the image of $\Phi_\Gamma$ is
contained in $\cD(\Gamma)$, the set of curvature and intersection data
allowed by Theorem \ref{tm:euclid}, while Lemma \ref{lm:rigidity-e}
shows that $\Phi_\Gamma$ is a local homeomorphism, and Lemma 
\ref{lm:compact} shows that $\Phi_\Gamma$ is proper. The arguments given
above in the hyperbolic case can also be used to show that each element
of $\cD(\Gamma)$ has at most one inverse image by $\Phi_\Gamma$. 

The only remaining point is that, for all $\Gamma$, any element of 
$\cD(\Gamma)$ has at least one inverse image by $\Phi_\Gamma$. This will
be proved using Theorem \ref{tm:hyperb} and the compactness stated in 
Lemma \ref{lm:compact}.

\paragraph{An approximation argument.}

Let $\kappa:\Gamma_2\rightarrow (-\infty, 2\pi)$ and $\theta:\Gamma_1
\rightarrow (0,2\pi)$ describe an element of $\cD(\Gamma)$, i.e. 
$\kappa$ and $\theta$ satisfy the hypothesis of Theorem \ref{tm:euclid}.

\begin{lemma} \label{lm:approx}
There exists sequences $(\kappa_n)_{n\in \N}$ and $(\theta_n)_{n\in \N}$,
with $\kappa_n:\Gamma_2\rightarrow (-\infty, 2\pi)$ and $\theta_n:
\Gamma_1\rightarrow (0,\pi)$, defined for $n\geq n_0$ for some $n_0\in \N$,
such that:
\begin{itemize}
\item for all $n\geq n_0$, $\kappa_n$ and $\theta_n$ satisfy the hypothesis
of Theorem \ref{tm:hyperb}, and therefore are the singular curvature and
intersection angles of a unique hyperideal circle pattern on $\Sigma$ with
a hyperbolic metric,
\item for all $f\in \Gamma_2$ and all $e\in \Gamma_1$, $\kappa_n(f)\rightarrow
\kappa(f)$ and $\theta_n(e)\rightarrow \theta(e)$ as $n\rightarrow \infty$.
\end{itemize}
\end{lemma}

\begin{proof}
Let $\kappa_n$ and $\theta_n$ be defined by:
$$ \forall f\in \Gamma_2, \kappa_n(f)=\kappa(f)+\frac{1}{n}~, ~~
\forall e\in \Gamma_1, \theta_n(e) = \theta(e) +\frac{1}{n}~. $$
We have supposed that $\kappa$ and $\theta$ satisfy the hypothesis of
Theorem \ref{tm:euclid}, in particular the sum of the $\kappa(f)$ over
all faces of $\Gamma$ is equal to $2\pi \chi(\Sigma)$. So, for all $n$,
the sum of $\kappa_n(f)$ over all faces of $\Gamma$ is strictly larger
than $2\pi\chi(\Sigma)$, so that condition (1) of Theorem \ref{tm:hyperb}
holds. 

Let $\Omega$ be an admissible domain in $(\Sigma,\Gamma)$. Still by our 
hypothesis, the second condition of Theorem \ref{tm:euclid} holds for 
$\Omega$:
$$ \sum_{e\in \Gamma_1, e\subset\dr\Omega} \theta(e) \geq (2\chi(\Omega) - 
m(\Omega))\pi - \sum_{f\in \Gamma_2, f\subset \Omega} \kappa(f)~, $$
with strict inequality except perhaps when $\Omega$ is a face of $\Gamma$.
If $\Omega$ is a face of $\Gamma$, then it has at least 3 edges, and
clearly this inequality is also valid for $\kappa_n$ and $\theta_n$,
for any $n\geq 1$. Otherwise, the inequality is strict, and it 
follows that, for $n$ large enough (depending on $\Omega$)
the inequality holds for $\kappa_n$ and $\theta_n$. Since there is a 
finite number of admissible domains in $(\Sigma, \Gamma)$, it follows that
for $n$ large enough the hypothesis of Theorem \ref{tm:hyperb} hold.
\end{proof}

\paragraph{Proof of Theorem 1.4.}

The proof of Theorem \ref{tm:euclid} now follows from a direct application
of Lemma \ref{lm:compact} to the hyperideal circle patterns associated
by Theorem \ref{tm:hyperb} to the functions $\kappa_n$ and $\theta_n$,
which are on hyperbolic surfaces which, after applying a sequence of 
homotheties so
that their area remains equal to $1$, converge to Euclidean surfaces.


\section{End of the proofs, remarks}
\label{se:end}

The first part of this section contains the proof of the main results
concerning circle patterns on surfaces with polygonal boundary, and the 
next contains a remark and some open
questions. 

\subsection{Circle patterns on surfaces with boundary}

We turn here to the proof of Theorem \ref{tm:framed-euclid} and Theorem
\ref{tm:framed-hyperb}.
The proofs of those two results follow the same pattern, so we give
the arguments mainly Theorem \ref{tm:framed-hyperb} only, leaving the (limited)
adapatation necessary for the Euclidean context to the interested reader. 
We consider now a hyperbolic surface $\Sigma$ with polygonal boundary,
along with a circle pattern $C$ -- as in the statement of Theorem 
\ref{tm:framed-hyperb} -- with incidence graph $\Gamma$ and extended
incidence graph $\Gamma'$. 

\paragraph{Doubling a circle pattern.}

We consider a closed surface with conical singularities, called $D(\Sigma)$,
which is the ``double'' of $\Sigma$; it is obtained by gluing two copies of
$\Sigma$ isometrically along their common boundary, the gluing being made
through the identity map. There is also a circle pattern in $D(\Sigma)$, 
which we call $D(C)$, which is the inverse image of $C$ under the canonical
projection from $D(\Sigma)$ to $\Sigma$, so that each circle of
$C$ lifts to the two circles in $D(\Sigma)$, 
corresponding to $C$ in the two copies of $\Sigma$. 

It is then a simple matter to check that the incidence graph of $D(C)$ is a
graph embedded in $D(\Sigma)$, which we call $D(\Gamma)$ here, which 
has: 
\begin{itemize}
\item two vertices for each vertex of $\Gamma$, i.e. 
  for each vertex of $\Gamma'$
  which is not on the boundary of the exterior face,
\item two faces for each face of $\Gamma$, and one for each face of $\Gamma'$
  which does not correspond to a face of $\Gamma$,
\item two edges for each edge of $\Gamma$, and one edge for each edge of
  $\Gamma'$ which is going from an  ``interior'' vertex of $\Gamma'$ to 
a ``boundary'' vertex of $\Gamma'$.
\end{itemize}
Moreover, the angle assigned to an edge of $D(\Gamma)$ is: 
\begin{itemize}
\item for edges corresponding to edges of $\Gamma$, 
the intersection angle between the two circles of $C$ corresponding 
to its endpoints,
\item for edges corresponding to edges of $\Gamma'$ going from an ``interior''
  to a ``boundary'' vertex, twice the angle on that edge (since the
  intersection angle between the circles corresponding to the endpoints is
  twice the angle between one of those circles and the corresponding boundary
  segment). 
\end{itemize}
Finally, the singular curvature assigned to a face of $D(\Gamma)$ is: 
\begin{itemize}
\item for faces corresponding to faces of $\Gamma$, identical to the singular
  curvature of that face, 
\item for faces corresponding to faces of $\Gamma'$ which are not faces of
  $\Gamma$, twice the angle assigned to the exterior edge of 
  that face in $\Gamma'$, since the
  singular angle of $D(\Sigma)$ at that point is twice the exterior angle of
  $\dr\Sigma$ at the corresponding point. 
\end{itemize}
This defines assignments $D(\kappa)$ and $D(\theta)$ of angles to the faces
and edges of $D(\Gamma)$, respectively. Note that their is a slight abuse
of notations here since $D(\kappa)$ depends not only on $\kappa$ but also
on $\theta$ (because the values of $D(\kappa)$ on faces corresponding to 
faces of $\Gamma'$ which are not faces of $\Gamma$ depend on the values of 
$\theta$ on the exterior edges of those faces).

The key point of the proofs is then:

\begin{lemma} \label{lm:equivalence}
$\Gamma$, $\kappa$ and $\theta$ satisfy the conditions of Theorem
\ref{tm:framed-hyperb} if and only if $D(\Gamma)$, $D(\kappa)$ and $D(\theta)$
satisfy the condition of Theorem \ref{tm:hyperb}.
\end{lemma}

\begin{proof}
First note that condition (1) of Theorem \ref{tm:hyperb} for the doubled 
surface is equivalent to condition (1) of Theorem \ref{tm:framed-hyperb}
for the original surface with boundary, because the curvature associated
to each face of $\Gamma$ counts once in each of the copies of $\Sigma$ which
are glued, the curvatures at the vertices of the gluing line are equal to
twice the exterior angles at the corresponding vertices of $\dr\Sigma$, 
while $\chi(D(\Sigma))=2\chi(\Sigma)$.

Suppose that $D(\Gamma)$, $D(\kappa)$ and $D(\theta)$
satisfy the condition of Theorem \ref{tm:hyperb}.
Let $\Omega$ be an admissible domain in $(\Sigma, \Gamma')$, and let
$D(\Omega)$ be the domain in $D(\Sigma)$ which is the inverse image of
$\Omega$ under the canonical projection from $D(\Sigma)$ to $\Sigma$.
We already know that equation (\ref{eq:cond-omega}) holds
for $D(\Omega)$, considered as an admissible domain in
$(D(\Sigma), D(\Gamma))$, so that:
$$ \sum_{\stackrel{e\in D(\Gamma)_1,}{e\subset \dr D(\Omega)}} \theta(e) = 
\pi(2\chi(D(\Omega)) - m(D(\Omega))) -
\sum_{\stackrel{f\in D(\Gamma)_2,}{f\subset D(\Omega)}} \kappa(f)~. $$
However the description of $D(\Gamma)$ -- and of the corresponding angles and
singular curvatures -- shows that: 
$$ \sum_{e\subset \dr D(\Omega)} \theta(e) = 
2 \sum_{e\in \Gamma_1} \theta(e) +
\sum_{\stackrel{e\in \Gamma'_1\setminus \Gamma_1,}
{e\subset \dr \Omega\setminus \dr
  \Sigma}} 2\theta(e)~, $$ 
where the first sum is over the edges of $\dr\Omega$ which are in 
$\Gamma$, while the second sum is over the edges of $\Gamma'$ which 
have one vertex on the exterior face of $\Gamma'$. The coefficient
$2$ in the first sum comes from the fact that each edge of $\Gamma$ 
corresponds to two edges in $D(\Gamma)$, while the same coefficient
in the second sum comes from the fact that the angle between the two
circles corresponding to the endpoints in $D(C)$ is twice the 
angle between the corresponding circle of $C$ and the adjacent 
segment of $\dr\Sigma$. In the same way: 
$$ \sum_{\stackrel{f\in D(\Gamma)_2,}{f\subset D(\Omega)}} \kappa(f) 
= 2\sum_{f\in \Gamma_2, f\subset
  \Omega} \kappa(f) + \sum_{\stackrel{e\in \Gamma'_1,}{e\subset \dr\Omega\cap 
\dr\Sigma}} 2\theta(e)~. $$ 
Finally, $\chi(D(\Omega)) = 2\chi(\Omega))-n(\Omega)$ 
and $m(D(\Omega))=2m(\Omega)$, and
it follows that (\ref{eq:cond-framed}) holds for $\Omega$.

Conversely, suppose that $\Gamma, \kappa$ and $\theta$ satisfy condition (2)
of Theorem \ref{tm:hyperb}, and let $\Omega$ be an admissible domain in 
$(D(\Sigma), D(\Gamma))$. Then $\Omega$ can be decomposed as the disjoint 
union of two domains $\Omega_+$ and $\Omega_-$, which are the intersections
of $\Omega$ with the two copies of $\Sigma$ which are glued in $D(\Sigma)$.
Then both $\Omega_+$ and $\Omega_-$ satisfy condition (2) of Theorem 
\ref{tm:framed-hyperb}, and the computation made above shows directly 
that $\Omega$ satisfies condition (2) of Theorem \ref{tm:hyperb}.
\end{proof}

\paragraph{Proof of Theorem 1.10.}

We can now show that Theorem 
\ref{tm:framed-euclid} (resp. Theorem \ref{tm:framed-hyperb})
is a consequence of Theorem \ref{tm:euclid} (resp. Theorem \ref{tm:hyperb}).

Consider first a framed hyperideal circle pattern $C$ on a compact surface
$\Sigma$ with boundary. Let $\Gamma$ be the incidence graph of $C$, and
let $\Gamma'$ be its extended incidence graph; let $\kappa$ and $\theta$
be the singular curvature and intersection data of $C$, so that 
$\kappa$ is a function defined on the set of faces of $\Gamma$ and 
$\theta$ is defined on the set of edges of $\Gamma'$. We consider again the
doubled surface $D(\Sigma)$, and the doubled circle pattern $D(C)$, which
has incidence graph $D(\Gamma)$, as 
above. By Theorem \ref{tm:hyperb}, $D(\Sigma)$ and $D(C)$ are uniquely
determined by $D(\Gamma)$ and by the intersection angles between the
circles and the singular curvature. This already shows that framed
hyperideal circle patterns on hyperbolic surfaces are uniquely determined
by their incidence graph and intersection angles and singular curvatures 
data. The same argument can be used for Euclidean circle patterns, based
on Theorem \ref{tm:euclid}.

To prove the existence of a circle pattern with given incidence
graph on a surface with boundary $\Sigma$ and with given intersection
angles and singular curvatures, consider the double $D(\Sigma)$ of
$\Sigma$, and the graph $D(\Gamma)$ embedded in $D(\Sigma)$, obtained
by the construction described above; $D(\Gamma)$ contains two copies of
$\Gamma'$ which are glued along their boundary faces (each boundary face
of one copy is identified with a boundary face of the other copy). 
There are natural intersection angles $D(\theta)$ and singular curvatures 
$D(\kappa)$ defined on $D(\Gamma)$ from $\kappa$ and $\theta$ (see above).
By Lemma \ref{lm:equivalence}, Theorem \ref{tm:hyperb} can be applied to
$(D(\Sigma), D(\Gamma))$ with the singular curvatures and intersection angles
data given by $D(\kappa)$ and $D(\theta)$.
Thus, the ``doubled'' data are 
realized uniquely as incidence graph, intersection angles and singular
curvature of a hyperideal circle pattern $D(\Sigma)$ on $D(\Sigma)$.

Moreover, this circle pattern (and the underlying metric) have the
same symmetry as $D(\Sigma)$ and $D(\Gamma)$, so $D(\Sigma)$ can be
``cut'' along a polygonal line to obtain two isometric surfaces with
polyhedral boundary, each with a framed hyperideal circle pattern, 
as required by Theorem \ref{tm:framed-hyperb}. The same arguments shows
Theorem \ref{tm:framed-euclid} from Theorem \ref{tm:euclid}.

\subsection{A remark and some open questions}

\paragraph{A alternate proof when $\kappa$ is non-negative.}

There is another possible proof of the main theorems, or at least of the
infinitesimal rigidity lemmas, which is based on a completely different
argument. However this proof only works under the hypothesis that the total
angle at the singular points is at most equal to $2\pi$, i.e. when the singular
curvature at the singuar points is positive. Moreover it only works for
strictly hyperideal circle patterns.

In those cases, the infinitesimal rigidity can be proved by first constructing
a hyperbolic cone-manifold with polyhedral boundary associated to a hyperideal
circle pattern, as in section 4, and then doubling it, to obtain a complete
finite volume cone-manifold, which has as singular locus a link (a disjoint
union of closed curves). The condition on the angles at the singular
points translates as the fact that the total angle around the singular curves
is less than $2\pi$, so that it is possible to apply a rigidity result of
Hodgson and Kerckhoff \cite{HK}, according to which those cone-manifolds can
not be deformed without changing the total angle around some of the singular
curves. 

It is quite striking that, in those examples, there are two ways to prove the
infinitesimal rigidity of those special circle patterns, one based on the
volume argument used here, and another based on the more analytic methods used
by Hodgson and Kerckhoff \cite{HK}. 

\paragraph{Spherical metrics.}

All the results presented here concern surfaces endowed with either Euclidean
or hyperbolic surfaces with conical singularities. However the basic result in
the theory, the Koebe theorem on circle packings, 
applies to the sphere. As already pointed out, the result of Bao and Bonahon
\cite{bao-bonahon} also translates as a statement on hyperideal circle
patterns on the sphere. Interesting recent results of Luo \cite{luo-sphere} 
also apply. It would be interesting
to know whether results similar to those proved here for Euclidean or
hyperbolic metrics holds for circle patterns on surfaces with
spherical metrics with conical singularities. 

The main problem there appears to be the infinitesimal rigidity of such circle
patterns, since the tools used in section 4 and 5 do not appear to work in the
spherical context. The
compactness results, however, could presumably be proved by the methods used
here.


\section*{Acknoledgements}

The results presented here owe much to many fruitful conversations with Boris
Springborn, who was strongly involved in parts of the research presented
here. I would also like to thank Igor Rivin for several conversations
concerning some techniques used in this text.

\bibliographystyle{alpha}
\def\cprime{$'$}

\end{document}